\documentclass{gtart_h}  

\def\ifplaintex{\expandafter\ifx\csname documentclass\endcsname\relax}

\ifplaintex 
\else
\fi

\expandafter\ifx\csname beginpicture\endcsname\relax
\expandafter\ifx\csname documentclass\endcsname\relax
\input pictex \else
\input prepictex \input pictex \input postpictex \fi\fi

\def\gt{{\mathsurround=0pt\it $\cal G\mskip-2mu$eometry \&\ 
$\cal T\!\!$opology}}        

\def\gtp{{\mathsurround=0pt\it $\cal G\mskip-2mu$eometry \&\ 
$\cal T\!\!$opology $\cal P\!$ublications}}  

\def\lognumber#1{\def\thelognumber{#1}}
\def\volumenumber#1{\def\thevolumenumber{#1}}
\def\papernumber#1{\def\thepapernumber{#1}}
\def\volumeyear#1{\def\thevolumeyear{#1}}

\def\pagenumbers#1#2{\def\startpage{#1}\def\finishpage{#2}}
\def\published#1{\def\publishdate{#1}}
\def\proposed#1{\def\theproposer{#1}}
\def\seconded#1{\def\theseconders{#1}}
\def\received#1{\def\receiveddate{#1}}
\def\revised#1{\def\reviseddate{#1}}
\def\accepted#1{\def\accepteddate{#1}}

\def\coverauthors#1{\def\thecoverauthors{#1}}

\def\asciiaddress#1{\def\theasciiaddress{#1}}
\def\asciiemail#1{\def\theasciiemail{#1}}

\long\def\asciiabstract#1{\long\def\theasciiabstract{#1}}

\let\\\par\let\thelognumber\relax
\let\thevolumenumber\relax\let\thepapernumber\relax
\let\thevolumeyear\relax\let\thesamplenumber\relax\let\startpage\relax
\let\finishpage\relax\let\publishdate\relax\let\receiveddate\relax
\let\reviseddate\relax\let\accepteddate\relax\let\theasciititle\relax
\let\theasciiauthors\relax\let\theasciiaddress\relax
\let\theasciiabstract\relax
\let\theasciiemail\relax\let\theshortauthors\relax\let\theshorttitle\relax
\let\thecoverauthors\relax

\long\def\maketitlep{   

\count0=\startpage

\gt\hfill      
\beginpicture
\setcoordinatesystem units <0.33truein, 0.33truein> point at 2.2 0.9
\setplotsymbol ({$\cal G$})
\plotsymbolspacing=9truept
\circulararc 315 degrees from 0 1 center at 0 0
\setplotsymbol ({$\cal T$})
\circulararc 315 degrees from 1 -1 center at 1 0
\endpicture
\break
{\small\ifx\thesamplenumber\relax 
Volume \else Sample
\fi\thevolumenumber\ (\thevolumeyear)
\startpage--\finishpage\nl
Published: \publishdate}
\vglue 0.5truein plus 0.4fil minus 0.1truein

{\parskip=0pt\leftskip 0pt plus 1fil\def\\{\par\smallskip}{\ifplaintex\large
\else\Large\fi\bf\thetitle}\par\medskip}   

\vglue 0pt plus 0.1fil 

{\parskip=0pt\leftskip 0pt plus 1fil\def\\{\par}{\sc\theauthors}
\par\medskip}

\vglue 0pt plus 0.1fil 

{\small\parskip=0pt\let\newline\\
{\leftskip 0pt plus 1fil\def\\{\par}{\sl\theaddress}\par}
\expandafter\ifx\theemail\relax    
\relax\else\vglue 5pt plus 0.02fil minus 2pt\def\\{\stdspace{\rm 
and}\stdspace} 
\cl{Email:\stdspace\tt\theemail}\fi
\ifx\theurl\relax                  
\relax\else\vglue 5pt plus 0.02fil minus 2pt\def\\{\stdspace{\rm 
and}\stdspace}
\cl{URL:\stdspace\tt\theurl}\fi\par}

\vglue 7pt plus 0.3fil minus 3pt

{\bf Abstract}
\vglue 5pt plus 0.1fil minus 2pt

\theabstract

\vglue 7pt plus 0.3fil minus 3pt

{\bf AMS Classification numbers}\quad Primary:\quad \theprimaryclass

Secondary:\quad \thesecondaryclass

\vglue 5pt plus 0.3fil minus 2pt

{\bf Keywords:}\quad \thekeywords

\vglue 10pt plus 0.5fil minus 5pt

{\small  Proposed: \theproposer\hfill Received: \receiveddate\nl
Seconded: \theseconders\hfill 
\ifx\reviseddate\relax                         
Accepted: \accepteddate                        
\else
Revised: \reviseddate                          
\fi}
\eject
}       

\let\maketitlepage\maketitlep
\let\maketitle\maketitlepage

\font\phead=cmsl9 scaled 950
\font=cmsl9 scaled 1050
\font\pnum=cmbx10 scaled 913
\font\lnum=cmbx10 
\font\pfoot=cmsl9 scaled 950
\font=cmsl9 scaled 1050
\ifplaintex
\headline{\vbox to 0pt{\vskip -4.5mm\line{\small\phead\ifnum
\count0=\startpage ISSN 1364-0380 (on line)
1465-3060 (printed) \hfill {\pnum\folio}\else\ifodd\count0\def\\{ }%
\ifx\theshorttitle\relax\thetitle\else\theshorttitle\fi\hfill{\pnum\folio}
\else\def\\{ and }{\pnum\folio}\hfill\ifx\theshortauthors\relax\theauthors
\else\theshortauthors\fi\fi\fi}\vss}}
\footline{\vbox to 0pt{\vglue 0mm\line{\small\pfoot\ifnum\count0=\startpage
\copyright\ \gtp\hfill\else
\gt, Volume \thevolumenumber\ (\thevolumeyear)\hfill\fi}\vss
}}
\else
\makeatletter
\def\@oddhead{{\small\lhead\ifnum\count0=\startpage ISSN 1364-0380 (on line)
1465-3060 (printed) \hfill {\lnum\number\count0}\else\ifodd\count0
\def\\{ }\ifx\theshorttitle\relax \thetitle \else\theshorttitle\fi\hfill
{\lnum\number\count0}\else\def\\{ and }{\lnum\number\count0}
\hfill\ifx\theshortauthors\relax 
\theauthors\else\theshortauthors\fi\fi\fi}}\def\@evenhead{\@oddhead}
\def\@oddfoot{\small\lfoot\ifnum\count0=\startpage\copyright\ \gtp\hfill\else
\gt, Volume \thevolumenumber\ (\thevolumeyear)\hfill\fi}
\def\@evenfoot{\@oddfoot}
\makeatother
\fi

\newwrite\gtoutfile
\long\gdef\makeheadfile{  
{\def\\{, }\def\s{ }
\immediate\openout\gtoutfile head.xxx
\immediate\write\gtoutfile{Proxy-for: \ifx\theasciiauthors\relax
\theauthors\else\theasciiauthors\fi\s<\ifx\theasciiemail\relax\theemail\else\theasciiemail\fi>}
\immediate\write\gtoutfile{\noexpand\\}
\immediate\write\gtoutfile{Authors: \ifx\theasciiauthors\relax
\theauthors\else\theasciiauthors\fi}
{\def\\{ }\immediate\write\gtoutfile{Title: \ifx\theasciititle\relax
\thetitle\else\theasciititle\fi}}
\immediate\write\gtoutfile{Subj-class: GT or SG or MG etc}
\immediate\write\gtoutfile{MSC-class: \theprimaryclass\ifx\thesecondaryclass\relax\else, \thesecondaryclass\fi}
\immediate\write\gtoutfile{Journal-ref: Geom. Topol. \thevolumenumber
(\thevolumeyear) \startpage-\finishpage}
\immediate\write\gtoutfile{Comments: Published by Geometry and Topology at}
\immediate\write\gtoutfile{\s\s http://www.maths.warwick.ac.uk/gt/GTVol\thevolumenumber/paper\thepapernumber.abs.html}
\immediate\write\gtoutfile{\noexpand\\}
\immediate\write\gtoutfile{}
\ifx\theasciiabstract\relax
\immediate\write\gtoutfile{\theabstract}\else
\immediate\write\gtoutfile{\theasciiabstract}\fi
\immediate\write\gtoutfile{}
\immediate\write\gtoutfile{\noexpand\\}
\immediate\write\gtoutfile{}
\immediate\closeout\gtoutfile}}  

\def\maketitlepage{\maketitlep\makeheadfile}
\let\maketitle\maketitlepage

\lognumber{338}

\received{10 July 2003}

\volumenumber{8}\papernumber{31}\volumeyear{2004}
\pagenumbers{1127}{1188}
\published{8 September 2004}
\revised{4 February 2004}
\accepted{1 September 2004}
\proposed{Joan Birman}
\seconded{David Gabai, Yasha Eliashberg}

\usepackage{enumerate,amssymb,amscd,amsmath,graphicx,psfrag}

\newcommand{\barrow}{\includegraphics[height=.10in]{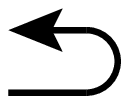}}

\newcommand{\sph}{S^2}
\newcommand{\tPhi}{\tilde{\Phi}}
\newcommand{\MIA}{\mathcal{M}\mathcal{I}\mathcal{A}}
\newcommand{\bs}{{\bar s}}
\newcommand{\pR}{{\partial\cR}}
\newcommand{\I}{^{-1}}
\renewcommand{\k}{\kappa}
\newcommand{\im}{$f\co I\barrow$\ }
\newcommand{\opti} 
{\raisebox{-.04cm}{$\stackrel{{\scriptscriptstyle 0}}{{\scriptscriptstyle
1}}$}}
\newcommand{\ol}[1]{\overline{#1}}

\newcommand{\pichere}[2]
 {\begin{center}\includegraphics[width=#1\textwidth]{#2}\end{center}}

\newcommand{\stretchpichere}[3]
{\begin{center}\includegraphics[width=#1\textwidth,height=#2]{#3}\end{center}}
 \newcommand{\lab}[3]{\psfrag{#1}[#3]{$\scriptstyle{#2}$}}

\newcommand{\cF}{{\mathcal{F}}}
\newcommand{\cR}{{\mathcal{R}}}
\newcommand{\cG}{{\mathcal{G}}}

\newcommand{\R}{\mathbb{R}}
\newcommand{\Z}{\mathbb{Z}}
\newcommand{\Q}{\mathbb{Q}}
\newcommand{\N}{\mathbb{N}}
\newcommand{\ofi}{\mathbb{I}}
\newcommand{\V}{\mathbb{V}}
\newcommand{\ofe}{\mathbb{E}}

\newcommand{\bt}{\operatorname{bt}}
\newcommand{\IE}{\operatorname{IE}}
\newcommand{\RE}{\operatorname{RE}}
\newcommand{\BP}{\operatorname{BP}}
\newcommand{\diam}{\operatorname{diam}}
\newcommand{\Int}{\operatorname{Int}}
\newcommand{\lhe}{\operatorname{lhe}}
\newcommand{\rhe}{\operatorname{rhe}}
\newcommand{\NBT}{\operatorname{NBT}}
\newcommand{\KS}{\operatorname{KS}}
\newcommand{\Width}{\operatorname{Width}}
\newcommand{\Mod}{\operatorname{Mod}}
\newcommand{\Area}{\operatorname{Area}}

\newtheorem{thm}{Theorem}
\newtheorem{lemma}[thm]{Lemma}

\theoremstyle{definition}

\newtheorem{defn}{Definition}
\newtheorem{defns}[defn]{Definitions}
\newtheorem{example}{Example}
\newtheorem{examples}[example]{Examples}
\newtheorem{alg}{Algorithm}

\theoremstyle{remark}

\newtheorem{remark}{Remark}

\begin{document}

\title{Unimodal generalized pseudo-Anosov maps}
\authors{Andr\'e de Carvalho\\Toby Hall}
\coverauthors{Andr\noexpand\'e de Carvalho\\Toby Hall}
\address{Departamento de Matem\'atica Aplicada, IME - USP\\Rua do
Mat\~ao 1010, Cidade Universit\'aria\\05508-090 S\~ao Paulo, SP,
Brazil\\{\rm and}\\Department of Mathematical Sciences, University of
Liverpool\\Liverpool L69 7ZL, UK} 
\asciiaddress{Departamento de Matematica Aplicada, IME-USP\\Rua do
Matao 1010, Cidade Universitaria\\05508-090 Sao Paulo, SP,
Brazi\\and\\Department of Mathematical Sciences, University of
Liverpool\\Liverpool L69 7ZL, UK}

\asciiemail{andre@ime.usp.br, T.Hall@liv.ac.uk}
\gtemail{\mailto{andre@ime.usp.br}\qua{\rm and}\qua\mailto{T.Hall@liv.ac.uk}}

\begin{abstract}
An infinite family of generalized pseudo-Anosov homeomorphisms of the
sphere $S$ is constructed, and their invariant foliations and singular
orbits are described explicitly by means of generalized train
tracks. The complex strucure induced by the invariant foliations is
described, and is shown to make~$S$ into a complex sphere. The
generalized pseudo-Anosovs thus become quasiconformal automorphisms of
the Riemann sphere, providing a complexification of the unimodal
family which differs from that of the Fatou/Julia theory.
\end{abstract}
\asciiabstract{%
An infinite family of generalized pseudo-Anosov homeomorphisms of the
sphere S  is constructed, and their invariant foliations and singular
orbits are described explicitly by means of generalized train
tracks. The complex strucure induced by the invariant foliations is
described, and is shown to make S into a complex sphere. The
generalized pseudo-Anosovs thus become quasiconformal automorphisms of
the Riemann sphere, providing a complexification of the unimodal
family which differs from that of the Fatou/Julia theory.}


\maketitlepage

\section{Introduction}
\subsection{Overview}

Pseudo-Anosov homeomorphisms were introduced by Thurston
in his classification of surface homeomorphisms up to isotopy. A
surface homeomorphism $\Phi\co S\barrow$ is pseudo-Anosov if it
preserves a transverse pair of measured foliations with finitely many
singularities, expanding one foliation uniformly by a
factor~$\lambda>1$ and contracting the other by a
factor~$1/\lambda$. Pseudo-Anosov maps can be described
combinatorially by train tracks --- a class of graphs embedded in the
surface with additional information on the vertices which specifies the
turns a train riding along the track can make at each vertex --- and
their endomorphisms. Using information obtained from an associated
transition matrix, the surface can be reconstructed by identifying sides
of Euclidean rectangles foliated by horizontal and vertical line
segments. The pseudo-Anosov expands the rectangles horizontally and
contracts them vertically, mapping them as dictated by the train track
map. Up to this point, the discussion is finite: train tracks are
graphs with finitely many edges (Thurston's theorem can be proved
algorithmically, the main step being to find a finite train track
invariant under a given isotopy class of homeomorphisms) and train
track maps are finite-to-one. It is possible, however, to forgo some
of the finiteness requirements --- the graphs remain finite, the
endomorphisms remain finite-to-one, but the additional information at
the vertices is allowed to be infinite --- but otherwise to go through
the construction of the maps as before. This leads to the construction
of {\em generalized pseudo-Anosov homeomorphisms}.  These are defined
similarly to pseudo-Anosov homeomorphisms, except that the invariant
foliations are permitted to have infinitely many singularities,
provided that they accumulate on only finitely many points.  The
purpose of this paper is to give a detailed construction and
description of an infinite family of generalized pseudo-Anosovs of the
sphere for which the underlying graph and graph map are the simplest
possible: an interval and a unimodal endomorphism (ie, a continuous
piecewise monotone map of the interval with exactly two monotone
pieces).

In~\cite{Ha}, a complete description of the family of pseudo-Anosov
maps with underlying unimodal interval endomorphisms was given. It was
shown that there is a countable family of such maps parameterized by a
rational number between 0 and 1/2, called {\em height}. Height turns
out to be a braid type invariant and this leads to the proof of weak
universality results for families of plane homeomorphisms passing from
trivial to chaotic dynamics as parameters are varied. Height also
plays a central role in this paper and, in turn, the results presented
here provide a geometric interpretation of it. The family of unimodal
generalized pseudo-Anosovs extends that of unimodal pseudo-Anosovs.
The height specifies the behaviour of the maps at infinity: given a
rational~$m/n$, there is an interval of kneading sequences of
height~$m/n$, whose associated generalized pseudo-Anosovs have the
same behaviour at infinity. The generalized pseudo-Anosov is a
pseudo-Anosov for exactly one kneading sequence in this interval.

This paper provides an explicit description of the
generalized train track associated to any periodic or preperiodic
kneading sequence, depending crucially on its height. The process of
constructing a generalized pseudo-Anosov from a generalized train
track map is similar to that of constructing a pseudo-Anosov from a
train track map, but requires more care because of the more intricate
nature of the identifications carried out on the sphere. The
topological tool used to guarantee that the identification space is
again a sphere is Moore's theorem about monotone upper semi-continuous
decompositions of the sphere.

The invariant foliations of a generalized pseudo-Anosov define a
complex structure on the sphere away from the accumulations of
singularities. It is shown that for the unimodal generalized
pseudo-Anosovs considered here, these accumulations are removable
singularities of the complex structure, so that the sphere is a
complex sphere, with the foliations being the horizontal and vertical
trajectories of an integrable quadratic differential having infinitely
many zeros and poles. The construction therefore provides a
complexification of unimodal maps as quasiconformal automorphisms of
the Riemann sphere, in contrast to the complexification arising via
the theory of Fatou/Julia, where one thinks of the unimodal map as the
real slice of an endomorphism of the Riemann sphere.

By a suitable normalization, the sphere of definition of the
generalized pseudo-Anosovs can be identified canonically with the
Riemann sphere, and hence the family of unimodal generalized
pseudo-Anosovs, which initially are constructed on abstract
topological spheres, can be regarded as a family of Teichm\"uller
mappings of the Riemann sphere, making it possible to consider taking
limits within the family. This is a necessary step in the problem of
constructing a completion of the set of all pseudo-Anosov
homeomorphisms of the sphere.

Section~\ref{sec:tim} describes the class of {\em thick interval
maps}. These are homeomorphisms of the sphere which provide the
starting point for the definition and construction of invariant
generalized train tracks, a process which is described in
Section~\ref{sec:gentt}. Section~\ref{sec:symbolics} provides a
summary of necessary results on unimodal maps and Smale's horseshoe,
and defines the subclass of {\em unimodal} thick interval maps which
are used in the remainder of the paper. In Section~\ref{sec:outside},
the {\em outside dynamics} of a unimodal map is defined and analysed:
intuitively, this provides a description of those orbits of a unimodal
map which are never lost `inside' the fold. The main contents of the
paper can be found in Sections~\ref{sec:unitt},~\ref{sec:genpA}
and~\ref{sec:comp}.  The invariant generalized train track for a given
unimodal thick interval map is described explicitly, a detailed
account is given of how the generalized train track map can be used to
construct the corresponding generalized pseudo-Anosov, and the complex
structure induced by the invariant foliations is analysed.

\smallskip\noindent {\bf Acknowledgements}\qua The authors are grateful
for the referee's careful reading of the paper and helpful
comments. This material is based upon work supported by the National
Science Foundation under Grant No.\ 0203975. The first author is
supported by FAPESP Grant No.\ 02/05072-5.

\subsection{Definitions and notation}

Let $X$ be a metric space. An {\em isotopy} on $X$ is a continuous map
$\psi\co X\times[0,1]$\break $\to X$ with the property that the slice maps
$\psi_t\co X\to X$ are homeomorphisms for $0\le t\le1$. A {\em
pseudo-isotopy} is defined similarly but it is only required that the
slice maps be homeomorphisms for $0\le t<1$. In particular, this means
that the map $\psi_1$ is a {\em near-homeomorphism}, ie, it can be
approximated arbitrarily closely by homeomorphisms. An isotopy or
pseudo-isotopy $\psi\co X\times[0,1]\to X$ is said to be {\em
supported} on a subset $U$ of $X$ if all of its slice maps agree
on~$X\setminus U$.

Symbolic dynamics on both $\{0,1\}^\N$ and $\{0,1\}^\Z$ will be used
in this paper. Elements of these spaces are regarded as semi-infinite
or bi-infinite sequences of $0$s and $1$s, and in the case of
$\{0,1\}^\Z$ a period is placed before the origin of the sequence
(ie, the image of $0\in\Z$): $\ldots s_{-2}s_{-1}\cdot
s_0s_1s_2\ldots\in\{0,1\}^\Z$. If $w\in\{0,1\}^n$ for some~$n$,
then the notation $w^\infty$ is used to indicate semi-infinite
repetition of~$w$ ($w^\infty=wwww\ldots\in\{0,1\}^\N$), while
$\overline{w}$ is used to indicate bi-infinite repetition
($\overline{w}=\ldots www\cdot www\ldots\in\{0,1\}^\Z$).

The notation $S^1$ and $S^2$ is reserved for the standard 1- and
2-dimensional spheres, and different symbols are used to denote
general topological spheres.

\begin{figure}[ht!]
\begin{center}
\pichere{0.6}{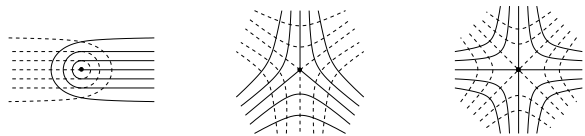}
\caption{Pronged singularities of the invariant foliations}
\label{fig:prongs}
\end{center}
\end{figure}

A homeomorphism $\Phi\co S\barrow$ of a smooth surface~$S$ is called a
{\em generalized pseudo-Anosov map} if there exist
\begin{enumerate}[a)]
\item a finite $\Phi$-invariant set $\Sigma$; 
\item a pair $(\cF^s,\mu^s)$, $(\cF^u,\mu^u)$ of transverse measured
foliations of $S\setminus\Sigma$ with countably many pronged
singularities (with local charts as depicted in
Figure~\ref{fig:prongs}), which accumulate on each point of $\Sigma$
and have no other accumulation points. The transverse measures are
required to be equivalent to Lebesgue measure on transversals;
\item a real number $\lambda>1$;
\end{enumerate}
such that
\[\begin{array}{ccc}
\Phi(\cF^s,\mu^s) &=& (\cF^s,\frac{1}{\lambda}\mu^s) \\
\Phi(\cF^u,\mu^u) &=& (\cF^u,\lambda\mu^u). 
\end{array}\]

In particular, a generalized pseudo-Anosov is a pseudo-Anosov map if
and only if there are only finitely many pronged singularities (ie,
$\Sigma=\emptyset$). Note that the measures necessarily have full
support and no atoms.

\section{Markov thick interval maps}
\label{sec:tim}
\subsection{Thick interval maps}
\label{subsec:tim}

Thick graph maps are a class of surface homeomorphisms which have been
described and used in several papers (for
example~\cite{BH,dCH1,FM,dC}). In this section a brief description of
thick interval maps (where the surface is the sphere and the graph is
an interval) is given. As the name suggests, a thick interval map is
essentially an interval endomorphism which has been thickened up and
made into a homeomorphism of the sphere, whose dynamics reflects that
of the underlying interval map.

Throughout the paper, $\sph$ will be thought of as the one-point
compactification of $\R^2$. The point at infinity will be denoted 
$\infty$ and the thick intervals defined below will always be assumed
not to contain $\infty$.

\begin{figure}[ht!]\small
\psfrag{1}[t]{$1$-junction}
\psfrag{2}[t]{$2$-junction}
\psfrag{l}[b]{leaves}
\begin{center}
\pichere{0.6}{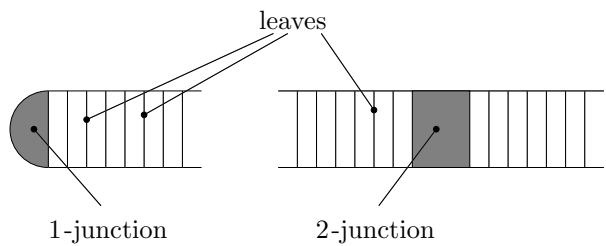}
\caption{Charts in a thick interval}
\label{fig:chart}
\end{center}
\end{figure}

\begin{defns}
A {\em thick interval} $\ofi\subset\sph$ is a closed topological $2$-disk 
partitioned into compact {\em decomposition elements}, such that
\begin{enumerate}[i)]
\item each decomposition element of $\ofi$ is either a {\em leaf}
homeomorphic to $[0,1]$, or a {\em junction} homeomorphic to a closed
$2$-disk;
\item the boundary in $\ofi$ of each junction consists of one or two
disjoint arcs: if there is one (respectively two) such arc(s) the
junction is called a 1-junction (respectively 2-junction); 
\item there are exactly two 1-junctions and finitely many 2-junctions;
\item each decomposition element is contained in a chart as depicted
in Figure~\ref{fig:chart}.
\end{enumerate}

If $\ofi$ is a thick interval, then the space obtained by collapsing
each decomposition element to a point is an interval, whose vertices
(the two endpoints, coming from the two 1-junctions, and finitely many
valence~$2$ vertices) correspond to the junctions of $\ofi$.  The
union of the junctions of $\ofi$ is denoted $\V$, and the
components of $\ofi\setminus\V$ are called {\em strips}: each strip
is therefore homeomorphic to $(0,1)\times[0,1]$.  The union of the
closures of the strips is denoted $\ofe$. Thus $\ofe\cap\V$ is
a union of closed arcs which are the boundary components (in $\ofi$)
of both the junctions and the strips of $\ofi$.
\end{defns}

The following definition of a thick interval map is more restrictive
than that used in other papers: conditions~iv) to~vi) have been added
in order to make the definition and construction of generalized train
tracks more straightforward. However, thick interval maps in the sense
of~\cite{dCH1} can be made to conform to the definition below by
simple isotopies and changes in the thick interval structure which do
not change the dynamics. In order to make the description more
explicit, fixed orientation-preserving coordinate maps
$h_s\co\overline{s}\to[0,1]\times[0,1]$ are introduced on the
closure of each strip $s$, with the property that the leaves of $s$
are of the form $h_s^{-1}(\{x\}\times[0,1])$ for $0<x<1$. 

\begin{defn}
\label{defn:tim}
A {\em thick interval map} is an orientation-preserving
homeomorphism $F\co(\sph,\ofi)\barrow$ such that:
\begin{enumerate}[i)]
\item $F(\ofi)\subset\Int(\ofi)$.
\item If $\gamma$ is a leaf of $\ofi$, then
$F(\gamma)$ is contained in a decomposition element, and
$\diam(F^n(\gamma))\to 0$ as $n\to\infty$. If~$J$ is a junction
of~$\ofi$, then~$F(J)$ is contained in a junction.
\item The point at infinity is a repelling fixed point 
whose basin contains $\sph\setminus\ofi$.
\item $F$ is linear with respect to the coordinates $h_s$: in each
connected component of $s_i \cap F^{-1}(s_j)$, where $s_i$, $s_j$ are
strips, $F$ contracts vertical coordinates uniformly by a factor
$\mu_{j}<1$ and expands horizontal coordinates uniformly by a factor
${\lambda}_{j} \ge 1$;
\item If $J, J'$ are junctions such that $F(J) \subset J'$,
then $F(\partial_\ofi J) \subset \partial_\ofi J'$;
\item If $J$ is a junction with $F^n(J)\subseteq J$ for some (least)
$n\ge1$, then $J$ has an attracting periodic point of least period $n$
in its interior whose basin contains $\Int(J)$.
\end{enumerate}
\end{defn}

\begin{remark}
Item iii) in the definition says that the dynamics of a thick
interval map in $\sph\setminus\ofi$ is easily understood and
uninteresting. 
\end{remark}

Let $F\co(\sph,\ofi)\barrow$ be a thick interval map, $\tilde{I}$ be
the interval obtained by collapsing each decomposition element of
$\ofi$ to a point, and $\tilde{\pi}\co\ofi\to\tilde{I}$ be the
canonical projection. Then~$F$ induces a continuous map $\tilde{f}\co
\tilde{I}\barrow$: however~$\tilde{f}$ is locally constant at
preimages of~$\tilde{\pi}(\V)$. It is therefore convenient to
collapse all intervals on which some iterate~${\tilde{f}}^n$ is
constant. Because closed intervals are collapsed to points, the
quotient of the interval $\tilde I$ is either a point or an
interval~$I$. In the latter case~$\tilde f$ induces a continuous
piecewise strictly monotone interval map~$f\co I\barrow$, which is
called the {\em quotient} of $F\co(\sph,\ofi)\barrow$. All thick
interval maps considered in this paper will have an interval map for
quotient (see Section~\ref{subsec:markov}). Note that
if~$\pi\co\ofi\to I$ denotes the canonical projection, then $\pi\circ
F|_\ofi=f\circ\pi$.

\begin{example}
\label{ex:hs}
The first example is Smale's horseshoe map which will be denoted 
$F_0\co(\sph,\ofi)\barrow$ here and in what follows. It is shown in
Figure~\ref{fig:gpA_hs}. The horseshoe has three fixed points, denoted
$x$, $x_0$, and $x_1$. The fixed point~$x$ is attracting
(condition~vi)), while $x_0$ and~$x_1$ are saddles by conditions~iv)
and~vi).
The quotient interval map --- the tent map
$f_0\co I\barrow $ --- is also shown in the figure.
\begin{figure}[ht!]
\begin{center}
\lab{tp}{\tilde\pi}{}
\lab{tf}{\tilde{f}_0}{}
\lab{p}{\pi}{}
\lab{T}{\ofi}{l}
\lab{F}{F_0}{}
\lab{f}{f_0}{}
\lab{x}{x}{}
\lab{y}{x_0}{}
\lab{z}{x_1}{}
\pichere{0.4}{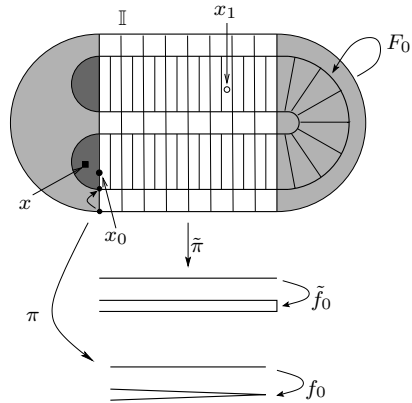}
\caption{The horseshoe map}
\label{fig:gpA_hs}
\end{center}
\end{figure} 
\end{example}

\begin{example}
Figure~\ref{fig:1001011thtrmap} depicts an example of a thick interval
map associated to a horseshoe periodic orbit: such thick interval maps
are the starting point for the construction of unimodal generalized
train tracks described in Section~\ref{sec:unitt}. The interval
endomorphism~$f$ is one with kneading sequence
$\k(f)=(1001011)^\infty$ (see Section~\ref{sec:symbolics}).
\begin{figure}[ht!]
\begin{center}
\lab{F}{F}{bl}
\lab{f}{f}{l}
\pichere{0.8}{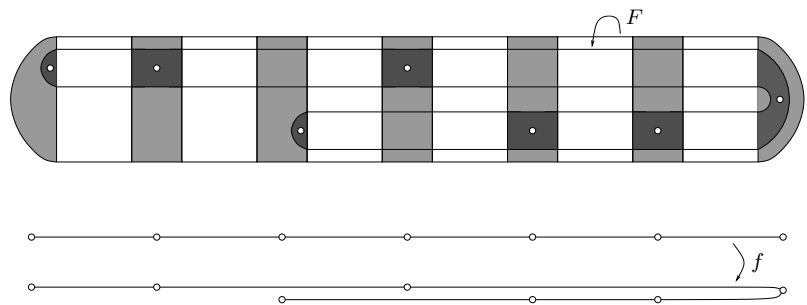}
\caption{A thick interval map associated to a horseshoe periodic
orbit}
\label{fig:1001011thtrmap}
\end{center}
\end{figure}
\end{example}

\subsection{The MIA property}
\label{subsec:markov}

First the basic concepts of the Perron-Frobenius theory for
non-negative integer matrices are described (see~\cite{Gant}).  Let
$M$ be a square matrix with non-negative integer entries, which is
not equal to the $1\times 1$ matrix $(1)$. $M$ is said to be {\em
reducible} if, by a permutation of the index set, it is possible to
put it in triangular block form:
\[ M = \left( \begin{array}{cc}
A & 0 \\ B & C \end{array} \right) , \] with $A$ and $C$
non-trivial square matrices.  Otherwise, $M$ is said to be {\em
irreducible}. An irreducible matrix~$M$ has a unique positive
eigenvector (up to scaling), and the associated eigenvalue $\lambda$,
called the {\em Perron-Frobenius eigenvalue} of~$M$, is simple and is
equal to the spectral radius of~$M$.  $M$ is {\em irreducible and
aperiodic} if there is a positive integer~$k$ such that every entry of
$M^k$ is positive. In this case the Perron-Frobenius eigenvalue
$\lambda$ satisfies $\lambda>1$, and is the only eigenvalue of $M$ on
the circle $\{z\in \mathbb{C}\,:\, |z|=\lambda\}$.

\medskip 
Let $F\co(S^2,\ofi)\barrow$ be a thick interval map. Associate a {\em
transition matrix} $M=\left(m_{ij}\right)$ to $F$ in the following
way: let $s_1,s_2,\ldots,s_n$ be the strips of~$\ofi$ and set
\[ m_{ij}=\mbox{number of times\ } F(s_j) \mbox{\ crosses\ }
s_i.\] 

\begin{defn}
A thick interval map is {\em MIA} if its associated transition matrix
is irreducible and aperiodic (the `M' stands for Markov).
\end{defn}

If~$F$ is MIA then the quotient interval~$I$ is not a point, and
indeed the projection~$\pi(s_j)$ of each strip of~$\ofi$ is a
non-trivial interval. For there are integers~$n_j$ such that
$F^{n_j}(s_j)$ crosses~$s_j$ at least twice, with horizontal expansion
and vertical contraction bounded away from~$1$. Thus the set of points
of~$s_j$ whose forward~$F^{n_j}$-orbits remain in~$s_j$ contains
$h_{s_j}^{-1}(C\times[0,1])$ for some Cantor set $C\subseteq[0,1]$,
and $\pi(h_{s_j}^{-1}(x_1,y_1))\not=\pi(h_{s_j}^{-1}(x_2,y_2))$ whenever
$x_1,x_2$ are distinct points of~$C$ and $y_1,y_2\in[0,1]$.

\section{Generalized train tracks}
\label{sec:gentt}

This section contains a description of generalized train tracks
associated to thick interval maps. It follows~\cite{dC}, where a
description of generalized train tracks associated to general 
thick graph maps is given.

Let $\ofi$ be a thick interval and $A\subset\ofi$ be a finite {\em
puncture} set, each of whose points lies in the interior of a junction
and such that no junction contains more than one point of $A$. For
each strip $s$ of $\ofi$, let $\gamma_s$ be (the image of) the arc
$t\mapsto h_s^{-1}(t,1/2)$, which joins the two boundary components of
$s$ in $\ofi$ and intersects each leaf of $s$ exactly once. Let $\RE$
denote the set of these arcs. The endpoints of the arcs $\gamma_s$ are
called {\em switches} and the set of switches is denoted~$L$.

\begin{defns}
A {\em generalized train track} $\tau\subseteq\ofi\setminus A$ is a
graph with vertex set~$L$ and countably many edges, each of which
intersects $\partial\V$ only at~$L$, such that
\begin{enumerate}[i)]
\item The edges of $\tau$ which intersect the interior of $\ofe$ are
precisely the elements of $\RE$, and
\item No two edges $e_1$, $e_2$ contained in a given junction $J$ are
{\em parallel}: that is, they do not bound a disk which contains no
point of $A$ or other edges.
\end{enumerate}

Two generalized train tracks $\tau$ and $\tau'$ are {\em equivalent},
denoted $\tau\sim\tau'$, if they are isotopic by an isotopy supported
on $\V\setminus A$.

The edges of $\tau$ which are contained in $\ofe$ (that is, the
elements of $\RE$) are called {\em real}, and the others (contained in
$\V$) are called {\em infinitesimal}.  Let $\IE$ denote the set of
infinitesimal edges of $\tau$.

A generalized train track~$\tau$ is {\em finite} if it has only
finitely many edges. An infinitesimal edge is called a {\em loop} if
its two endpoints coincide, and a {\em bubble} if in addition it
bounds an open disk which is disjoint from~$\tau$. A bubble of $\tau$
is {\em homotopically trivial} if this disk contains no point of $A$,
and is {\em homotopically non-trivial} otherwise.

Clearly $\tau$ is determined by its infinitesimal edges, and may be
written $\tau=\tau(\IE)$ when the thick interval $\ofi$ and the set
$A$ are clear from the context.
\end{defns}

Note that $\tau$ is not required to be connected: while the
generalized train tracks constructed below will always be connected,
disconnected ones are needed during the construction.  

\begin{defns}
A homotopy $\{\alpha_t\}$ of a path $\alpha\co[0,1]\to X$ is said to
be {\em relative to} $U\subset X$ if, for each $s\in[0,1]$,
$\alpha_0(s)\in U$ implies $\alpha_t(s)\in U$ for all $t$.
Let~$[\alpha]$ be a homotopy class of paths in $\sph\setminus A$
relative to $\partial_\ofi\V$, with endpoints
in~$\partial_\ofi\V$. Then $[\alpha]$ is {\em carried} by a
generalized train track $\tau$ if it can be realized by an edge-path
in $\tau$ with alternating real and infinitesimal edges.
\end{defns}

Generalized train tracks are normally drawn in such a way as to
suggest that if a real edge $\gamma_s$ and an infinitesimal edge $e$
share a common endpoint, then their union is a smooth (branched if the
endpoints of $e$ coincide) $1$-manifold. With this intuition,
$[\alpha]$ is carried by $\tau$ if it can be realized by a smooth path
in $\tau$ (or, even more intuitively, by a train running along
$\tau$). It is more convenient, however, to express this smoothness
combinatorially as above.

The next aim is to define the image of $\tau$ under a thick interval
map $F$. In order that this image should itself be a generalized train
track, it is necessary to apply pseudo-isotopies to $F(\tau)$ so as to
squash $F(\tau)\cap\ofe$ onto the real edges, and amalgamate
pairs of parallel edges.

Let $(\ofi,A)$ be a thick interval together with a puncture set. 
On each strip $s$ of~$\ofi$ define $\psi_s\co \bs
\times [0,1] \to \bs$ using the coordinates~$h_s$ by \[ \psi_s(x,y,t) =
(x, (1-t)y+t/2)\] (thus $\psi_s(\cdot,1)$ maps $\bs$ onto
$\gamma_s$). Extend these maps to a pseudo-isotopy $\psi_0\co
\sph\times[0,1]\to \sph$ in the following way: first extend the
$\psi_s$ to mutually disjoint disk neighbourhoods $U_s\supset\bs$, with
$U_s\subset \sph\setminus A$, so that they are isotopies on
$U_s\setminus\bs$, the identity on $\partial U_s$, and send $\V$
into $\V$; then extend them to be the identity elsewhere.

If $\tau$ is a graph in $\ofi\setminus A$ which satisfies the
definition of a generalized train track except that there are finitely
many pairs $(e,e')$ of parallel edges, define a pseudo-isotopy
$\psi_{e,e'}$ supported on a neighbourhood in $\V$ of the closed
disk~$\Delta$ bounded by $e\cup e'$ with the property that
$\psi_{e,e'}(e,1)=\psi_{e,e'}(e',1)$ but no points outside~$\Delta$
are identified by $\psi_{e,e'}(\cdot,1)$. Since this pseudo-isotopy
introduces no new pairs of parallel edges, it is possible to define a
pseudo-isotopy $\psi_1\co\sph\times[0,1]\to\sph$ (whose dependence on
$\tau$ is suppressed) by composing successive pseudo-isotopies on each
pair of parallel edges, with the property that $\psi_1(\tau,1)$ is a
generalized train track.

In this paper, the puncture set $A$ will always be taken to be the set
of attracting periodic points of the thick interval map~$F$ (so in
particular $F(A)=A$). Conditions~iv) and~vi) in the definition of a
thick interval map ensure that the points of~$A$ are contained in the
interiors of the junctions of $\ofi$, with at most one point of~$A$ in
each junction.

\begin{defns}
Let $F\co(\sph,\ofi,A)\barrow$ be a thick interval map, where $A$
is the set of attracting periodic points of $F$, and let
$\tau\subseteq\ofi\setminus A$ be a generalized train track. Since $F$
restricts to an embedding $(\V,A)\barrow$ and the underlying
interval endomorphism is piecewise monotone, there can be only
finitely many pairs of parallel edges in $\psi_0(F(\tau),1)$, and the
{\em image} $F_*(\tau)$ of $\tau$ under $F$ is defined as
$F_*(\tau)=\psi(\tau)$, where the near-homeomorphism $\psi\co
\sph\barrow$ is given by $\psi(x)=\psi_1(\psi_0(F(x),1),1)$.

A generalized train track $\tau$ is {\em $F$-invariant} if
$F_*(\tau)$ is equivalent to~$\tau$.

If $\tau$ is $F$-invariant, then by definition there is a
homeomorphism $H\co\sph\barrow$, isotopic to the identity by an
isotopy supported on $\V\setminus A$, such that
\mbox{$H(\psi(\tau))=\tau$}.  The {\em train track map} $\phi\co\tau\barrow$
associated to $F\co(\sph,\ofi,A)\barrow$ is the restriction of
$H\circ\psi$ to $\tau$, well-defined up to homotopy relative to the
vertices of $\tau$.
\end{defns}

\medskip\medskip

The following straightforward result guarantees the existence of
invariant generalized train tracks, and its proof provides a method
for constructing them.

\begin{thm}
\label{thm:ttexist}
Let $F\co(\sph,\ofi,A)\barrow$ be a thick interval map, where $A$ is
the set of attracting periodic points of $F$. Then there exists an
$F$-invariant generalized train track $\tau\subseteq\ofi\setminus A$.
\end{thm}

\begin{proof}
Let $\tau_0=\tau(\emptyset)$ be the generalized train track with no
infinitesimal edges. By definition,
$\tau_0\subseteq\tau_1=F_*(\tau_0)$. Let $\tau_2'=F_*(\tau_1)$. Then
$\tau_1$ need not be a subset of~$\tau_2'$, since~$\tau_2'$ is only
defined up to equivalence, but there is a generalized train track
$\tau_2\sim\tau_2'$ with $\tau_1\subseteq\tau_2$. Continuing in this
way, construct a nested sequence $(\tau_n)$ of generalized train
tracks with $F_*(\tau_n)\sim\tau_{n+1}$ for all~$n$. Then
$\tau=\bigcup_{n\ge 0}\tau_n$ is a generalized train track which
satisfies $F_*(\tau)\sim\tau$.
\end{proof}

The $F$-invariant generalized train track~$\tau$ constructed in this
proof is minimal, in the sense that it contains only those edges which
arise as images of its real edges: more precisely, if~$\tau'$ is also
$F$-invariant then there is a subset~$\IE$ of the infinitesimal edges
of $\tau'$ such that $\tau\sim\tau(\IE)$. For this reason~$\tau$ is
referred to as {\em the $F$-invariant generalized train track}.

\begin{example}
\label{ex:hstrtr}
The invariant generalized train track for the horseshoe map $F_0$ is
shown in Figure~\ref{fig:hstrtr}. The set $A$ consists of the fixed
point $x$ contained in the left 1-junction of $\ofi$. No bubble
encloses it and it is not shown in the figure. It is instructive to
construct this train track starting from $\tau(\emptyset)$ as described in
the proof of Theorem~\ref{thm:ttexist}.
\begin{figure}
\lab{t}{\tau}{}
\lab{f}{\phi}{b}
\pichere{0.4}{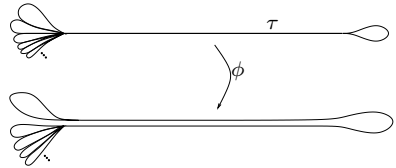}
\caption{The invariant train track for the horseshoe}
\label{fig:hstrtr}
\end{figure}
\end{example}

The invariant train track $\tau$ and the train track map $\phi$ can be
thought of as more careful 1-dimensional representations of the thick
interval $\ofi$ and the thick interval map $F\co
(S^2,\ofi,A)\barrow$\,\,. Whereas \im does not pay attention to junctions
--- they are collapsed to points --- the map $\phi\co\tau\barrow$
gives a careful account of the behaviour of the images of strips under
iterates of $F$ inside the junctions.

\begin{figure}[ht!]
\begin{center}
\pichere{0.6}{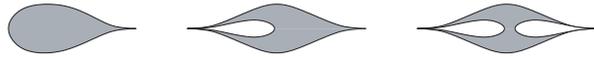}
\caption{Examples of $n$-gons for $n=1,3,4$}
\label{fig:ngons}
\end{center}
\end{figure}

\begin{defn}
An {\em infinitesimal polygon} of $\tau$ is a component of
$\sph\setminus\tau$ bounded by finitely many infinitesimal edges. It is
called an {\em $n$-gon} if it is bounded by $n$ infinitesimal edges
(see Figure~\ref{fig:ngons}, in which the $n$-gons are shaded).
\end{defn}

\begin{remark}
Bigons (2-gons) can only occur if they contain a point of $A$
(otherwise the two edges bounding the bigon would be parallel).
\end{remark}

\section{Unimodal maps, symbolic dynamics and the horseshoe}
\label{sec:symbolics}

This section contains a summary of the theory of unimodal maps,
symbolic dynamics, and the Smale horseshoe which will be used later.
Detailed accounts can be found in~\cite{Devaney,MT,dMvS,KH} (unimodal
maps) and~\cite{Devaney,Ha} (horseshoe). 

\subsection{Unimodal maps}
\label{subsec:unimodal}
Let $I=[a,b]$. For the purposes of this paper, a {\em unimodal map} on
$I$ is a continuous surjection $f\co I\barrow$ with $f(b)=a$, for
which there exists $c\in(a,b)$ such that $f$ is strictly increasing on
$[a,c]$ and strictly decreasing on $[c,b]$ (and in
particular~$f(c)=b$). The point $c$ is the {\em critical} or {\em
turning point} of~$f$, and $b=f(c)$ is its {\em critical value}.

The conditions that $f$ be surjective and that $f(b)=a$ are not
standard. However, a unimodal map~$f$ which doesn't satisfy them has
trivial dynamics outside of its {\em dynamical interval}
$[f^2(c),f(c)]$, and hence there is no loss of generality, and a gain
of convenience, in adding these requirements. Note also that some
authors don't require the monotonicity on~$[a,c]$ and~$[c,b]$ to be
strict, and would refer to a map defined as above as {\em strictly
unimodal}.

Symbolic dynamics for unimodal maps is introduced 
by defining the {\em itinerary} of a point $x\in I$ to be the
sequence
\[i(x)=s_0s_1s_2\ldots\,\in\{0,1\}^\N\] 
given by
\begin{equation*}
s_j=
\begin{cases}
0& \text{if $f^j(x)<c$}\\
1& \text{if $f^j(x)\ge c$.}
\end{cases}
\end{equation*}
This defines a map $i\co I\to\{0,1\}^\N$ with the property that
$i\circ f=\sigma\circ i$, where $\sigma\co\{0,1\}^\N\barrow$ is
the {\em shift map} given by
$\sigma(s_0s_1s_2\ldots)=s_1s_2s_3\ldots$.

\begin{remark}
In other contexts it is useful to define the itinerary $i(x)$ to lie
in $\{0,C,1\}^\N$, where $i_j(x)=C$ if $f^j(x)=c$. This distinction
will not be necessary here.
\end{remark}

\smallskip
The {\em unimodal order} $\preceq$ is a total order (but not a
well-ordering) on $\{0,1\}^\N$ which reflects the usual ordering of
points in~$I$. If~$s$ and~$t$ are distinct elements of~$\{0,1\}^\N$,
then $s\prec t$ if and only if $\sum_{j=0}^n s_j$ is even, where~$n$
is least such that $s_n\not=t_n$. For general $s,t\in\{0,1\}^\N$,
define $s\preceq t$ if and only if $s\prec t$ or $s=t$.

The unimodal order is defined precisely in order that $i(x)\prec
i(y)\implies x<y$ for all $x,y\in I$.  It follows immediately that
$x<y\implies i(x)\preceq i(y)$: the possibility of equality in this
partial converse cannot be excluded in general.  A unimodal map $f$ is
said to {\em have no homtervals} if the full converse holds,
ie, $x<y\iff i(x)\prec i(y)$ (a {\em homterval} is a non-trivial
interval $J$ with $c\not\in f^n(J)$ for all $n\ge0$). 

The itinerary of the critical value of a unimodal map $f\co
I\barrow$ plays a particularly important role: the {\em kneading
sequence} $\k(f)$ of a unimodal map $f\co I\barrow$ is the
itinerary $\k(f)=i(b)$. Since $\sigma(\k(f))=i(a)$, it follows that
$\sigma(\k(f))\preceq i(x)\preceq\k(f)$ for all $x\in I$, and in
particular $\sigma(\k(f))\preceq\sigma^n(\k(f))\preceq\k(f)$ for all
$n\ge0$. This statement characterizes those elements of $\{0,1\}^\N$
which are the kneading sequence of some unimodal map, as expressed by
the following definition and theorem (which is a translation into the
language of this paper of Theorem~12.1 of~\cite{MT}).

\begin{defn}
An element $s$ of $\{0,1\}^\N$ is {\em a kneading sequence} if
$\sigma(s)\preceq \sigma^n(s)$\break $\preceq s$ for all $n\in\N$.
\end{defn}

\begin{thm}
An element $s$ of $\{0,1\}^\N$ is equal to $\k(f)$ for some unimodal
map~$f$ if and only if it is a kneading sequence.
\end{thm}

\begin{defn}
An MIA thick interval map $F\co(S^2,\ofi)\barrow$ is {\em unimodal} 
if its quotient interval endomorphism \im is unimodal. In this case,
define the {\em kneading sequence} $\kappa(F)$ of $F$ to be
$\kappa(f)$.
\end{defn}

The quotient $f\co I\barrow$ of a unimodal thick interval map has no
homtervals. For since $F\co(S^2,\ofi)\barrow$ is MIA, there is
some~$n$ such that $f^n$ expands each subinterval~$J$ of~$I$ for which
$c\not\in\bigcup_{i=0}^{n-1} f^i(J)$, the expansion being bounded away
from~$1$ with respect to the projections of the coordinate maps~$h_s$
on the strips of~$\ofi$.

Because thick intervals as defined here have only finitely many
junctions, the orbits of the vertices of $I$ under $f$ are finite. In
particular, the kneading sequence of a unimodal thick interval map is
always periodic or preperiodic.  The set of periodic or preperiodic
kneading sequences whose associated transition matrices are
irreducible and aperiodic will be denoted $\MIA$. Notice that given
a kneading sequence in $\MIA$, there is a natural construction of a
unimodal MIA thick interval map with that kneading sequence:
contruct a piecewise affine interval map with the given kneading
sequence and thicken it.

\subsection{Symbolic dynamics for the horseshoe}

This section contains a brief summary of the application of symbolic
dynamics to the horseshoe map of Example~\ref{ex:hs}:
see~\cite{Devaney} for a more detailed description.  The set
$\Lambda=\bigcap_{n\in\Z}F_0^n(\overline{s})$ (where~$s$ is the strip
of~$\ofi$) is a Cantor set, and an itinerary homeomorphism
$i\co\Lambda\to\{0,1\}^\Z$ is defined by setting
\begin{equation*}
i_j(z)=
\begin{cases}
0& \text{if $F_0^j(z)\in V_0$}\\
1& \text{if $F_0^j(z)\in V_1$}
\end{cases}
\end{equation*}
where $V_0,V_1$ are the left and right connected components of
$F_0\I(s)\cap s$ respectively.  The itinerary homeomorphism conjugates
$F_0|_\Lambda$ with the shift map $\sigma\co\{0,1\}^{\Z}\barrow$,
defined by
\[\sigma(\ldots s_{-2}s_{-1}\cdot s_0s_1\ldots)= 
\ldots s_{-1}s_{0}\cdot s_1s_2 \ldots\]
Let $f_0\co I\barrow$ be the quotient interval map (which is conjugate
to a full tent map). Then the invariant Cantor set $\Lambda$ of
$F_0$ inside $s$ projects to the whole interval $I$, and this
projection establishes a 1-1 correspondence between the $F_0$-periodic
orbits in $\Lambda$ and the $f_0$-periodic orbits in~$I$.

Note that for all $z\in\Lambda$, $i(\pi(z))$ is obtained from $i(z)$
by deleting all symbols before the origin. The correspondence between
periodic orbits is reflected in the correspondence between itineraries
in the obvious way: the itinerary of a period $n$ periodic point of
$F_0$ is of the form $\overline{w}$ for some length~$n$ word~$w$, and
the itinerary of its projection is $w^\infty$. This correspondence
will be invoked without further comment in the remainder of the paper.

A periodic orbit $P$ of $F_0$ of (least) period $n$ is described by
its {\em code} $c_P\in\{0,1\}^n$, which is given by the first $n$
symbols of the itinerary $i(p)$ of its rightmost point~$p$: thus, for
example, the period~5 orbit which contains the point with itinerary
$\overline{01001}$ has code $10010$. A word $w\in\{0,1\}^n$ is
therefore the code of a period $n$ horseshoe orbit if and only if it
is maximal in the sense of the following definition.
\begin{defn}
$w\in\{0,1\}^n$ is {\em maximal} if
$\sigma^i(w^\infty)\prec w^\infty$ for $1\le i<n$.
\end{defn}
It follows that if $w$ is the code of a periodic orbit of~$F_0$, then 
$w^\infty$ is a kneading sequence.

\subsection{Braid types} 
Let $F\co S^2\barrow$ be a homeomorphism (which in this paper will
always be orient\-ation-preserving), and $A$ be a finite
$F$-invariant subset of~$S^2$. Then the {\em braid type} $\bt(A,F)$ is
defined to be the isotopy class of~$F$ relative to~$A$ up to
topological change of coordinates~\cite{Boy}: that is,
$\bt(A,F)=\bt(B,G)$ if and only if there is an orientation-preserving
homeomorphism $h\co S^2\barrow$ with $h(A)=B$ such that $h\circ F\circ
h^{-1}$ is isotopic to~$G$ relative to~$B$. Using Thurston's
classification theorem for surface homeomorphisms~\cite{Thurston},
braid types can be classified as finite order, reducible, or
pseudo-Anosov.

The braid type of a periodic orbit~$P$ of the horseshoe map $F_0\co
S^2\barrow$ is defined to be $\bt(P\cup\{\infty\},F_0)$. Thus
horseshoe periodic orbits can also be classified as finite order,
reducible, or pseudo-Anosov.

\subsection{Height}
The description of unimodal generalized train tracks in
Section~\ref{sec:unitt} depends upon the notion of the {\em height} of
an element~$s$ of $\MIA$. The height is
defined using words~$c_q\in\{0,1\}^{n+1}$ associated to each rational
$q=m/n\in(0,1/2]$. Motivation for the definition can be found
in~\cite{Ha}.

\begin{defn}
Given $q=m/n\in\Q\cap(0,1/2]$, define a word $c_q\in\{0,1\}^{n+1}$
as follows. Let $L_q$ be the straight line in $\R^2$ from $(0,0)$
to $(n,m)$. For $0\le i\le n$, let $s_i=1$ if $L_q$ crosses some line
$y=\mbox{integer}$ for $x\in(i-1,i+1)$, and $s_i=0$ otherwise. Then
$c_q=s_0s_1\ldots s_n$.
\end{defn}

\begin{figure}[ht!]
\lab{0}{0}{}
\lab{1}{1}{}
\lab{2}{(0,0)}{r}
\lab{3}{(10,3)}{l}
\begin{center}
\pichere{0.4}{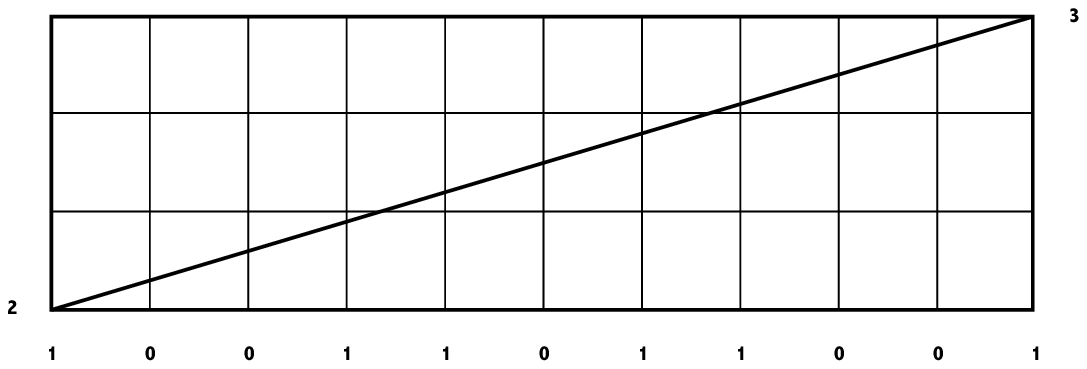}
\end{center}
\caption{$c_{3/10}=10011011001$}
\label{fig:linebox2}
\end{figure}

\begin{example}
Figure~\ref{fig:linebox2} shows that $c_{3/10}=10011011001$. 
\end{example}

The words $c_q$ are manifestly palindromic.  Their general form is
indicated by the examples in Table~\ref{tab}, in which the column
headings and row headings denote the numerator and denominator of~$q$
respectively. The $n-2m+1$ zeros are partitioned `as even-handedly as
possible' into $m$ subwords (possibly empty), separated by $11$.

\begin{table}[ht!]
\begin{center}
\begin{tabular}{|c||c|c|c|c|}
\hline
& $1$ & $2$ & $3$ & $4$ \\
\hline\hline
$3$ & $1001$ & & & \\ \hline
$4$ & $10001$ & & & \\ \hline
$5$ & $100001$ & $101101$ && \\ \hline
$6$ & $1000001$ & & & \\ \hline
$7$ & $10000001$ & $10011001$ & $10111101$ & \\ \hline
$8$ & $100000001$ & & $101101101$ & \\ \hline
$9$ & $1000000001$ & $1000110001$ & & $1011111101$ \\ \hline
$10$ & $10000000001$ & & $10011011001$ & \\ \hline
$11$ & $100000000001$ & $100001100001$ & $100110011001$ &
$101101101101$ \\ \hline
\end{tabular}
\caption{Examples of the words $c_q$}
\label{tab}
\end{center}
\end{table}

The following description of the words $c_q$ is easily
shown to be equivalent: given $q=m/n$, define integers $\kappa_i(q)$
for $1\le i\le m$ by
\begin{equation}
\label{eq:kappa}
\kappa_i(q)=\left\{
\begin{array}{ll}
\lfloor 1/q\rfloor -1 & \mbox{ if }i=1\\
\lfloor i/q\rfloor - \lfloor (i-1)/q\rfloor - 2 & \mbox{ if }2\le
i\le m
\end{array}
\right.
\end{equation}
($\lfloor x\rfloor$ denotes the greatest integer which does not
exceed $x$).
Then \[c_q=10^{\kappa_1(q)}110^{\kappa_2(q)}11\ldots
110^{\kappa_m(q)}1.\]

\begin{example}
Let $q=3/7$, so $\kappa_1(q)=\lfloor 7/3\rfloor-1=1$,
$\kappa_2(q)=\lfloor 14/3\rfloor-\lfloor 7/3\rfloor-2=4-2-2=0$, and
$\kappa_3(q)=\lfloor 21/3\rfloor-\lfloor 14/3\rfloor-2=7-4-2=1$. Thus
$c_{3/7}=10111101$.
\end{example}

 The next lemma~\cite{Ha} motivates the definition of height: in
particular, it will imply that the height function
$q\co\{0,1\}^\N\to[0,1/2]$ is decreasing with respect to the
unimodal order $\prec$ on $\{0,1\}^\N$ and the usual order on
$[0,1/2]$.

\begin{lemma}
\label{lem:order}
For each $q\in\Q\cap(0,1/2]$, the word $c_q1$ is
maximal. Moreover, if $q,r\in\Q\cap(0,1/2]$ with $q<r$ then
$(c_r1)^\infty\prec(c_q1)^\infty$.
\end{lemma}

\begin{defns}
Let $c\in\{0,1\}^\N$. Then the {\em height} $q(c)\in[0,1/2]$ of $c$ is
given by
\[q(c)=\inf\{q\in\Q\cap(0,1/2]\,:\,q=1/2\text{ or
}(c_q1)^\infty\prec c\}.\] 
If $P$ is a horseshoe periodic orbit of period $n\ge 2$ with code
$c_P$, then the {\em height} $q(P)$ of $P$ is given by
$q(P)=q(c_P^\infty)$.
\end{defns}

The height of a horseshoe periodic orbit is a braid type
invariant. The height of any element of $\MIA$ is rational, and can be
computed using an algorithm described in~\cite{Ha}.

The next result describes the kneading sequences with given rational
height~\cite{Ha}.

\begin{defn}
For each $m/n\in(0,1/2)$, define $w_{m/n}\in\{0,1\}^{n-1}$ to be the
word obtained by deleting the last two symbols of $c_{m/n}$, and
$\widehat{w}_{m/n}\in\{0,1\}^{n-1}$ to be the reverse of $w_{m/n}$.
\end{defn}

\begin{thm}
\label{thm:heightinterval}
Let $s\in\MIA$. Then $q(s)\in[0,1/2)\cap\Q$, and
\begin{enumerate}[\rm a)]
\item $q(s)=0$ if and only if $s=10^\infty$.
\item If $m/n\in(0,1/2)$, then $q(s)=m/n$ if and only if
\[(w_{m/n}1)^\infty\preceq s\preceq
w_{m/n}01(1\widehat{w}_{m/n})^\infty.\]
\end{enumerate}
\end{thm}

Note that $w_{m/n}01(1\widehat{w}_{m/n})^\infty$ can be written more
concisely as $c_{m/n}(1\widehat{w}_{m/n})^\infty$, but the former
expression is more suggestive in calculations, as in the proof of
Lemma~\ref{lem:kiq} below. In fact, although it isn't immediately
apparent from the definitions, this kneading sequence is preperiodic
to $(w_{m/n}1)^\infty$, ie, there is an integer~$j$ with
$\sigma^j((1\widehat{w}_{m/n})^\infty)=(w_{m/n}1)^\infty$.

The endpoints of the intervals of kneading sequences of given height
will be important in the remainder of the paper, as will the kneading
sequences $(c_{m/n}1)^\infty$ used to define the height. The acronym
NBT in the following definitions stands for `no bogus transitions',
and reflects the original motivation of these ideas.

\begin{defns}
Let $m/n\in(0,1/2)$. Then write $\NBT(m/n)=(c_{m/n}1)^\infty$,
$\lhe(m/n)=(w_{m/n}1)^\infty$,
$\rhe(m/n)=c_{m/n}(1\widehat{w}_{m/n})^\infty$, and $\KS(m/n)$
for the set of kneading sequences $s\in\MIA$ with $\lhe(m/n)\preceq
s\preceq\rhe(m/n)$ (ie, the set of kneading sequences of height $m/n$).
\end{defns}

\begin{example}
Let $m/n=2/7$, so $c_{2/7}=10011001$, $w_{2/7}=100110$, and
$\widehat{w}_{2/7}=011001$. Then $\NBT(2/7)=(100110011)^\infty$,
$\lhe(2/7)=(1001101)^\infty$, and
$\rhe(2/7)=10011001(1011001)^\infty=10(0110011)^\infty$. An element
$s$ of $\MIA$ lies in $\KS(2/7)$ (ie, has height $2/7$) if and only
if
\[(1001101)^\infty\preceq s\preceq 10(0110011)^\infty.\]
\end{example}

Notice that by Theorem~\ref{thm:heightinterval},
\[\MIA=\{10^\infty\}\cup\bigcup_{m/n\in(0,1/2)}\KS(m/n).\]

\begin{lemma}
\label{lem:kiq}
Let $s\in\MIA$ have height $q=m/n>0$. Then either
$s=(w_q1)^\infty$ or $s$ has $c_q$ as an initial word. In particular,
if $s$ is periodic then either \mbox{$s=\lhe(m/n)$} (period~$n$), or
$s=\NBT(m/n)$ (period~$n+2$), or the period of $s$ is at least $n+3$.
\end{lemma}
\begin{proof}
By~Theorem~\ref{thm:heightinterval}, $s$ has $w_q$ as an initial
word. If the next symbol is $0$, then since $w_q$ is odd and $s\preceq
w_q01(1\widehat{w}_q)^\infty$, it follows that $s=w_q01\ldots$,
ie, $s=c_q\ldots$. If the next symbol is $1$, then either
$s=(w_q1)^\infty$, or let $k\ge1$ be the greatest integer such that
$s=(w_q1)^kd$ for some $d\in\{0,1\}^\N$. Then $d\preceq s$ (since
$s$ is a kneading sequence), but $d\succeq(w_q1)^\infty$ (by
Theorem~\ref{thm:heightinterval}), and hence $d=w_q1\ldots$,
contradicting the definition of $k$. The proof of the final statement
follows readily, and can be found in~\cite{Ha} (Theorem~3.5).
\end{proof}

\section{The outside dynamics of a unimodal map}
\label{sec:outside}
\begin{figure}
\lab{a}{a}{}\lab{b}{b}{}\lab{c}{c}{}\lab{p}{p}{}\lab{f}{f}{}\lab{fa}{f(b)}{}
\lab{fb}{f(a)}{}\lab{fc}{f(c)}{}\lab{fp}{f(p)}{}
\lab{Ga}{\gamma}{}\lab{ha}{\hat a}{}\lab{hb}{\hat b}{}\lab{pu}{p_u}{}
\lab{Tf}{\theta}{}\lab{Tha}{\theta(\hat a)}{}\lab{Thb}{\theta(\hat b)}{}
\lab{TGa}{\theta(\gamma)}{}\lab{Tpu}{\theta(p_u)}{}
\pichere{.9}{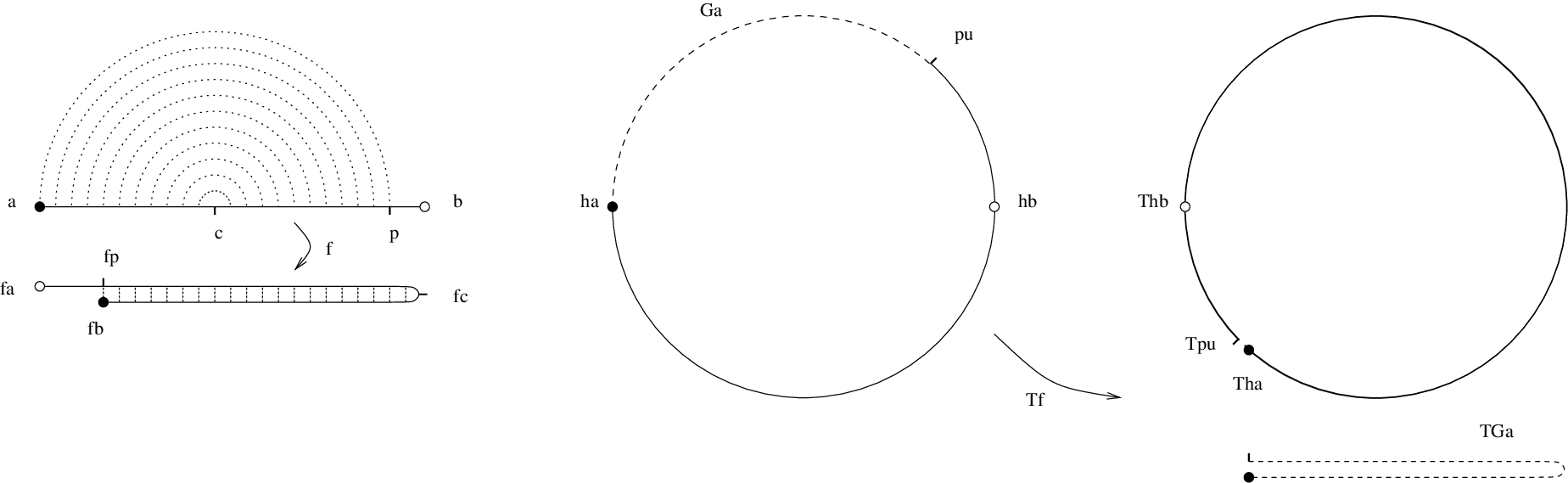}
\caption{A unimodal map and the induced outside map}
\label{fig:outside}
\end{figure}

The unimodal maps considered in this paper are destined to be
thickened into thick interval maps, and as such have an implicit
two-dimensional structure: the interval $[a,b]$ is thought of as being
folded at the critical point $c$, and laid down over itself in such a
way that points to the right of $c$ end up above points to the
left. Thus all of the points below the interval and some of the points
above it (namely those whose image is to the left of $f(a)$) remain
outside the interval, whereas the other points above the interval are
trapped in the fold (see the left hand side of
Figure~\ref{fig:outside} for clarification). The aim of this section
is to formalize this intuitive idea, and to analyze the dynamics of
points whose entire orbits remain outside the interval. This will play
an important role in the construction of generalized pseudo-Anosovs,
when it is necessary to `sew up the outside boundary'.

The outside of the interval is represented by a circle $S^1$ obtained
by gluing together two copies of $[a,b]$ at their endpoints, and the
unimodal map $f\co [a,b]\barrow$ induces a map $\theta\co S^1\to
S^1\cup[a,b]$, reflecting the two possible fates of the image of a
point on the outside: to remain outside, or to be folded inside the
interval.

Let $f\co I=[a,b]\barrow$ be a unimodal map with critical point $c$.
Let $S^1$ be the unit circle in~$\R^2$, coordinatize both the upper
and lower halves of $S^1$ with coordinates in $I$ in such a way that
$(-1,0)$ has coordinate $a$ and $(1,0)$ has coordinate $b$, and let
$\pi\co S^1\to I$ take each point of $S^1$ to the point of $I$ given
by its coordinate. Thus $\pi^{-1}(a)$ and $\pi^{-1}(b)$ each contain
a single point, denoted $\hat a$ and $\hat b$ respectively, while for
an interior point $x$ of $I$, $\pi^{-1}(x)$ contains two points,
denoted $x_u$ (in the upper half circle) and $x_l$ (in the lower half
circle).

Let $p\in I$ be the point with $p>c$ and $f(p)=f(a)$. The
(discontinuous) function $\theta\co S^1\to S^1\cup I$ is defined by
\begin{eqnarray*}
\theta(\hat{a})&=&f(a)_l\\
\theta(\hat{b})&=&\hat{a}\\
\theta(x_u)&=&f(x)\quad\text{ if $x<p$}\\
\theta(x_u)&=&f(x)_l\quad\text{if $x\ge p$}\\
\theta(x_l)&=&f(x)_l\quad\text{if $x<c$}\\
\theta(c_l) &=&\hat{b}\\
\theta(x_l)&=&f(x)_u\quad\text{if $x>c$},
\end{eqnarray*}
reflecting the intuitive notion of the action of $f$ on the outside of
the interval (see Figure~\ref{fig:outside}), and in particular
satisfying $\pi\circ\theta=f\circ\pi$. Note that the set of
points mapped by $\theta$ into $I$ is an open interval
$\gamma=\{x_u\,:\,a<x<p\}$, and that the complement of $\gamma$ is
mapped strictly monotonically onto $S^1$, with both endpoints mapping
to $f(a)_l$.

The following result describes the relationship between the dynamics
on the outside and the height of $\k(f)$ which will be used later. If
$x,y\in S^1$ then the notation $(x,y)$ always denotes the open
interval with endpoints $x$ and $y$ which is disjoint from $\gamma$
(so $x,y\not\in\gamma$), while $(x,y)_\gamma$ denotes the interval
contained in $\gamma$ (so $x,y\in\overline{\gamma}$). Similar notation
is used for half-open and closed intervals: intervals which intersect
but are not contained in $\gamma$ are not used.

\begin{thm}
\label{thm:outside}
Let $f\co [a,b]\barrow$ be the quotient of a unimodal thick interval
map with $s=\k(f)\in\MIA$, and let $\theta\co S^1\to S^1\cup I$ be the
induced map on the outside.  Let $q(s)=m/n>0$. Then
$\theta^i(\hat{a})\not\in\gamma$ for $0\le i<n$, while
$\theta^n(\hat{a})\in\overline{\gamma}$, with
\begin{enumerate}[\rm i)]
\item $\theta^n(\hat a)=\hat a \iff s=\lhe(m/n)$
\item $\theta^n(\hat a)=c_u \iff s=\NBT(m/n)$
\item $\theta^n(\hat a)=p_u \iff s=\rhe(m/n)$
\item $\theta^n(\hat a)\in (\hat a,c_u)_\gamma \iff
\lhe(m/n)\prec s\prec\NBT(m/n)$
\item $\theta^n(\hat a)\in (c_u,p_u)_\gamma \iff
\NBT(m/n)\prec s\prec\rhe(m/n)$. 
\end{enumerate}
Let $\Lambda=\bigcap_{n=0}^\infty \theta^{-n}(S^1)$ be the set of
points whose $\theta$-orbits remain in $S^1$. Then
$\Lambda$ contains exactly one periodic orbit, which has
period $n$ and whose points are permuted as by a rigid rotation of the
circle through angle $2\pi m/n$. If $\lhe(m/n)\prec s\prec \rhe(m/n)$
then there are no other points in $\Lambda$, while if $s=\lhe(m/n)$
(respectively $s=\rhe(m/n)$) then $\Lambda$ is the union of this
periodic orbit and the set $\{\theta^{-i}(p_u)\,:\,i\ge0\}$
(respectively $\{\theta^{-i}(\hat{a})\,:\,i\ge0\}$).
\end{thm}

Notice in particular that the periodic orbit in $\Lambda$ is the orbit
of $\hat{a}$ in the case when $s=\lhe(m/n)$, and is the orbit of $p_u$ when
$s=\rhe(m/n)$ (since $\theta(p_u)=\theta(\hat{a})$). 

Three lemmas are used in the proof. The integers $\k_i(q)$ in the
statement of the first are as defined by~(\ref{eq:kappa}).

\begin{lemma}
\label{lem:compare}
Let $q=m/n\in\Q\cap(0,1/2)$. For each integer $r$ with $1\le r\le m$,
the word
\[10^{\k_r(q)}1^20^{\k_{r+1}(q)}1^2\ldots 1^20^{\k_m(q)}1\]
 disagrees with the word \[10^{\k_1(q)-1}1^20^{\k_2(q)}1^2\ldots
1^20^{\k_m(q)}1\] within the shorter of their lengths, and is greater
than it in the unimodal order. 
\end{lemma}

\begin{proof}
(See also~\cite{Stars}, Lemma~63.) If $r=1$ then the result is
obvious, so suppose that $1<r\le m$.  If the two words didn't
disagree, then it would follow that $\kappa_1(q)=\kappa_r(q)+1$ and
that $\kappa_m(q)=\kappa_{m-r+1}(q)$, contradicting the fact that
$c_q$ is palindromic.

Observe that formula~(\ref{eq:kappa}) gives, for each $s$ with
$1\le s\le m+1-r$,
\[(\kappa_1(q)-1)+\sum_{i=2}^s\kappa_i(q)=\left\lfloor\frac{s}{q}\right\rfloor-2s,\]
and
\begin{eqnarray*}
\sum_{i=r}^{r+s-1}\kappa_i(q)
&=&\left\lfloor\frac{r+s-1}{q}\right\rfloor-\left\lfloor
   \frac{r-1}{q}\right\rfloor-2s\\
&\ge&\left\lfloor \frac{r+s-1}{q}-\frac{r-1}{q}\right\rfloor
   -2s\\
&=&\left\lfloor\frac{s}{q}\right\rfloor-2s.
\end{eqnarray*}
Hence at the point where they first disagree,
$10^{\k_r(q)}1^20^{\k_{r+1}(q)}1^2\ldots 1^20^{\k_m(q)}1$ has a longer
block of $0$s than $10^{\k_1(q)-1}1^20^{\k_2(q)}1^2\ldots
1^20^{\k_m(q)}1$, and so is greater in the unimodal order.
\end{proof}

\begin{lemma}
\label{lem:lem2}
Let $s\in\MIA$ have height $q=q(s)>0$, and suppose $s=c_qd$ for some
$d\in\{0,1\}^\N$. Then $\sigma(d)\succeq\sigma^2(s)$ and $d\preceq
1\sigma^2(s)$.
\end{lemma}

\begin{proof}
Suppose first that the first symbol of $d$ is $0$. In this case it is
obvious that $d\preceq 1\sigma^2(s)$, so, writing $d=0e$, it is only
necessary to show that
\mbox{$e\succeq\sigma^2(s)=\widehat{w}_q0e$}. Now
$\sigma^n(s)=10e\preceq s$, since $s$ is a kneading sequence, and
$s=10\ldots$, so $\sigma^{n+2}(s)=e\succeq \sigma^2(s)$ as required.

If $d=1e$ for some $e{\in}\{0,1\}^\N$ then
$\sigma(d)\succeq\sigma^2(s)$ clearly implies that \mbox{$d\preceq
1\sigma^2(s)$}, so again it is required to show that
$e\succeq\sigma^2(s)=\widehat{w}_q1e$.

Assume for a contradiction that $e\prec \widehat{w}_q1e$. Now
$s=c_q1e\preceq c_q(1\widehat{w}_q)^\infty$ by
Theorem~\ref{thm:heightinterval}, and so
$e\succeq(\widehat{w}_q1)^\infty$ (since $c_q1$ is an odd word). Thus
\[(\widehat{w}_q1)^\infty\preceq e\prec \widehat{w}_q1e,\]
and so $e$ has $\widehat{w}_q1$ as an initial word. Let $k\ge 1$ be
the greatest integer such that $e=(\widehat{w}_q1)^kf$ for some
$f\in\{0,1\}^\N$ (such a greatest integer exists since if
$e=(\widehat{w}_q1)^\infty$ then it is not true that
$e\prec\widehat{w}_q1e$). Since $\widehat{w}_q1$ is an even word, this
gives
\[(\widehat{w}_q1)^\infty\preceq f\prec \widehat{w}_q1f,\]
so $f=\widehat{w}_q1\ldots$, contradicting the definition of $k$.
\end{proof}

\begin{lemma}
\label{lem:goodorbit}
Let $s\in\MIA$ have height $q=q(s)>0$. Then every $t\in\{0,1\}^\N$ on
the $\sigma$-orbit of $(w_q1)^\infty$ of the form $t=1^k0\ldots$,
where $k$ is odd, satisfies $t\succeq 1\sigma^2(s)$.
\end{lemma}
\begin{proof}
If $s=c_qd$ for some $d\in\{0,1\}^\N$, then
$d\preceq 1\sigma^2(s)=1\widehat{w_q}d$ by Lemma~\ref{lem:lem2}, so by
induction using the fact that $1\widehat{w_q}$ is an even word,
$d\preceq (1\widehat{w_q})^nd$ for all $n\ge 0$, and hence $d\preceq
(1\widehat{w_q})^\infty$ and
$1\sigma^2(s)=1\widehat{w_q}d\preceq(1\widehat{w_q})^\infty$.  If
$s=(w_q1)^\infty$, then $1\sigma^2(s)\preceq (1\widehat{w_q})^\infty$
by direct comparison. It therefore suffices to show that each such~$t$
satisfies $t\succeq (1\widehat{w_q})^\infty$. Recall that
$\rhe(m/n)$ is preperiodic to~$\lhe(m/n)$,
ie, $(1\widehat{w_q})^\infty$ lies on the $\sigma$-orbit of
$(w_q1)^\infty$.

Writing $q=m/n$ and writing $\k_i$ for $\k_i(q)$, 
\begin{eqnarray*}
(w_q1)^\infty&=&(10^{\k_1}1^20^{\k_2}1^2\ldots 1^2
0^{\k_{m-1}}1^20^{\k_m-1}1)^\infty\quad\text{ and}\\
(1\widehat{w_q})^\infty&=&(10^{\k_1-1}1^20^{\k_2}1^2\ldots 1^2
0^{\k_{m-1}}1^20^{\k_m}1)^\infty,
\end{eqnarray*}
and hence it is required to show that
\[(10^{\k_m-1}1^20^{\k_1}1^20^{\k_2}1^2\ldots1^20^{\k_{m-1}}1)^\infty
\succeq
(10^{\k_1-1}1^20^{\k_2}1^2\ldots1^20^{\k_{m-1}}1^20^{\k_m}1)^\infty,\]
and
\[(10^{\k_j}1^2\!\ldots 1^20^{\k_m-1}1^20^{\k_1}1^2\!\ldots 1^2 0^{\k_{j-1}}1)^\infty
\succeq
(10^{\k_1-1}1^20^{\k_2}1^2\!\ldots 1^20^{\k_{m-1}}1^20^{\k_m}1)^\infty\]
for $1\le j\le m-1$.

For the former, since $\k_m-1=\k_1-1$, it is equivalent to show that
\[(10^{\k_1}1^20^{\k_2}1^2\ldots1^20^{\k_{m-1}}1^20^{\k_m-1}1)^\infty
\succeq
(10^{\k_2}1^2\ldots1^20^{\k_{m-1}}1^20^{\k_m}1^20^{\k_1-1}1)^\infty,\]
which is immediate since the left hand side is $(w_q1)^\infty$, a
kneading sequence, and the right hand side is a shift of it.

For the latter, Lemma~\ref{lem:compare} gives that either
$10^{\k_j}1^2\ldots1^20^{\k_m-1}1$ disagrees with
$10^{\k_1-1}1^2\ldots 1^20^{\k_{m-j+1}}1$ and is greater than it
(which establishes the result), or the two are equal. If they are
equal, then removing this even word from the front of each side of the
inequality leaves $(w_q1)^\infty$ on the left hand side and a shift
of it on the right hand side, and the result follows.
\end{proof}

\begin{proof}[Proof of Theorem~\ref{thm:outside}]
Recall that $f$ has no homtervals, being the quotient of an MIA
unimodal thick interval map, and hence a point $x\in I$ can be
specified uniquely by its itinerary $i(x)$.

Let $x\in(a,b)\setminus\{c\}$ have itinerary
$i(x)=i_0i_1i_2\ldots$. Then $\theta(x_l)=f(x)_l$ if $i_0=0$, and
$\theta(x_l)=f(x)_u$ if $i_0=1$. Similarly, $\theta(x_u)=f(x)_l$ if
$x\ge p$ (ie, if $i(x)\succeq 1\sigma^2(s)$), while $x_u\in\gamma$
if $x\prec 1\sigma^2(s)$. Hence the least $n>0$ such that
$\theta^n(\hat{a})\in\overline{\gamma}$ is the least $n>0$ such that
either $\sigma^n(i(a))=i(a)$, or $\sigma^{n-k-1}(i(a))$ is of the form
$01^{k}\ldots$ for some odd $k$ with $\sigma^n(i(a))\preceq
1\sigma^2(s)$.

Suppose first that $\lhe(m/n)\prec s\prec\rhe(m/n)$, so in particular
$s=c_qd$ for some $d\in\{0,1\}^\N$ by Lemma~\ref{lem:kiq}. Thus 
\begin{eqnarray*}
i(a)&=&0^{\k_1}1^20^{\k_2}1^2\ldots
1^20^{\k_m}1d\quad\text{ and}\\
1\sigma^2(s)&=&10^{\k_1-1}1^20^{\k_2}1^2\ldots 1^20^{\k_m}1d
\end{eqnarray*}
(where $\k_i=\k_i(q)$). Since $10^{\k_j}1^2\ldots 1^20^{\k_m}1d\succ
1\sigma^2(s)$ for $2\le j\le m$ by Lemma~\ref{lem:compare}, it follows
that $\theta^i(\hat{a})\not\in\overline\gamma$ for $0<i<n$. On the
other hand, $\sigma^n(i(a))=d\preceq 1\sigma^2(s)$ by
Lemma~\ref{lem:lem2}, and hence $\theta^n(\hat{a})\in\overline\gamma$
as required.

To show the different cases for $\theta^n(\hat{a})$, it is only
necessary to translate the statements on the left hand sides into
statements about $i(f^n(a))=d$, and then convert these to equivalent
statements about $s=c_qd$. For instance, 
\begin{eqnarray*}
\theta^n(\hat{a})\in(\hat{a},c_u)_\gamma &\iff& d\prec i(c)=1s \\
&\iff&  s\prec c_q1s\\
&\iff& s\prec(c_q1)^\infty =\NBT(m/n).
\end{eqnarray*}
Cases ii) and v) follow similarly. 

Lemma~\ref{lem:goodorbit} translates directly into the statement that
$\theta$ has a periodic orbit $\tilde{P}$ above the period~$n$ orbit $P$
of~$f$ containing the point with itinerary $(w_q1)^\infty$ (the points
of $\tilde{P}$ on the upper half-circle are exactly those with
itineraries having $1^k0$ as an initial word for some odd~$k$). This
itinerary is known to correspond to a periodic orbit of rotation
type $m/n$.

Now let $x_l$ be the point of~$\tilde{P}$ with 
\[i(x)=(0^{\k_1}1^20^{\k_2}1^2\ldots 1^20^{\k_m-1}1^2)^\infty,\]
and consider the interval $J_0=[\hat{a},x_l)$. Then $J_0$ contains no
fixed points of~$\theta^n$: for the leftmost point $y_l$ on the
corresponding periodic orbit would have $i(a)\prec i(y)\prec i(x)$ and
hence $i(y)=0^{\k_1}1^2\ldots1^20^{\k_m-1}\opti1\ldots$, so either
$y=x$ or \mbox{$i(y)=(0^{\k_1}1^2\ldots1^20^{\k_m-1}01)^\infty$}, which
contains an isolated $1$ and therefore does not correspond to an orbit
of $\theta$. So $\theta^n$ is continuous on $J_0$, fixes $x_l$, maps
$\hat{a}$ into $\gamma$, and has no fixed point in $J_0$, and hence
every point of $J_0$ falls into $\gamma$ under iteration of
$\theta^n$, ie, $\Lambda\cap J_0=\emptyset$.

A similar argument shows that if $y_u$ is the point of $\tilde{P}$ with 
\[i(y)=(10^{\k_1-1}1^20^{\k_2}1^2\ldots 1^20^{\k_m}1)^\infty\] and
$J_1=[p_u,y_u)$, then $\Lambda\cap J_1=\emptyset$.  Thus
$J=\theta(J_0)\cup\theta(J_1)=(f(x)_l,f(y)_l)$ is disjoint from
$\Lambda$. In particular, the endpoints of $J$ are consecutive points
of $\tilde{P}$, and so
$\bigcup_{i=0}^{n-1}\theta^i(J)=S^1\setminus\tilde{P}$, establishing
that $\Lambda=\tilde{P}$ as required.

The two special cases $s=\lhe(m/n)$ and $s=\rhe(m/n)$ can be treated
similarly, but more straightforwardly since explicit expressions
for~$s$ are available. The reason that $\Lambda\neq \tilde{P}$ in
these cases is that $\tilde{P}$ contains the unique point
$\theta(\hat{a})$ of~$S^1$ that has two $\theta$-preimages, $\hat{a}$
and $p_u$.
\end{proof}

\section{Invariant unimodal generalized train tracks}
\label{sec:unitt}

In this section, the invariant generalized train tracks corresponding
to elements of~$\MIA$ are described
explicitly. This is achieved by a relatively straightforward analysis
of the construction of invariant generalized train tracks given in the
proof of Theorem~\ref{thm:ttexist}: the important point for what
follows is the way in which the structure of the train track is
governed by the height of the kneading sequence.

The case of preperiodic kneading sequences is somewhat more
complicated than the periodic case, and as such the two are treated
separately.

\subsection{The periodic case}
Let $s=\k(f)\in\MIA$ be the kneading sequence of a unimodal map $f$
whose critical point~$c$ is periodic of period~$N$, and let
$F\co(\sph,\ofi,A)\barrow$ be the associated thick interval map,
where~$A$ is the set of attracting periodic points of~$F$, consisting
of a single periodic orbit whose points are in natural correspondence
with the points of the $f$-orbit of~$c$ (see
Figure~\ref{fig:1001011thtrmap} for an example in the case
$s=(1001011)^\infty$). Let $\tau\subseteq\ofi\setminus A$ be the
$F$-invariant generalized train track: the aim of this section is to
describe~$\tau$ explicitly. Label the junctions of $\ofi$ with
integers $1$ to $N$ from left to right: thus junctions $1$ and
$N$ are $1$-junctions, while junctions $2$ to $N-1$ are
$2$-junctions.

Because of the convention adopted in
Section~\ref{subsec:unimodal} for the itinerary of a point whose orbit
contains $c$, the kneading sequence is given by $s=w^\infty$, where
$w$ is a word of length~$N$ whose final symbol is~$1$.

It will be seen (Theorem~\ref{thm:pertt} below) that only four
different basic configurations of infinitesimal edges can occur in any
given junction of $\ofi$. These are as follows: note that the position
of the puncture in the junction relative to the infinitesimal edges is
also specified (indicated by a small circle in the figures).

\begin{figure}[ht!]
\begin{center}
\pichere{0.5}{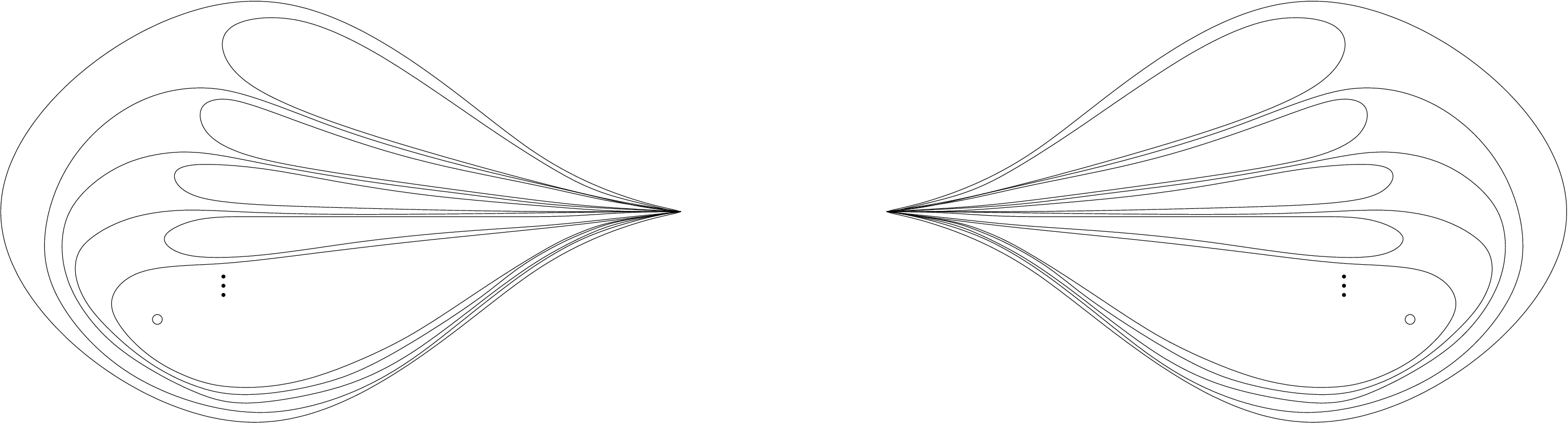}
\end{center}
\caption{$W^+$ and $W^-$}
\label{fig:w+-}
\end{figure}

\begin{figure}[ht!]
\begin{center}
\pichere{0.2}{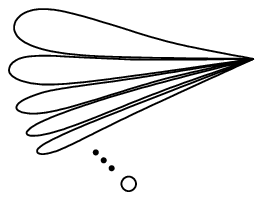}
\end{center}
\caption{$S^+$}
\label{fig:bubblebunch}
\end{figure}

\begin{figure}[ht!]

\begin{center}
\pichere{0.3}{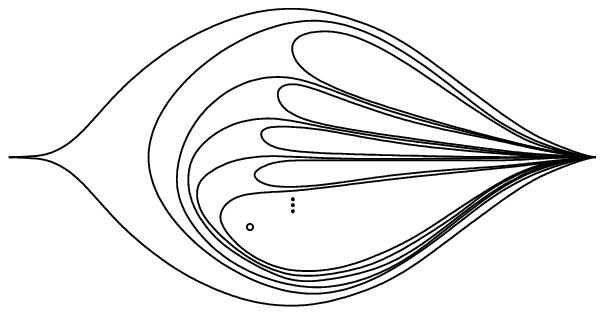}
\end{center}
\caption{$V_3^+$}
\label{fig:v3}
\end{figure}

\begin{description}
\item[$\BP$] A single bubble containing the puncture.
\item[$W$] A configuration of infinitely many loops of which
infinitely many are bubbles (all homotopically trivial), and
infinitely many are not, as depicted in Figure~\ref{fig:w+-}. There
are two different versions, $W^+$ and $W^-$ which are mirror images of
each other.
\item[$S^+$] A semi-infinite bouquet of homotopically trivial bubbles
as depicted in Figure~\ref{fig:bubblebunch}.
\item[$V_3$] A single bigon containing $W$. There are two versions,
$V_3^+$ and $V_3^-$, containing $W^+$ and $W^-$ respectively (see
Figure~\ref{fig:v3}: $V_3^-$ is the mirror image of $V_3^+$).
\end{description}


The infinitesimal edges of the train tracks described in this section
are specified as follows. For each junction of $\ofi$, there is given
one of the above symbols describing the basic configuration of
infinitesimal edges. The junction is said to be {\em type $\BP$}, {\em
type $W^+$}, {\em type $W^-$}, {\em type $S^+$}, {\em type $V_3^+$} or
{\em type~$V_3^-$} correspondingly.  For each $2$-junction, there is
in addition a symbol~$L$ or~$R$, specifying whether the loop(s) in
the basic configuration are attached to the left or right switch in
that junction. The $2$-junction is said to be {\em type L} or {\em
type R} correspondingly.

In each $2$-junction not of type~$V_3^\pm$, there is one additional
infinitesimal edge, whose endpoints are the two switches on the
boundary of the junction. This edge passes above (respectively below)
all other infinitesimal edges in the junction if the junction is type
R (respectively type L).

The reader seeking clarification as to how this notation is used can
look ahead to Examples~\ref{ex:pertt} and the accompanying figures.

There are two special cases: the {\em finite order} case
$s=\lhe(m/n)$ (where each junction is of type~$S^+$), and the {\em
NBT} case $s=\NBT(m/n)$ (where each junction is of type $\BP$, and
hence $\tau$ is finite). In all other cases, the invariant generalized
train track contains junctions of types $W^\pm$ and $V_3^\pm$.

In all cases, it can be determined which $2$-junctions are type $L$
and which are type $R$ by applying the following simple algorithm. The
idea is straightforward: the infinitesimal edges in junction~$1$ are
clearly attached to the right hand edge of that junction. This
information is propagated around the orbit, with the type changing
after each symbol~$1$ (corresponding to a `flip'), and remaining
unchanged after each symbol~$0$.

\begin{alg}
\label{alg:percase}
Let $s\in\MIA$ be a period $N$ kneading sequence, and let $\pi=\pi_s\in
S_N$ denote the induced permutation on the points of the critical
orbit of a corresponding unimodal map. Then the partition
$\{2,\ldots,N-1\}=L\cup R$ is determined inductively as follows:
$\pi(1)\in R$, and for each $r$ with $2\le r\le N-2$, $\pi^r(1)$ lies
in the same set as $\pi^{r-1}(1)$ if $\pi^{r-1}(1)<\pi\I(N)$, and
in the other set if $\pi^{r-1}(1)>\pi\I(N)$.
\end{alg}

\begin{example}
The easiest way to determine these sets is to write down the first $N$
symbols of $s$, and to place $L$s and $R$s above them.  The $LR$
sequence starts with an $R$ at the third symbol, and changes after
each $1$. 

For example, let $s=(10011001011)^\infty$, so that
$\pi=(1\,3\,7\,10\,2\,5\,9\,4\,8\,6\,11)$. Then
Algorithm~\ref{alg:percase} gives
\[\begin{array}{ccccccccccc}
11 & 1 & 3 & 7 & 10 & 2 & 5 & 9 & 4 & 8 & 6\\
&&R&R&L&R&R&R&L&L&R\\
1&0&0&1&1&0&0&1&0&1&1.
\end{array}\]
Thus $L=\{4,8,10\}$ and $R=\{2,3,5,6,7,9\}$.
\end{example}

The next result describes the invariant generalized train track
corresponding to any periodic kneading sequence in $\MIA$.

\begin{thm}
\label{thm:pertt}
Let $s\in\MIA$ be a periodic kneading sequence of period~$N$ with
height $m/n>0$. Then the invariant generalized train
track corresponding to $s$ has infinitesimal edges as
follows. Junction $i$ is of type $L$ or $R$ according as $i\in L$ or
$i\in R$. The basic configurations of infinitesimal edges are:
\begin{enumerate}[\rm a)]
\item If $s=\lhe(m/n)$, then all junctions are of type $S^+$.
\item If $s=\NBT(m/n)$, then all junctions are of type $\BP$.
\item Otherwise, let $\pi\in S_N$ denote the induced permutation on
the points of the periodic critical orbit of a unimodal map with
kneading sequence $s$, and let $\epsilon\in\{+,-\}$ be given by
$\epsilon=+$ if $\pi^{-1}(N)\in R$, and $\epsilon=-$ if
$\pi^{-1}(N)\in L$. Then for $0\le r\le N-1$, the junction $\pi^r(N)$
is of type $W^\epsilon$ if $r\le n+1$, and of type $V_3^\epsilon$
if $r>n+1$.
\end{enumerate}
\end{thm}

The proof is delayed until some explanatory examples have been presented.

\begin{examples}
\label{ex:pertt}
The three examples correspond to the three cases of
Theorem~\ref{thm:pertt}.
\begin{enumerate}[a)]
\item Let $s=(101)^\infty=\lhe(1/3)$, so that $\pi=(1\,2\,3)$. 
Algorithm~\ref{alg:percase} gives
\[\begin{array}{ccc}
3&1&2\\
&&R\\
1&0&1
\end{array}\] 
Thus the invariant generalized train track is described by
\[(S^+\,\,;\,\,S^+,R\,\,;\,\,S^+)\]
(see Figure~\ref{fig:101itt}: note that the infinitesimal edge with
endpoints the left and right switches in junction~2 passes above the
other infinitesimal edges in the junction, since the junction is of
type R. In this and subsequent figures, the punctures have been
labelled to clarify the dynamics.)

\begin{figure}[ht!]
\begin{center}
\lab{f}{\phi}{}
\lab{1}{1}{}
\lab{2}{2}{}
\lab{3}{3}{}
\pichere{0.6}{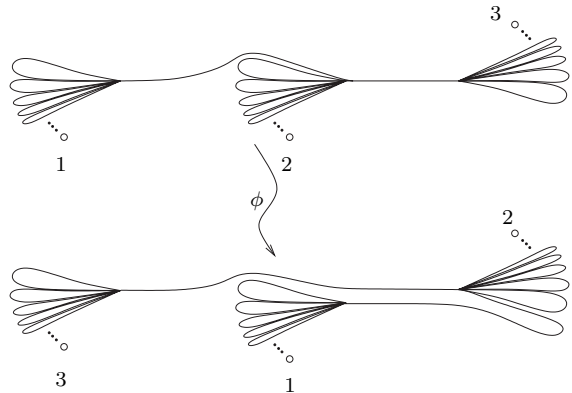}
\end{center}
\caption{The generalized train track corresponding to
$s=\lhe(1/3)$ and its image}
\label{fig:101itt}
\end{figure}

\item Let $s=(10011)^\infty=\NBT(1/3)$, so 
$\pi=(1\,2\,4\,3\,5)$. Algorithm~\ref{alg:percase} gives
\[\begin{array}{ccccc}
5&1&2&4&3\\
&&R&R&L\\
1&0&0&1&1
\end{array}\] 
Thus the invariant generalized train track is described by
\[(\BP\,\,;\,\,\BP,R\,\,;\,\,\BP,L\,\,;\,\,\BP,R\,\,;\,\,\BP)\]
(see Figure~\ref{fig:10010tt}).
\begin{figure}[ht!]
\begin{center}
\lab{f}{\phi}{}
\lab{1}{1}{}
\lab{2}{2}{}
\lab{3}{3}{}
\lab{4}{4}{}
\lab{5}{5}{}
\pichere{0.8}{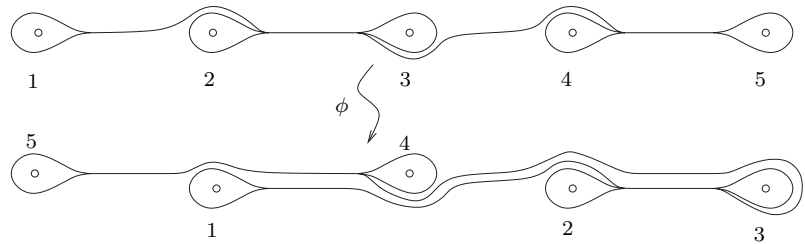}
\end{center}
\caption{The generalized train track corresponding to
$s=\NBT(1/3)$ and its image}
\label{fig:10010tt}
\end{figure}

\item Let $s=(1001011)^\infty$, so $q(s)=1/3$ and
$\pi=(1\,3\,6\,2\,5\,4\,7)$. Algorithm~\ref{alg:percase} gives
\[\begin{array}{ccccccc}
7&1&3&6&2&5&4\\
&&R&R&L&L&R\\
1&0&0&1&0&1&1
\end{array}\] 
Thus $\epsilon=+$, since $\pi^{-1}(7)=4\in R$. The junctions $\pi^r(7)$
with $0\le r\le 4$ are $7$, $1$, $3$, $6$, and $2$, so these junctions
are of type $W^+$, while junctions $5$ and $4$ are of type $V_3^+$.
Thus the invariant generalized train track is described by
\[(W^+\,\,;\,\,W^+,L\,\,;\,\,W^+,R\,\,;\,\,V_3^+,R\,\,;\,\,
V_3^+,L\,\,;\,\,W^+,R\,\,;\,\,W^+)\]
(see Figure~\ref{fig:1001011}).

\begin{figure}[ht!]
\begin{center}
\lab{f}{\phi}{}
\lab{1}{1}{}
\lab{2}{2}{}
\lab{3}{3}{}
\lab{4}{4}{}
\lab{5}{5}{}
\lab{6}{6}{}
\lab{7}{7}{}
\stretchpichere{.95}{5cm}{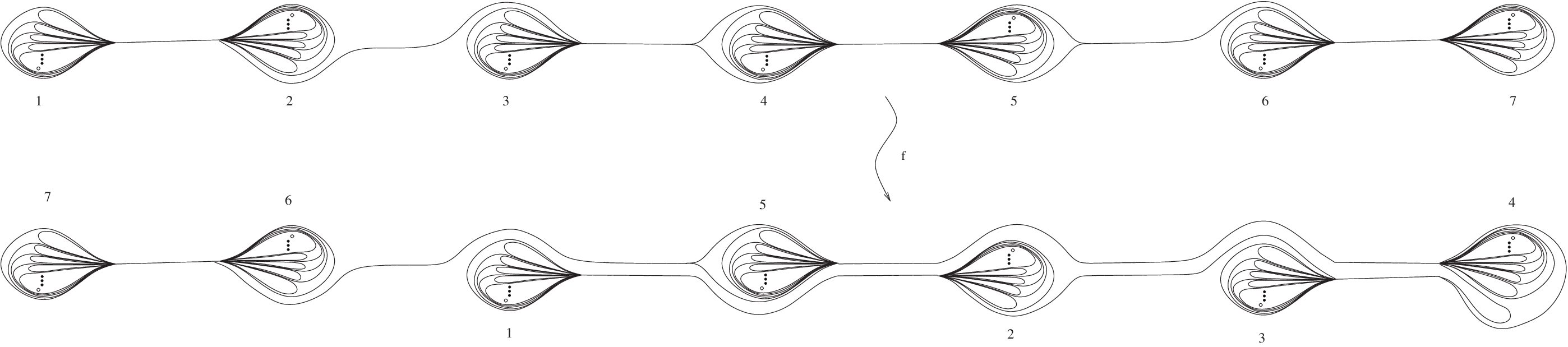}
\end{center}
\caption{The generalized train track corresponding to
$s=(1001011)^\infty$ and its image}
\label{fig:1001011}
\end{figure}
\end{enumerate}
\end{examples}

\begin{proof}[Proof of Theorem~\ref{thm:pertt}]
Let $\tau$ be the invariant generalized train track, as constructed
using the algorithm of Theorem~\ref{thm:ttexist}. Let
$A\subseteq\{2,\ldots,N-1\}$ (respectively
$B\subseteq\{2,\ldots,N-1\}$) be the set of $2$-junctions which
contain an infinitesimal edge of $\tau$ joining the two switches in
the junction and passing above (respectively below) the
puncture. There are two ways such edges can arise during the
construction. First, as images of real edges: this requires that $r\in
A$ (respectively $r\in B$) if $r\ge\pi(1)$ and
$\pi^{-1}(r)<\pi^{-1}(N)$ (respectively $r\ge\pi(1)$ and
$\pi^{-1}(r)>\pi^{-1}(N)$). Second, as images of other such
infinitesimal edges: this requires that if $\pi^{-1}(r)<\pi^{-1}(N)$ then
$r$ lies in $A$ (respectively $B$) if $\pi^{-1}(r)$ lies in $A$
(respectively $B$); while if $\pi^{-1}(r)>\pi^{-1}(N)$ then
$r$ lies in $B$ (respectively $A$) if $\pi^{-1}(r)$ lies in $A$
(respectively $B$).

To express this symbolically, partition $\{1,\ldots,N\}=O\cup
I$, with $r\in O$ if $r<\pi^{-1}(N)$ (ie, if $r$ corresponds to the
symbol $0$), and $r\in I$ otherwise. Then $A$ and $B$ are the smallest
subsets of $\{2,\ldots,N-1\}$ such that
\begin{enumerate}[i)]
\item If $r\ge\pi(1)$, then $r\in A$ if $\pi^{-1}(r)\in O$, and $r\in
B$ if $\pi^{-1}(r)\in I$.
\item If $\pi^{-1}(r)\in O$ then $r\in A$ if $\pi^{-1}(r)\in A$, and
$r\in B$ if $\pi^{-1}(r)\in B$. If $\pi^{-1}(r)\in I$ then $r\in
B$ if $\pi^{-1}(r)\in A$, and $r\in A$ if $\pi^{-1}(r)\in B$.
\end{enumerate}
By condition~ii), $A\cap B$ is a final segment of the sequence
$(\pi^j(1))_{j=1}^{N-2}$. By conditions~i) and~ii), the least $j$ such
that $\pi^j(1)\in A\cap B$ is the least $j$ with $\pi^j(1)>\pi(1)$ and
$\pi^{j-1}(1)\in B$.  The main step of the proof is to establish that
this least $j$ is equal to $n+1$. In particular, by
Lemma~\ref{lem:kiq}, $A\cap B=\emptyset$ if and only if either
$s=\lhe(m/n)$ or $s=\NBT(m/n)$.

By Lemma~\ref{lem:kiq}, either $s=\lhe(m/n)$ or $s=c_{m/n}d$ for
some $d\in\{0,1\}^\N$. Assume that the latter holds: the proof that
$A\cap B=\emptyset$ in the former case is similar. Now
since $1\in O$, condition~i) gives $\pi(1)\in A$, and the assignments
to $A$ and $B$ given by condition~ii) are as follows:

\[s=
\begin{array}{ccccccccccccc}
 & A & A & B & A & A & B & & A & B & A & A & \\
1 & 0^{\kappa_1} & 1 & 1 & 0^{\kappa_2} & 1 & 1 & \ldots & 1 & 1 &
0^{\kappa_m} & 1 & d
\end{array}
\]
Using
$i(\pi(1))=0^{\k_1-1}1^20^{\k_2}1^2\ldots1^20^{\k_m}1d=\sigma^2(s)$
and applying Lemma~\ref{lem:compare} gives $j\ge n+1$. Since the first
symbol of $d$ corresponds to an element of $B$, it only remains to
show that $\sigma(d)\succeq \sigma^2(s)$. This is given by
Lemma~\ref{lem:lem2}.

Notice that for $r\not\in A\cap B$, comparing condition~ii) with
Algorithm~\ref{alg:percase} gives that $r\in A$ if and only if $r\in
R$.

Now consider how loops can arise during the construction of
$\tau$. Again, there are two ways. First, as the image of an
infinitesimal edge joining the two switches of junction $\pi^{-1}(N)$:
this gives rise to a homotopically trivial bubble in junction~$N$
(respectively a loop containing the puncture in junction~$N$) if the
edge passes above (respectively below) the puncture in
junction~$\pi^{-1}(N)$. Second, as the image of another loop.

If $s=\NBT(m/n)$, ie,
\[s=
\begin{array}{cccccccccccccc}
 & A & A & B & A & A & B & & A & B & A & A & B &\\
1 & 0^{\kappa_1} & 1 & 1 & 0^{\kappa_2} & 1 & 1 & \ldots & 1 & 1 &
0^{\kappa_m} & 1 & 1 & \ldots,
\end{array}
\]
then any infinitesimal edge joining the switches of junction
$\pi^{-1}(N)$ passes below the puncture, and hence gives rise to a
loop in junction~$N$ which contains the puncture. Statement~b)
follows easily. Similarly, if $s=\lhe(m/n)$, ie,
\[s=
\begin{array}{ccccccccccccc}
 & A & A & B & A & A & B & & A & B & A & A & \\
1 & 0^{\kappa_1} & 1 & 1 & 0^{\kappa_2} & 1 & 1 & \ldots & 1 & 1 &
0^{\kappa_m-1} & 1 & \ldots,
\end{array}
\]
then any infinitesimal edge joining the switches of junction
$\pi^{-1}(N)$ passes above the puncture, and hence gives rise to a
homotopically trivial bubble in junction~$N$. Statement~a) follows.

In any other case, $\pi^{-1}(N)\in A\cap B$, so that both
homotopically trivial bubbles and homotopically non-trivial loops are
created in junction~$N$, leading under iteration of the construction
to a configuration of infinitesimal edges of type~$W$ in each
junction. The orientation of these infinitesimal edges bubbles depends
on whether $\pi^{-1}(N)\in R$ or $\pi^{-1}(N)\in L$ (which determines
whether a new homotopically trivial bubble arises below or above the
existing loops in junction~$N$ during the
construction). Junction~$r$ is of type $V_3^\pm$ if $r\in A\cap B$,
but of type $W^\pm$ if $r$ is in only one of $A$ and $B$. Statement~c)
follows.
\end{proof}

\subsection{The preperiodic case}

The preperiodic case is somewhat more complicated than the periodic
case, though there are no new ideas introduced: as such, the treatment
is more informal and some of the details are omitted. The most
substantial modifications arising because the critical point is not
periodic are:
\begin{enumerate}[a)]
\item There is a unique loop which arises during the construction as
the image of an infinitesimal edge which is not a loop, and this loop
is a homotopically trivial bubble. All other loops are images of this
one, and hence all loops are homotopically trivial bubbles.
\item Only finitely many infinitesimal edges occur in those junctions
corresponding to the strictly preperiodic part of the orbit of the
critical point.
\item The parity of the periodic part of the kneading sequence
plays an important role. If there are an odd number of $1$s, then
there is a flip each time around the periodic orbit, and hence bubbles
are attached to both the left and right switches of $2$-junctions.
\item In the periodic case, the configuration of infinitesimal edges
changes (from $W$ to $V_3$) at a point in the orbit given by the
denominator of the height. In the preperiodic case, separate cases
arise according as this point lies in the preperiodic or periodic part
of the orbit.
\end{enumerate}

\begin{figure}[ht!]
\begin{center}
\pichere{0.3}{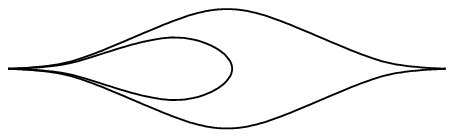}
\end{center}
\caption{$V_0$}
\label{fig:v0}
\end{figure}

\begin{figure}[ht!]
\begin{center}
\lab{v1-}{V_1^-}{}\lab{v1+}{V_1^+}{}
\pichere{0.6}{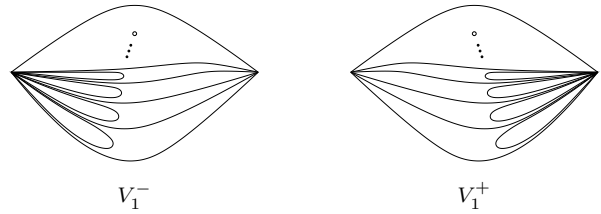}
\end{center}
\caption{$V_1^-$ and $V_1^+$}
\label{fig:v1}
\end{figure}

\begin{figure}[ht!]
\begin{center}
\lab{v2-}{V_2^-}{}
\pichere{0.3}{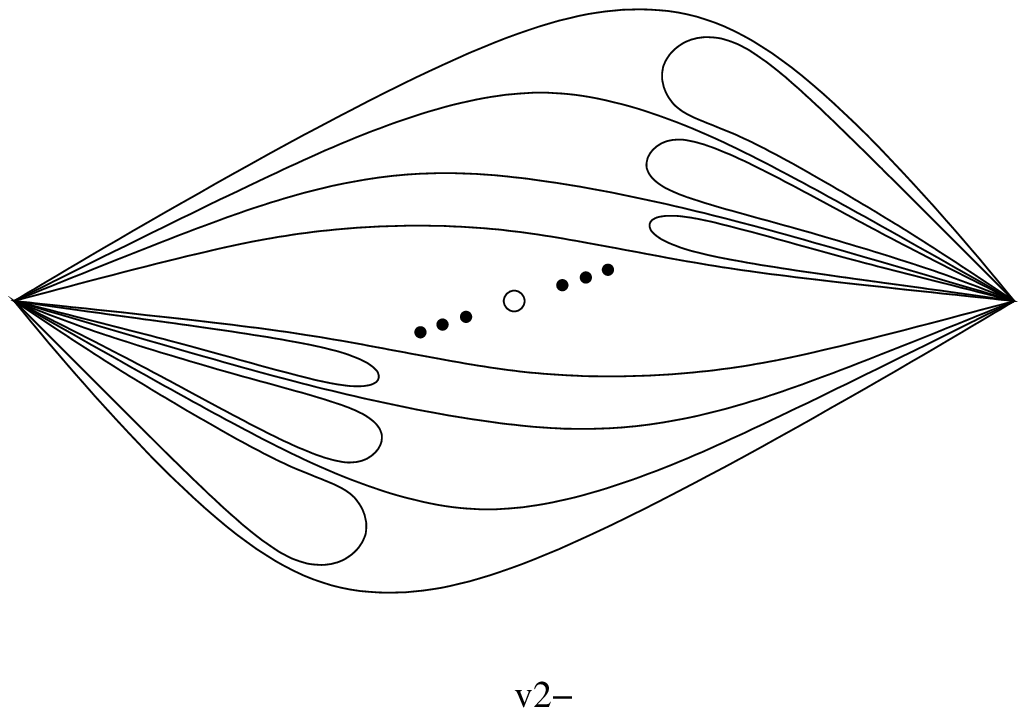}
\end{center}
\caption{$V_2^-$}
\label{fig:v2}
\end{figure}

\begin{figure}[ht!]
\begin{center}
\pichere{0.6}{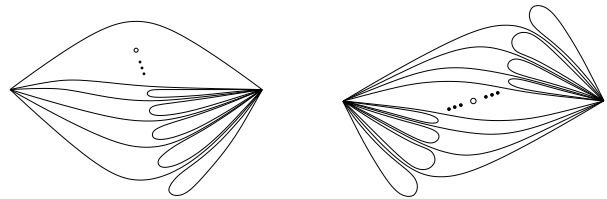}
\end{center}
\caption{Extra bubbles: the configurations $V_1^+\,B$
and $V_2^-\,B^2$}
\label{fig:extrabubbles}
\end{figure}

The basic configurations of infinitesimal edges are entirely different
from the periodic case. There are two distinct types: those which
arise in junctions corresponding to the preperiodic part of the orbit
($B$ and $V_0$), which contain only finitely many infinitesimal edges
and have no associated puncture; and those which arise in periodic
junctions, which contain infinitely many infinitesimal edges and have
an associated puncture. The configurations are:

\begin{description}
\item[$B$] A single (homotopically trivial) bubble.
\item[$V_0$] A bigon containing a single bubble, as depicted in
Figure~\ref{fig:v0}.
\item[$S^-$] The mirror image of $S^+$.
\item[$V_1$] A decorated bigon as depicted in
Figure~\ref{fig:v1}. There are two versions, $V_1^+$ and $V_1^-$,
which are mirror images of each other. In addition, an additional
external bubble may be present in the position indicated in
Figure~\ref{fig:extrabubbles}. In this case, the configuration is
denoted $V_1^\pm\,B$.
\item[$V_2$] A decorated bigon as depicted in
Figure~\ref{fig:v2}. There are two versions, $V_2^+$ and $V_2^-$:
$V_2^+$ is the mirror image of $V_2^-$. In addition, there may be 
two additional external bubbles in the positions indicated in
Figure~\ref{fig:extrabubbles}. In this case, the configuration is
denoted $V_2^\pm\,B^2$.
\end{description}

For each junction of $\ofi$, one of the above symbols is given
describing the basic configuration of infinitesimal edges: the
junction is said to be {\em type $V_0$}, {\em type $V_1^-$}, {\em type
$V_2^+\,B^2$}, etc. Each $2$-junction which is not of type $V_2$ is
either {\em type~$L$} or {\em type~$R$}, according as the bubbles are
attached to the left or right switch of the junction. For each
$2$-junction of type $B$ or $S^-$, there is one additional
infinitesimal edge, whose endpoints are the two switches on the
boundary of the junction. This edge passes above (respectively below)
all other infinitesimal edges in the junction if the junction is type
$R$ (respectively type $L$).

Throughout this subsection $s$ denotes a strictly preperiodic kneading
sequence, written $s=vw^\infty$, where~$v$ is taken to be as short a
word as possible (so the final symbol of $v$ is not equal to the final
symbol of $w$). The lengths of $v$ and $w$ will be denoted $k$ and $l$
respectively (so $k,l\ge 1$), and the size of the orbit of~$s$ under the
shift map is denoted $N=k+l$. The $N$ points of this orbit are labelled with
the integers $1$ through $N$ according to their relative unimodal
ordering (so point $r$ lies in junction $r$). The shift map $\sigma$
thus induces a map $\rho\co\{1,\ldots,N\}\barrow$. Note that $\rho$ is
not a permutation: the point~$N$ (corresponding to the first symbol of
$s$) has no preimage, and one point (corresponding to the first symbol
of $w$) has two preimages.

The following algorithm is the analogue of Algorithm~\ref{alg:percase}
for the preperiodic case. Notice that $L\cup R$ need not be a
partition of $\{2,\ldots,N-1\}$ (reflecting the possibility that the
word~$w$ may be odd): if $r\in L\cap R$, then bubbles are
attached to both switches in junction~$r$ (configuration $V_2$).

\begin{alg}
\label{alg:prepercase}
The subsets $L$ and $R$ of $\{2,\ldots,N-1\}$ are determined
inductively as follows.  Let $\rho^2(N)\in R$, and for each $r$ with
$3\le r\le N+l-1$, let $\rho^r(N)$ lie in the same set as
$\rho^{r-1}(N)$ if $s_{r-1}=0$, and in the other set if $s_{r-1}=1$.
\end{alg}

If $w$ is an even word, then the final $l$ steps in this algorithm are
unnecessary, and $L\cap R=\emptyset$.

\begin{examples}
\label{ex:preperLR}
\begin{enumerate}[a)]
\item Let $s=1000(101)^\infty$. Then $v=1000$, $w=101$, and
\[\rho=\left(
\begin{array}{ccccccc}
1 & 2 & 3 & 4 & 5 & 6 & 7\\
2 & 4 & 5 & 6 & 6 & 3 & 1
\end{array}
\right)\]
(ie, the $\rho$-orbit of $7$ is $7124(635)^\infty$). Thus
\[
\begin{array}{ccccccccl}
7 & 1 & 2 & 4 & & 6 & 3 & 5 & \\
 &  & R & R & & R & L & L & ,\\
1 & 0 & 0 & 0 & & 1 & 0 & 1 &
\end{array}
\]
and so $L=\{3,5\}$ and $R=\{2,4,6\}$.

\item Let $s=100101(10)^\infty$. Then $v=100101$, $w=10$, and
\[\rho=\left(
\begin{array}{cccccccc}
1 & 2 & 3 & 4 & 5 & 6 & 7 & 8\\
4 & 5 & 6 & 7 & 6 & 3 & 2 & 1
\end{array}
\right)\]
(ie, the $\rho$-orbit of $8$ is $814725(63)^\infty$). Thus
\[
\begin{array}{ccccccccccccl}
8 & 1 & 4 & 7 & 2 & 5 & & 6 & 3 & & 6 & 3 &\\
  &   & R & R & L & L & & R & L & & L & R &,\\
1 & 0 & 0 & 1 & 0 & 1 & & 1 & 0 & & 1 & 0 &
\end{array}
\]
and so $L=\{2,3,5,6\}$ and
$R=\{3,4,6,7\}$. Note that since $w$ is odd, the two repeating
points $3$ and $6$ of the orbit lie in both $L$ and $R$.

\end{enumerate}

\end{examples}

For preperiodic kneading sequences, there is only one special case:
that in which $s=\rhe(m/n)$.

\begin{thm}
\label{thm:prepertt}
Let $s=vw^\infty\in\MIA$ be a strictly preperiodic kneading sequence
of height $m/n>0$. Let $k$ and $l$ be the lengths of $v$
and $w$ (chosen as small as possible). Let $\rho$ be the induced map
on $\{1,\ldots,N=k+l\}$. Then the invariant generalized train track
corresponding to $s$ has infinitesimal edges as follows. If $i\in
L\setminus R$ (respectively $i\in R\setminus L$) then junction $i$ is
of type $L$ (respectively~$R$). The basic configurations of
infinitesimal edges are:
\begin{enumerate}[\rm a)]
\item If $s=\rhe(m/n)$, then the junction $\rho^r(N)$ is of type $B$
for $0\le r<k$, and of type $S^-$ for $k\le r<N$.
\item Otherwise, let $\epsilon\in\{+,-\}$ be given by
$\epsilon=+$ if $\rho^{k-1}(N)\in R$, and $\epsilon=-$ if
$\rho^{k-1}(N)\in L$. Then for each $r$ with $0\le r<N$:
\begin{enumerate}[\rm i)]
\item If $r\le n+1$, then junction $\rho^r(N)$ is of type
\begin{itemize}
\item $B$, if $r\le k$.
\item $V_1^\epsilon\,B$ or $V_2^\epsilon\,B^2$ (according as
$w$ is even or odd), if $r>k$.
\end{itemize}
\item If $r>n+1$, then junction $\rho^r(N)$ is of type
\begin{itemize}
\item $V_0$, if $r\le k$.
\item $V_1^\epsilon$ or $V_2^\epsilon$ (according as
$w$ is even or odd), if $r>k$.
\end{itemize}
\end{enumerate}
\end{enumerate}
\end{thm}

\begin{figure}[ht!]
\begin{center}
\lab{f}{\phi}{}
\lab{1}{1}{}
\lab{2}{2}{}
\lab{3}{3}{}
\lab{4}{4}{}
\lab{5}{5}{}
\pichere{0.8}{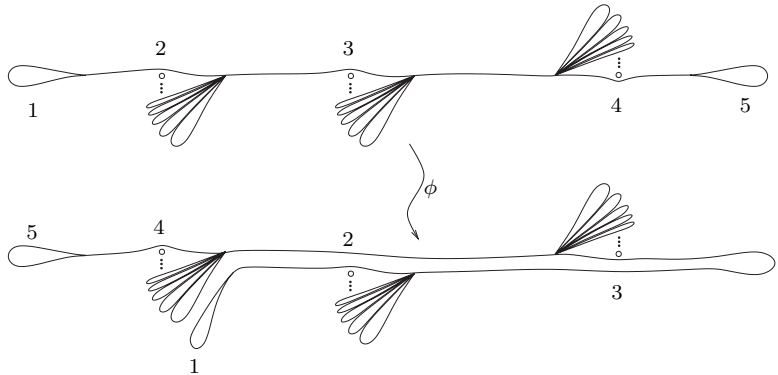}
\end{center}\vspace{-1.3cm}
\caption{The generalized train track corresponding to
$s=\rhe(1/3)$ and its image}
\label{fig:10b011}
\end{figure}

\begin{figure}[ht!]
\begin{center}
\lab{f}{\phi}{}
\lab{1}{1}{}
\lab{2}{2}{}
\lab{3}{3}{}
\lab{4}{4}{}
\lab{5}{5}{}
\lab{6}{6}{}
\lab{7}{7}{}
\pichere{.95}{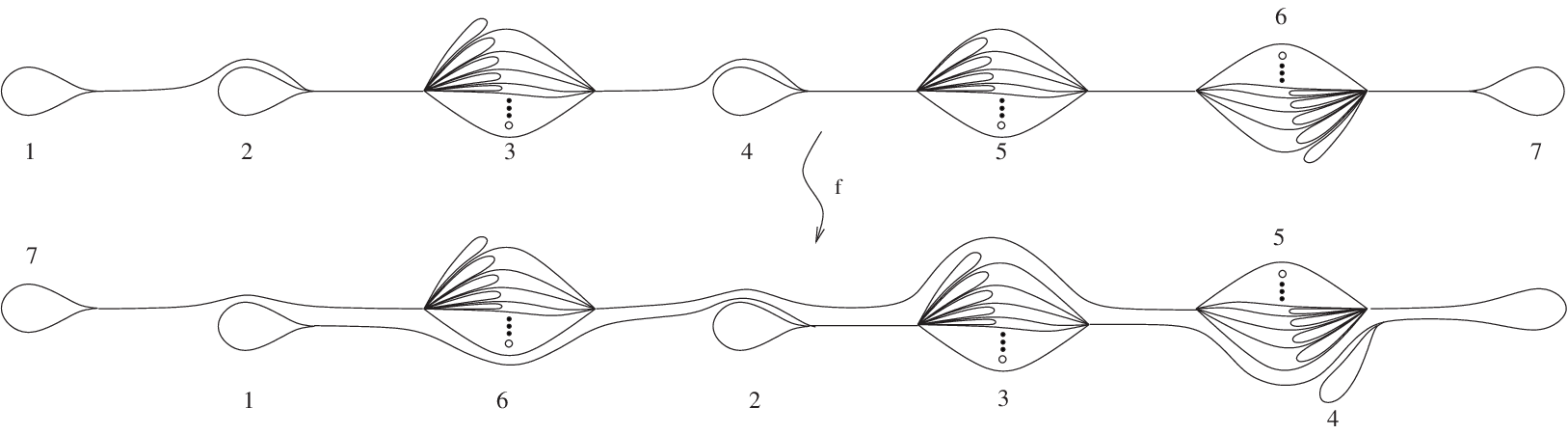}
\end{center}
\caption{The generalized train track corresponding to
$s=1000(101)^\infty$ and its image}
\label{fig:1000b101}
\end{figure}

\begin{figure}[ht!]
\begin{center}
\lab{f}{\phi}{}
\lab{1}{1}{}
\lab{2}{2}{}
\lab{3}{3}{}
\lab{4}{4}{}
\lab{5}{5}{}
\lab{6}{6}{}
\lab{7}{7}{}
\lab{8}{8}{}
\pichere{.95}{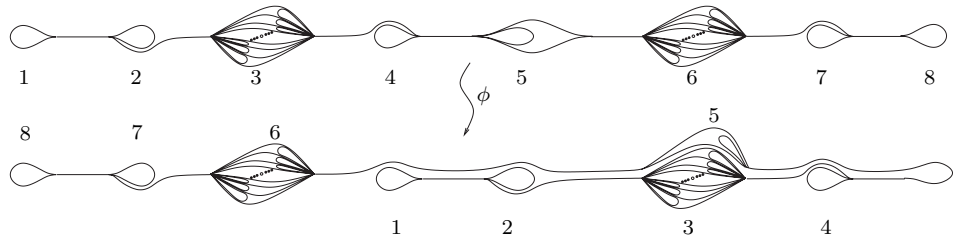}
\end{center}
\caption{The generalized train track corresponding to
$s=100101(10)^\infty$ and its image}
\label{fig:100101b10}
\end{figure}

\begin{examples}
The first example illustrates the case $s=\rhe(m/n)$. The other two
examples show the two possibilities $w$ even and $w$ odd.
\begin{enumerate}[a)]
\item Let $q=1/3$ and $s=\rhe(q)$ so
$s=1001(101)^\infty=10(011)^\infty$. Thus $v=10$, $w=011$, and
\[\rho=\left(
\begin{array}{ccccc}
1 & 2 & 3 & 4 & 5 \\
2 & 3 & 4 & 2 & 1 
\end{array}
\right)\]
(ie, the $\rho$-orbit of $5$ is $51(234)^\infty$). Thus
\[
\begin{array}{ccccccl}
5 & 1 & & 2 & 3 & 4 & \\
  &   & & R & R & L & ,\\
1 & 0 & & 0 & 1 & 1 &
\end{array}
\]
and so $R=\{2,3\}$ and $L=\{4\}$.

Since $k=2$ and $l=3$, junctions $N=5$ and $\rho(N)=1$ are of
type~$B$, while the other junctions are of type $S^-$.

Thus the invariant generalized train track is described by
\[(B\,\,;\,\,S^-,R\,\,;\,\,S^-,R\,\,;\,\,S^-,L\,\,;\,\,B)\]
(see figure~\ref{fig:10b011}: here all junctions have been labelled to
clarify the dynamics).

\item Let $s=\overset{7}{1}
\overset{1}{0}
\overset{2}{0}
\overset{4}{0}
(
\overset{6}{1}
\overset{3}{0}
\overset{5}{1})^\infty$. Then 
$L=\{3,5\}$ and $R=\{2,4,6\}$ (see Example~\ref{ex:preperLR}~a)).

Since $k=4$ and $l=3$, $\rho^{k-1}(N)=\rho^3(7)=4\in R$, and hence
$\epsilon=+$. $s$ has height $q(s)=1/4$, so $n=4$. Thus junctions
$N=7$, $\rho(N)=1$, $\rho^2(N)=2$, and $\rho^3(N)=4$ are of type $B$;
junctions $\rho^4(N)=6$ and $\rho^5(N)=3$ are of type $V_1^+\,B$; and
junction $\rho^6(N)=5$ is of type $V_1^+$.

Thus the invariant generalized train track is described by
\[(B\,\,;\,\,B,R\,\,;\,\,V_1^+\,B,L\,\,;\,\,B,R\,\,;\,\,
V_1^+,L\,\,;\,\,V_1^+\,B,R\,\,;\,\,B)\]
(see Figure~\ref{fig:1000b101}).

\item Let $s=\overset{8}{1} \overset{1}{0} \overset{4}{0}
\overset{7}{1} \overset{2}{0} \overset{5}{1} ( \overset{6}{1}
\overset{3}{0})^\infty$.  
Then $L=\{2,3,5,6\}$ and $R=\{3,4,6,7\}$ (see
Example~\ref{ex:preperLR}~b)).

Since $k=6$ and $l=2$, $\rho^{k-1}(N)=\rho^5(8)=5\in L$, and hence
$\epsilon=-$. $s$ has height $q(s)=1/3$, so $n=3$. Thus junctions
$N=8$, $\rho(N)=1$, $\rho^2(N)=4$, $\rho^3(N)=7$, and $\rho^4(N)=2$
are of type~$B$; junction $\rho^5(N)=5$ is of type $V_0$; and
junctions $\rho^6(N)=6$ and $\rho^7(N)=3$ are of type~$V_2^-$.

Thus the invariant generalized train track is described by
\[(B\,\,;\,\,B,L\,\,;\,\,V_2^-\,\,;\,\,B,R\,\,;\,\,
V_0,L\,\,;\,\,V_2^-\,\,;\,\,B,R\,\,;\,\,B)\]
(see Figure~\ref{fig:100101b10}).

\end{enumerate}
\end{examples}

\begin{proof}[Sketch proof of Theorem~\ref{thm:prepertt}]
The proof is very similar to that of Theorem~\ref{thm:pertt}. Defining
sets~$A$ and~$B$ as there, and noting that $s=c_{m/n}d$ for some
$d\in\{0,1\}^\N$, the same analysis shows that $A\cap B=\emptyset$ if
and only if $s=\rhe(m/n)$, and in all other cases $A\cap
B=\{\rho^j(N)\,:\,n+2\le j\le N-1\}$. This is again the essential part
of the proof.

When $s=\rhe(m/n)$, it is then straightforward to check case~a). In
other cases, it is immediate that for $r\le k$ (ie, in the
preperiodic part of the orbit) $r$ has type $V_0$ if $r\in A\cap B$,
and type $B$ otherwise.

All of the bubbles arise in the construction as images of the bubble
in junction~$N$ (which itself is created as the image of part of a
real edge). In particular, they are all homotopically trivial. The
combinatorics of the different cases can be checked routinely.
\end{proof}

\section{Unimodal generalized pseudo-Anosov maps}
\label{sec:genpA}

In this section the construction of generalized pseudo-Anosov maps
starting from the invariant train tracks of Section~\ref{sec:unitt} is
described. 

\subsection{Preliminaries}
Let $F\co(\sph,\ofi,A)\barrow$ be a unimodal MIA thick interval map,
where $A$ is the set of attracting periodic points of $F$. Let
$f\co I\barrow$ be the quotient unimodal map, and $\phi\co\tau\barrow$
be the associated generalized train track map. Let
$e_1,e_2,\ldots,e_{N-1}$ denote the real edges of $\tau$, and
$e_N,e_{N+1},\ldots$ denote the infinitesimal edges. Define the
transition matrix~$M=(m_{ij})$ by
\[ m_{ij}=\mbox{number of times\ } \phi(e_j) \mbox{\ crosses\ } e_i.\]
Thus~$M$ is an infinite matrix, unless $\k(F)=\NBT(m/n)$ for some
$m/n\in(0,1/2)$.

From the construction of the invariant generalized train track in
Section~\ref{sec:unitt}, it is immediate that each~$\phi(e_j)$ is a
finite length edge path; that each~$e_i$ is contained in the edge path
$\phi(e_j)$ for only finitely many~$j$; and that each infinitesimal
edge is mapped by~$\phi$ homeomorphically onto an infinitesimal
edge. These observations yield the following facts about the strucure
of the transition matrix:
\begin{enumerate}[a)]
\item $M$ is of the form
\[M=\left( \begin{array}{cc}
               A & 0 \\
               B & \Pi
             \end{array} \right), \]
where $A$ is an $(N-1)\times(N-1)$ matrix recording transitions between
real edges, $\Pi$ records transitions between infinitesimal edges, and
$B$ records transitions from real to infinitesimal edges.
\item All entries of $\Pi$ are either $0$ or $1$. Each column of $\Pi$
has exactly one non-zero entry, and each row has only finitely many
non-zero entries.
\item $B$ has only finitely many non-zero entries in each column.
\end{enumerate}

In particular, if~$M$ is regarded as an operator acting on~$l^1$,
then~$M$ is bounded with
$\|M\|_1\le\max_j\left\{\sum_i|m_{ij}|\right\}<\infty$.  Let $\lambda$
be the Perron-Frobenius eigenvalue of $A$, and~$Y\in\R^{N-1}$ the
associated positive eigenvector. Then~$Y$ can be extended to a
positive eigenvector $y=(Y\,Y')\in l^1$ of $M$, called the
Perron-Frobenius eigenvector of~$M$, by setting
\[
Y' =
\frac{1}{\lambda}(B+\frac{1}{\lambda}\Pi 
B+\frac{1}{\lambda^2}\Pi^2B+\ldots)Y.
\]
Note that $\|\Pi\|_1\le 1$ so that the series above converges in
$l^1$, implying that $y=(Y\,Y')\in l^1$.
Note also that the matrices $\Pi^kB$ represent transitions from a
real edge to an infinitesimal edge under~$\phi$, and then to another
infinitesimal edge under a further~$k$ iterations of~$\phi$. By the
construction of~$\tau$ described in the proof of
Theorem~\ref{thm:ttexist}, every infinitesimal edge of~$\tau$ is the
image, under some iterate of~$\phi$, of some infinitesimal edge
of~$F_*(\tau(\emptyset))$, and each of these infinitesimal edges is
crossed by $\phi(e)$ for some real edge~$e$. Thus for all $i\ge 1$
there exists $k\ge1$ such that row~$i$ of~$\Pi^kB$ is
non-zero. Since~$Y$ is a strictly positive vector, this establishes
that every entry of~$Y'$ is also positive.

\begin{defn}
A sequence $(y_i)$ of positive real numbers is said to satisfy the
{\em switch conditions for $\tau$} if, for each switch $q$ of $\tau$,
\[ y_{i_0}= \sum y_i + 2\sum y_j, \]
where $e_{i_0}$ is the real edge with endpoint $q$, and the first and
second sums range over the sets of indices of infinitesimal edges
having one and both endpoints at~$q$ respectively.
\end{defn}

The proof of the following lemma is essentially the same as that of
the corresponding fact for finite matrices in Section~3.4
of~\cite{BH}.

\begin{lemma} 
\label{lem:BH}
Let the transition matrix~$M$ of
$\phi\co\tau\barrow$ have Perron-Frobenius eigenvector
$y=(y_1\,y_2\,\ldots)$. Then $y$ satisfies the switch conditions for
$\tau$. 
\end{lemma}
\begin{proof}
For $r,s,k\ge 1$ denote by $m^{(k)}_{rs}$ the number of times that
$\phi^k(e_s)$ crosses $e_r$.  Thus $m^{(k)}_{rs}=(M^k)_{rs}$. The fact
that~$y$ is an eigenvector of $M^k$ with eigenvalue~$\lambda^k$ thus
gives that for any~$k$ and any $r$,
\[y_r=\frac{1}{\lambda^k}\sum_{s=1}^\infty m^{(k)}_{rs}y_s.\]
Now if $q$ is a switch of $\tau$ which is an endpoint of the real edge
$e_{i_0}$, and if~$I$ and~$J$ are the sets of indices of infinitesimal
edges having one and two endpoints at $q$ respectively, then for
any~$s$,
\[\left(\sum_{i\in I}m^{(k)}_{is} +2\sum_{j\in J}m^{(k)}_{js}\right) - 2 \le
m^{(k)}_{i_0s} \le \left(\sum_{i\in I}m^{(k)}_{is} +2\sum_{j\in
J}m^{(k)}_{js}\right) + 2.\] This is because, except at its
endpoints, $\phi^k(e_s)$ must cross both~$e_{i_0}$ and an edge indexed
in~$I\cup J$ each time it intersects~$q$. Multiplying by $y_s$,
summing over~$s$, and dividing by $\lambda^k$ gives
\[\left|\sum_{i\in I}y_i+2\sum_{j\in
J}y_j-y_{i_0}\right|\le\frac{2\sum_{s=1}^\infty y_s}{\lambda^k},\]
which gives the result on letting $k\to\infty$.
\end{proof}

Now let $X\in\R^{N-1}$ be an eigenvector of $A^T$ associated to the
Perron-Frobenius eigenvalue $\lambda$. If $x=(X\,X')\in l^\infty$ is a
completion of~$X$ to an eigenvector of the adjoint~$M^*$ then~$X'=0$,
since it satisfies the equation $(\lambda I-\Pi^*)X'=0$, and
$\|\Pi^*\|_\infty\le 1<\lambda$.

\subsection{The construction}
\label{subsec:construction}
The construction will be illustrated by means of a running example, where
$F\co (S^2,\ofi,A)\barrow$ is a unimodal MIA thick interval map
with $\k(F)=(1001011)^\infty$, as depicted in
Figure~\ref{fig:1001011thtrmap}. The associated generalized train
track map $\phi\co\tau\barrow$ is depicted in
Figure~\ref{fig:1001011}. 

Let $x=(x_1,x_2,\ldots,x_{N-1},0,0,\ldots)$ and $y=(y_1,y_2,\ldots)$
be the eigenvectors of $M^*$ and $M$ associated to the
Perron-Frobenius eigenvalue~$\lambda$ of~$A$ as described above. For
definiteness these vectors are 
scaled so that $\|x\|_2=1$ and $\sum x_iy_i=1$, although this
normalization can be chosen arbitrarily. In the running
example, to~$3$ decimal places:
\[A=\left( \begin{array}{cccccc}
               0&0&0&0&0&1 \\
               0&0&0&0&1&0 \\
               1&0&0&0&1&0 \\
               1&0&0&1&0&0 \\ 
               0&1&0&1&0&0 \\ 
               0&0&1&1&0&0             
             \end{array} \right)\!,\qua\lambda=1.686, \qua
 X=\left(\begin{array}{c}
         0.543\\ 0.104\\ 0.191\\ 0.724\\ 0.175\\ 0.322
         \end{array}\right)\!,\qua
 Y=\left(\begin{array}{c}
         0.368\\0.291\\0.509\\0.536\\0.490\\0.620
         \end{array}\right)
\]
Associate a Euclidean rectangle~$R_i$ of width~$x_i$ and height~$y_i$,
foliated by horizontal and vertical line segments, to each real
edge~$e_i$ of~$\tau$. The idea of the construction is as follows: a
topological sphere~$S$ is constructed by identifying the vertical
sides of these rectangles in the manner dictated by the infinitesimal
edges of~$\tau$ and then identifying the horizontal sides in such a
way that the train track map~$\phi$ defines a homeomorphism
$\Phi\co S\barrow$. The foliations of the rectangles project to
singular foliations of~$S$ with transverse measures induced by
the Euclidean metric: these foliations are preserved by~$\Phi$, the
horizontal and vertical foliations being stretched and contracted by
factors~$\lambda$ and~$1/\lambda$ respectively, so that~$\Phi$ is a
generalized pseudo-Anosov map.

When~$\tau$ is finite (ie, when~$\k(F)=\NBT(m/n)$ for some~$m/n$),
this construction is identical to that of~\cite{BH}, and produces a
pseudo-Anosov homeomorphism. When~$\tau$ is infinite, however, the
construction involves infinitely many sets of identifications, and it
is therefore necessary to describe them carefully and to justify that
the quotient space is a sphere. The tool used to do this is Moore's
theorem~\cite{Moo} on monotone upper semi-continuous (m.u.s.c.)
decompositions of the sphere.

\begin{defns}
A partition or {\em decomposition} $\cG$ of a topological space into
compact sets is {\em upper semi-continuous} if for every decomposition
element $\zeta\in\cG$ and every open set $U\supset\zeta$, there
exists an open set $V\subset U$ with $\zeta\subset V$ such that every
$\zeta'\in\cG$ with $\zeta'\cap V\neq\emptyset$ has $\zeta'\subset
U$. The decomposition is {\em monotone} if its elements are connected.
\end{defns}

The proof of the following lemma is routine.

\begin{lemma}
\label{lem:muscid}
A decomposition~$\cG$ of a compact metric space~$X$ is upper
semi-continuous if and only if whenever $x_n\to x$ and $y_n\to y$ are
convergent sequences in~$X$ such that $x_n$ and $y_n$ belong to the
same element of~$\cG$ for all~$n$, then~$x$ and~$y$ belong to the same
element of~$\cG$.\qed
\end{lemma}

\begin{thm}[Moore]
\label{thm:moore}
Let $\cG$ be a monotone upper semi-continuous decomposition of a
topological sphere $S$ such that no decomposition element
separates~$S$. Then the quotient space obtained by collapsing each
decomposition element to a point is a topological sphere.
\end{thm}

\smallskip
The side identifications will be realized by {\em arc bands} with
which m.u.s.c.\ decompositions of the sphere will be
constructed. These arc bands are of two types:
\begin{itemize}
\item A {\em rectangular arc band~$\eta$ of size~$y$} is the image
of an embedding $\psi\co[0,1]\times[0,y]\to S^2$ which restricts to
a Euclidean isometry on $\{0\}\times[0,y]$ and on $\{1\}\times[0,y]$,
foliated by the leaves $\psi([0,1]\times\{t\})$ for $0\le t\le y$. The
{\em vertices} of~$\eta$ are the points $\psi(0,0)$, $\psi(0,y)$,
$\psi(1,0)$, and $\psi(1,y)$.
\item A {\em semi-circular arc band~$\eta$ of size~$y$} is the image
of an embedding
\[\psi\co\{(s,t)\in\R^2\,:\,t\in[-y,y],\,s\in[0,\sqrt{y^2-t^2}]\}\to
S^2\]
which restricts to a Euclidean isometry on $\{0\}\times[-y,y]$,
foliated by the leaves 
\[\psi\left(\{(\sqrt{r^2-t^2},t)\,:\,t\in[-r,r]\}\right)\]
for $0\le r\le y$. The {\em vertices} of~$\eta$ are the points
$\psi(0,-y)$ and $\psi(0,y)$, and its {\em centre} is $\psi(0,0)$.
\end{itemize}
The decomposition of an arc band~$\eta$ into the leaves of its
foliation is clearly m.u.s.c.

\subsubsection{Identifications of vertical sides}
\begin{figure}[ht!]
\begin{center}
\lab{R1}{R_1}{}\lab{R2}{R_2}{}\lab{R3}{R_3}{}\lab{R4}{R_4}{}\lab{R5}{R_5}{}
\lab{R6}{R_6}{} 
\pichere{0.8}{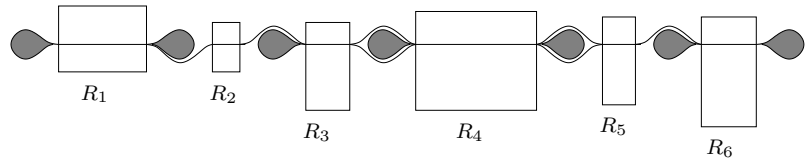}
\caption{The rectangles $R_i$ placed along the real edges of $\tau$}
\label{fig:constr1}
\end{center}
\end{figure}

Place isometric copies of the rectangles~$R_i$ in~$S^2$ along the real
edges of~$\tau$, as illustrated for the running example in
Figure~\ref{fig:constr1}.

\begin{figure}[ht!]
\begin{center}
\lab{ei}{e_i}{}
\lab{ej}{e_j}{}
\lab{ek}{e_k}{}
\lab{etai}{\eta_i}{}\lab{etaj}{\eta_j}{}\lab{etak}{\eta_k}{}
\pichere{0.5}{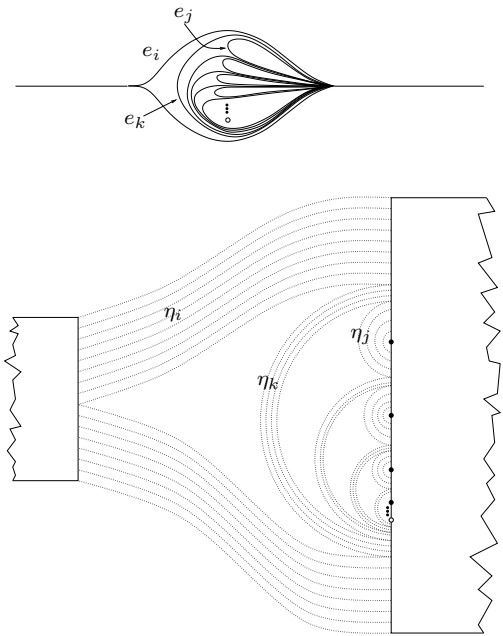}
\caption{The arc bands in the junction between $R_3$ and $R_4$}
\label{fig:constr2}
\end{center}
\end{figure}

For each infinitesimal edge~$e_i$ of~$\tau$, place an arc
band~$\eta_i$ of size~$y_i$ in~$S^2$, joining the segments of
length~$y_i$ on the vertical sides of the rectangles corresponding to
the real edges adjacent to the endpoints of~$e_i$. If~$e_i$ is a
bubble then these segments are adjacent on the same vertical side of
the same rectangle, and~$\eta_i$ is a semi-circular arc band;
otherwise, $\eta_i$ is rectangular. The arc bands are placed on the
sphere in the configuration dictated by~$\tau$, and are chosen to be
mutually disjoint except possibly at their vertices. By
Lemma~\ref{lem:BH}, the arc bands incident on any given vertical side
of a rectangle precisely cover that side. When there are infinitely
many infinitesimal edges in a given junction, then some care is needed
in the placement of these arc bands to ensure that Moore's theorem can
be applied (despite the fact that the arc bands only serve to indicate
how the vertical sides of the rectangles should be identified, and
from this point of view the details of their placement is irrelevant).
\begin{description}
\item[$S^\pm$] Choose the semi-circular arc bands with diameters
converging to~$0$, which is possible since their sizes converge
to~$0$.
\item[$W^\pm$, $V_3^\pm$] Choose the semi-circular arc bands with
diameters converging to~$0$. It is then possible to choose the
rectangular arc bands so that their diameters converge to~$0$ (see
Figure~\ref{fig:constr2} for the arc bands in the junction
between~$R_3$ and~$R_4$ in the running example). Note that each
bounded complementary component of the union of the arc bands in the
junction is a $3$-gon.
\item[$V_1^\pm$, $V_1^\pm\,B$] Choose the semi-circular arc bands with
diameters converging to~$0$. The rectangular arc bands can then be
chosen so that they converge in the Hausdorff topology to a single
arc~$v$ joining the two switches of the junction, and separating all
but one of the arc bands from the puncture. Note that the bounded
complementary components of the union of the arc bands comprise one
$2$-gon (containing the puncture), and infinitely many $3$-gons.
\item[$V_2^\pm$, $V_2^\pm\,B^2$] Choose the semi-circular arc bands
with diameters converging to~$0$, and the rectangular arc bands so
that the sequence of arc bands above (respectively below) the puncture
converges to an arc $v^+$ (respectively~$v^-$) joining the two
switches, and $v^+$ and $v^-$ are disjoint except at their
endpoints. Again, the bounded complementary components of the union of
the arc bands comprise one $2$-gon containing the puncture, and
infinitely many $3$-gons.
\end{description}

Let $R$ denote the closure of $\bigcup_{i<N}R_i\cup\bigcup_{i\ge
N}\eta_i$. The complement of $R$ in $S^2$ consists of one unbounded
component~$C_\infty$ (containing $\infty$), finitely many
open~$2$-gons containing punctures, and countably many open $3$-gons.
Define a decomposition of the sphere with elements of four types:
\begin{enumerate}[i)]
\item the closure of each bounded complementary component of $R$;
\item the connected components of the intersection of the closure of
$C_\infty$ with the closure of the union of the arc bands;
\item the arcs in each arc band except those contained in elements
of~i) or~ii);
\item single points in the rectangles and in~$C_\infty$ which are not
contained in elements of~i),~ii), or~iii).
\end{enumerate}

Decomposition elements of type~ii), which would not otherwise be
decomposition elements of type~iii), arise in the following
circumstances: in $2$-junctions of type $B$, $BP$, $W^\pm$, and
$V_1^\pm\,B$ there is one such element which is the union of two arcs
with a common endpoint; in $2$-junctions of type $V_2^\pm\,B^2$, there
are two such elements, each the union of two arcs with common
endpoints; and in junctions of type $S^\pm$ there is one such element,
which is the closure of the union of infinitely many arcs.

\begin{lemma}
\label{lem:musc1}
This decomposition of $\sph$ is m.u.s.c.\ and none of its elements
separates $\sph$.
\end{lemma}
\begin{proof}
It is clear that the decomposition elements don't separate the
sphere. Thus it is required to show (using Lemma~\ref{lem:muscid})
that if $\zeta_n$ is a sequence of decompositions elements, and $x_n\to
x$, $y_n\to y$ are convergent sequences in $\sph$ with
$x_n,y_n\in\zeta_n$ for all $n$, then $x$ and $y$ lie in the same
decomposition element.

If all but finitely many of the~$\zeta_n$ are single points
then~$x=y$, so passing to subsequences it can be assumed that none of
the~$\zeta_n$ is a single point. Since there are only finitely many
decomposition elements of type~ii), each compact, it can be assumed
that none of the~$\zeta_n$ is of type~ii). If infinitely many of
the~$\zeta_n$ belong to a given arc band then the result follows by
the upper semi-continuity of each arc band, so it can be assumed that
no two of the~$\zeta_n$ belong to the same arc band. If infinitely
many of the~$\zeta_n$ are semi-circular arcs, then~$x=y$ since the
diameters of the semi-circular arc bands tend to~$0$, so it can be
assumed that no~$\zeta_n$ is a semi-circular arc. If infinitely many
of the~$\zeta_n$ are the closure of the same complementary component
of~$R$, then the result follows by the compactness of these
closures. Thus there are only two cases left to consider: where
$(\zeta_n)$ is a sequence of arcs from distinct rectangular arc bands,
and where $(\zeta_n)$ is a sequence of distinct closures of
complementary components of~$R$.

In either case, it can be assumed that the rectangular arc bands or
complementary components are all in a single junction of type
$V_1^\pm$ or $V_2^\pm$, since the diameters of rectangular arc bands
and complementary components in junctions of other types tend
to~$0$. Then by construction,~$x$ and~$y$ lie on the limiting arcs~$v$
or~$v^\pm$ in that junction, which are contained in the closure of the
complementary component of~$R$ containing the periodic point in that
junction.
\end{proof}

\begin{figure}[ht!]
\begin{center}
\pichere{0.9}{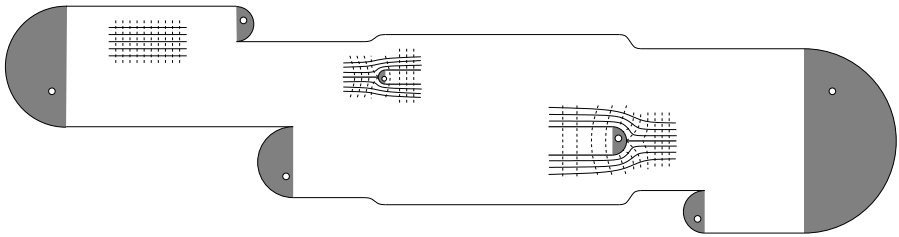}
\caption{The quotient disk $\cR$ and its foliation}
\label{fig:constr3}
\end{center}
\end{figure}

\begin{figure}[ht!]
\begin{center}
\pichere{0.4}{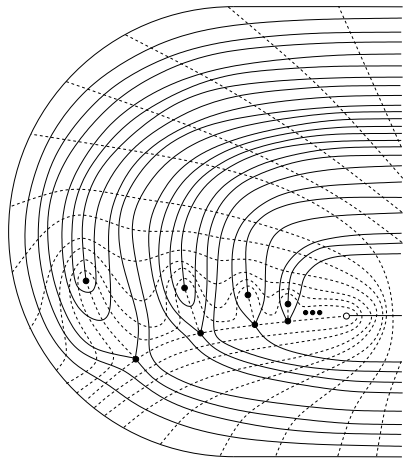}
\caption{Detailed view of the shaded regions of Figure~\ref{fig:constr3}}
\label{fig:bbquot}
\end{center}
\end{figure}

It follows by Theorem~\ref{thm:moore} that the quotient space obtained
by collapsing each decomposition element to a point with the induced
quotient topology is a topological sphere~$\tilde{S}$. The projection
of the complement of~$C_\infty$ is a closed topological disk which is
denoted $\cR$. The horizontal line segments in the rectangles
project to a singular foliation of $\cR$ at all but finitely many
points (the accumulations of singularities). There is a $3$-pronged
singularity at each point of the quotient corresponding to a bounded
complementary $3$-gon, and a $1$-pronged singularity at each point
corresponding to the centre of a semi-circular arc band. The
disk~$\cR$ and parts of its foliation in the case of the running
example is depicted in Figure~\ref{fig:constr3}, and details of the
shaded regions of Figure~\ref{fig:constr3} are shown in
Figure~\ref{fig:bbquot}.


The train track map $\phi\co\tau\barrow$ induces a map
$\tPhi\co\cR\barrow$: unstable leaves (coming from horizontal segments)
are stretched by the factor $\lambda$, stable leaves (coming from
vertical segments) are contracted by the factor $1/\lambda$, the
(projections of) rectangles are mapped as dictated by the action of
$\phi$ on the real edges of $\tau$ and the identifications are mapped
according to the action of $\phi$ on the infinitesimal edges of
$\tau$. 

\medskip
\subsubsection{Identifications of horizontal sides}
\label{sec:hident}
The identifications of the horizontal sides of $\cR$ are carried out
in such a way as to restore the injectivity of $\tPhi$ on
$\partial\cR$. In order to describe them, it is therefore necessary to
understand the structure of $\partial\cR$ and the action of $\tPhi$ on
it. This was the purpose of the analysis of the outside map of the
unimodal map~$f$ in Section~\ref{sec:outside}: the dynamics of
$\tPhi|_{\partial\cR}$ is modelled by the dynamics of the outside map,
in a sense made precise by the following lemma.

\begin{lemma}
\label{lem:outsides}
Let $\theta\co S^1\to S^1\cup I$ be the outside map of~$f$. Then
there is a homeomorphism $h\co S^1\to\partial\cR$ with the property
that $\theta(x)\in S^1$ if and only if $\tPhi(h(x))\in\partial\cR$,
and in this case $h(\theta(x))=\tPhi(h(x))$.
\end{lemma}
\begin{proof}
An analysis of the different possible configurations of infinitesimal
edges in the junctions of~$\ofi$, as given by Theorems~\ref{thm:pertt}
and~\ref{thm:prepertt}, shows that for any kneading sequence
$s\in\MIA$ the two corners on the left of $R_1$ are identified, the
two corners on the right of $R_{N-1}$ are identified, and for $1\le
j<N-1$ the top (respectively bottom) right corner of $R_j$ is
identified with the top (respectively bottom) left corner of
$R_{j+1}$. The corners of the rectangles therefore define $2N-2$
points on $\partial\cR$, which correspond naturally with the $2N-2$
points on $S^1$ defined by the critical orbit of~$f$ (namely $\hat{a}$
and $\hat{b}$, together with~$x_l$ and~$x_u$ for each other point~$x$
in the critical orbit).

Define $h\co S^1\to\partial\cR$ to map each of these $2N-2$ points
on $S^1$ to the corresponding point on $\partial\cR$, and to map each
interval of~$S^1$ between two such points affinely onto the
corresponding interval of $\partial\cR$. Then $h$ is a homeomorphism,
the points $x\in S^1$ with $\theta(x)\not\in S^1$ are precisely those
for which $\tPhi(h(x))\not\in\partial\cR$, and away from these points
both $\theta$ and $\tPhi$ are affine on each of the $2N-2$ intervals,
and send corresponding endpoints to corresponding endpoints.
\end{proof}

\begin{figure}[ht!]
\begin{center}
\lab{eta0}{\eta_0}{}
\lab{eta-1}{\eta_{-1}}{}
\lab{eta-2}{\eta_{-2}}{}
\lab{gamma}{\gamma}{}
\lab{P(G)}{\tPhi(\gamma)}{b}
\lab{P(c)}{\tPhi(c)}{}
\lab{c}{c}{}
\lab{ha}{a}{}
\lab{p}{p}{}
\lab{q}{q=\tPhi(a)=\tPhi(p)}{bl}
\pichere{0.8}{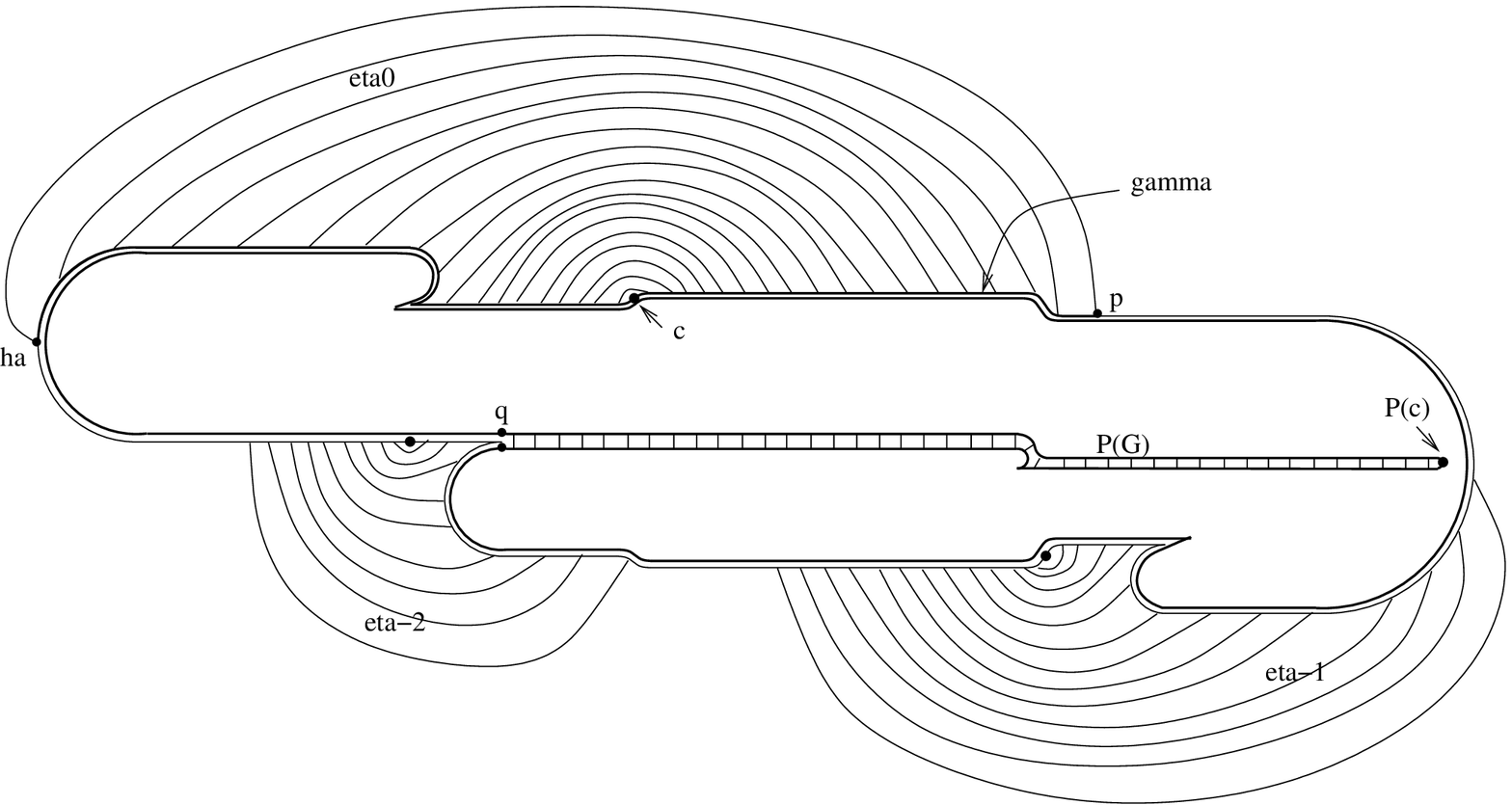}
\caption{$\cR$ and its image under $\tPhi$}
\label{fig:constr4}
\end{center}
\end{figure}

Abusing notation, let $\gamma\subseteq\partial\cR$,
$a,c,p\in\partial\cR$ be defined as $h(\gamma)$, $h(\hat{a})$,
$h(c_u)$, $h(p_u)$ respectively (see Figure~\ref{fig:constr4}). 

\begin{thm}
\label{thm:outside2}
Let $s=\k(F)$, and suppose that $q(s)=m/n\in(0,1/2)\cap\Q$. Then
$\tPhi|_{\partial\cR}$ has a unique periodic orbit, which has
period~$n$ and rotation number~$m/n$. The boundary $\pR$ has the
structure of an $n$-sided polygon, with vertices the points
$\tPhi^j(a)$ for $1\le j\le n$. Moreover, $v=\tPhi^n(a)\in\ol{\gamma}$
satisfies:
\begin{enumerate}[\rm i)]
\item $v=a$ (so that $a$ is periodic under $\tPhi$) if and only
if $s=\lhe(m/n)$
\item $v=c$ if and only if $s=\NBT(m/n)$
\item $v=p$ (so that $p$ is periodic under $\tPhi$) if and only if 
$s=\rhe(m/n)$ 
\item $v$ is strictly between~$a$ and~$c$ if and only if
$\lhe(m/n)\prec s\prec\NBT(m/n)$ 
\item $v$ is strictly between~$c$ and~$p$ if and only if
$\NBT(m/n)\prec s\prec\rhe(m/n)$
\end{enumerate}
\end{thm}
\begin{proof}
The statements about the periodic orbit on $\pR$ and the position
of~$v$ are direct translations of Theorem~\ref{thm:outside}. An
analysis of the different possible configurations of infinitesimal
edges in the junctions of~$\ofi$ shows that there are singularities on
$\pR$ corresponding to each junction with configuration $S^\pm$, and
cusps on $\pR$ corresponding to each $2$-junction with configuration
$BP$, $B$, $W^\pm$, $V_1^\pm B$, or $V_2^\pm B^2$. By
Theorems~\ref{thm:pertt} and~\ref{thm:prepertt}, there are exactly $n$
such junctions, and these junctions are precisely those containing the
first $n$ images under~$F$ of the left hand $1$-junction.
\end{proof}

\begin{remark}
Observe that if $s$ is not one of the endpoints of $\KS(m/n)$, then
the vertices of $\pR$ correspond exactly to the vertices of the
unbounded complementary component of~$\tau$.
\end{remark}

\smallskip
The m.u.s.c.\ decomposition giving the identifications on the
horiztonal edges can now be described: as before, it is constructed
using rectangular and semi-circular arc bands joining points on the
horizontal edges which will be identified. These arc bands will be
labelled $\eta_j$ for $j\le 0$ (in fact $\eta_j$ will often denote a
union of two arc bands), with the arcs of $\eta_j$ joining precisely
those points of $\pR$ which are identified by $\tPhi^{1-j}$ but not by
$\tPhi^{-j}$: thus $\tPhi$ will project to a homeomorphism of the
quotient space. In particular, $\eta_0$ is a semi-circular arc band
centred at~$c$ whose arcs join pairs of points of $\ol{\gamma}$ which
are identified by $\tPhi$. This is illustrated for the running example
in Figure~\ref{fig:constr4}, where the image of $\gamma$ has been
depicted slightly separated for clarity.  Notice that $\tPhi|_{\pR}\I$
is discontinuous at $q=\tPhi(a)=\tPhi(p)$ and is not surjective since
it misses $\gamma$, but away from~$q$ it is a continuous injection
(the injectivity is what makes it possible to construct the arc
bands~$\eta_j$ in such a way that they intersect $\pR$ with disjoint
interiors).

The preimages $\gamma_j=\tPhi^j(\gamma)$ are connected for $-n+1\le j<0$,
since $\ol{\gamma_{-n+1}}$ is the first preimage of $\ol{\gamma}$
which contains the point of discontinuity~$q$ by
Theorem~\ref{thm:outside2}. Thus the arc bands $\eta_j$ for $-n+1\le
j<0$ are semi-circular, centred on $c_j=\tPhi^j(c)$, and can be chosen
to be disjoint.

For $j\le -n$ there are three different cases, according as $s$ is an
endpoint of~$\KS(m/n)$, is equal to~$\NBT(m/n)$, or is neither. As in
the statement of Theorem~\ref{thm:outside2}, let~$v=\tPhi^n(a)$ denote
the vertex of $\pR$ which lies in $\ol{\gamma}$.

\begin{enumerate}[a)] 
\item {\em $s$ is an endpoint of $\KS(m/n)$:} Then $v$ coincides with
an endpoint of $\gamma$ and is periodic of period~$n$. Thus $q$ is one
of the endpoints of $\gamma_{-n+1}$ so the preimages $\gamma_j$ are
connected for all $j<0$. The arc bands $\eta_j$ are semi-circular
centred at $c_j$ for all $j\leq 0$. Since their sizes tend to~$0$,
they can be chosen to be disjoint away from their vertices and to have
diameters tending to~$0$ as $j\to-\infty$. They cover $\pR$ by
Theorem~\ref{thm:outside}, accumulate on the periodic orbit on $\pR$,
and their vertices lie on the orbit of preimages of this periodic
orbit (described in
Theorem~\ref{thm:outside}). Figure~\ref{fig:constr5} depicts these arc
bands in the case $s=\lhe(1/3)$: the white dots depict the points of
the period~$3$ orbit on~$\pR$.
\begin{figure}[ht!]
\begin{center}
\lab{0}{\eta_0}{}
\lab{-1}{\eta_{-1}}{}
\lab{-2}{\eta_{-2}}{}\lab{-3}{\eta_{-3}}{}\lab{-4}{\eta_{-4}}{}
\lab{-5}{\eta_{-5}}{}
\pichere{0.5}{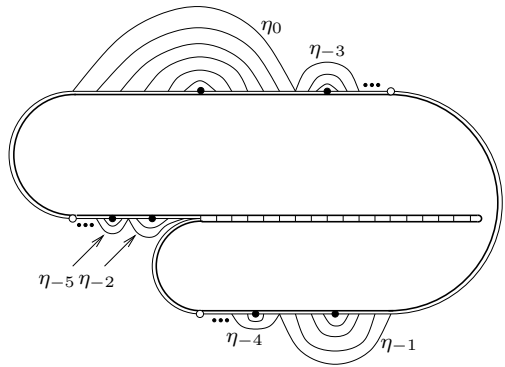}  
\caption{The arc bands in case a): here $s=(101)^\infty=\lhe(1/3)$}
\label{fig:constr5}
\end{center}
\end{figure}  
The m.u.s.c.\ decomposition consists of the following elements:
\begin{enumerate}[i)]
\item the closure of the complement of
$\cR\cup\bigcup_{j=0}^{-\infty}\eta_j$, which is a closed topological
disk bounded by infinitely many arcs, containing
the periodic orbit on its boundary and the point at infinity in its
interior;
\item the arcs of each arc band $\eta_j, j\leq 0$, except for those
contained in elements of~i);
\item the points of $\Int(\cR)$.
\end{enumerate}
\item {\em $s=\NBT(m/n)$:} Then $v=c$, so $q=\tPhi^{-n+1}(c)$ and the
point of discontinuity of $\tPhi|_\pR\I$ coincides with the centre of
$\eta_{-n+1}$. Thus $\ol{\gamma_{-n}}$ is the union of two arcs, and
each point in each arc is paired with exactly one point in the
other. Inductively, the same is true for $\ol{\gamma_j}$ for all
$j\leq -n$: thus the $\eta_{j}$ for $j\le -n$ can be chosen to be
rectangular arc bands, one of whose boundaries coincides with the outer
boundary of $\eta_{j+n}$ and the other with the inner boundary
of~$\eta_{j-n}$. The sizes of these arc bands decrease exponentially,
and they cover $\pR$ by Theorem~\ref{thm:outside}. They are depicted
in Figure~\ref{fig:constr6} for the case $s=\NBT(1/3)$.  The m.u.s.c.\
decomposition consists of the following elements:
\begin{figure}[ht!]
\begin{center}
\lab{c}{c=v}{}
\pichere{0.7}{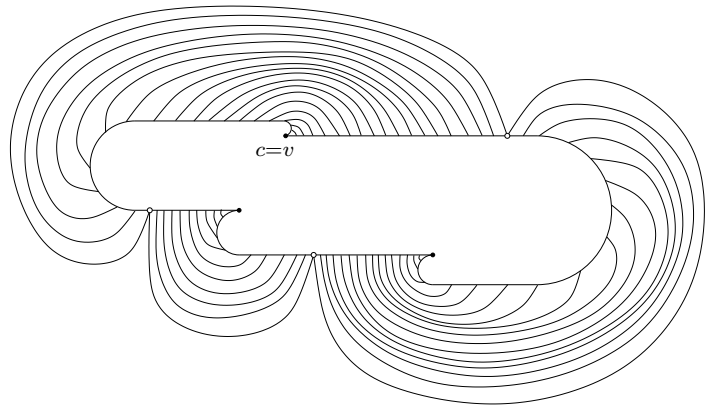}
\caption{The arc bands in the NBT case: here
$s=(10011)^\infty$}
\label{fig:constr6}
\end{center}
\end{figure}
\begin{enumerate}[i)]
\item the closure of the complement of
$\cR\cup\bigcup_{j=0}^{-\infty}\eta_j$, which is a closed topological
disk bounded by $n$ arcs, having the points
of the periodic orbit on $\pR$ for vertices, and containing the point
at infinity in its interior;
\item the arcs of each arc band $\eta_j, j\leq 0$ (in this case, no
such arc is contained in an element of~i));
\item the points of $\Int(\cR)$.
\end{enumerate}
\item {\em $s$ is neither $\NBT(m/n)$ nor an endpoint of $\KS(m/n)$:}
As in case~b), $\gamma_j$ has two connected components for $j\le
-n$. However, because $v\neq c$, one of the components of each such
$\gamma_j$ is divided into two subarcs: the points of one subarc are
paired with other points in the same subarc, while those in the other
are paired with points in the other component of $\gamma_j$.  In
Figure~\ref{fig:constr7},
$\gamma_{-n+1}=\alpha_1\cup\alpha_2\cup\alpha_3$, points of
$\tPhi\I(\alpha_1)$ are paired to points of $\tPhi\I(\alpha_3)$
whereas points of $\tPhi\I(\alpha_2)$ are paired to points of
$\tPhi\I(\alpha_2)$ itself. Thus $\eta_j$ is the union of a
semi-circular arc band and a rectangular arc band; the diameters of
the semi-circular bands tend to~$0$ as $j\to -\infty$, as do the sizes
of the rectangular bands. The arc bands cover $\pR$ by
Theorem~\ref{thm:outside}.

\begin{figure}[ht!]
\begin{center}
\lab{P}{\tPhi\I}{}
\lab{c}{c}{}\lab{v}{v}{}\lab{q}{q}{}
\lab{a}{\alpha_1}{}\lab{b}{\alpha_2}{}
\lab{k}{\alpha_3}{}
\lab{d}{\tPhi\I(\alpha_1)}{}\lab{e}{\tPhi\I(\alpha_2)}{}
\lab{f}{\tPhi\I(\alpha_3)}{}
\lab{0}{\eta_0}{}\lab{n}{\eta_{-n}}{}
\lab{c(n-1)}{c_{n-1}}{}
\pichere{0.8}{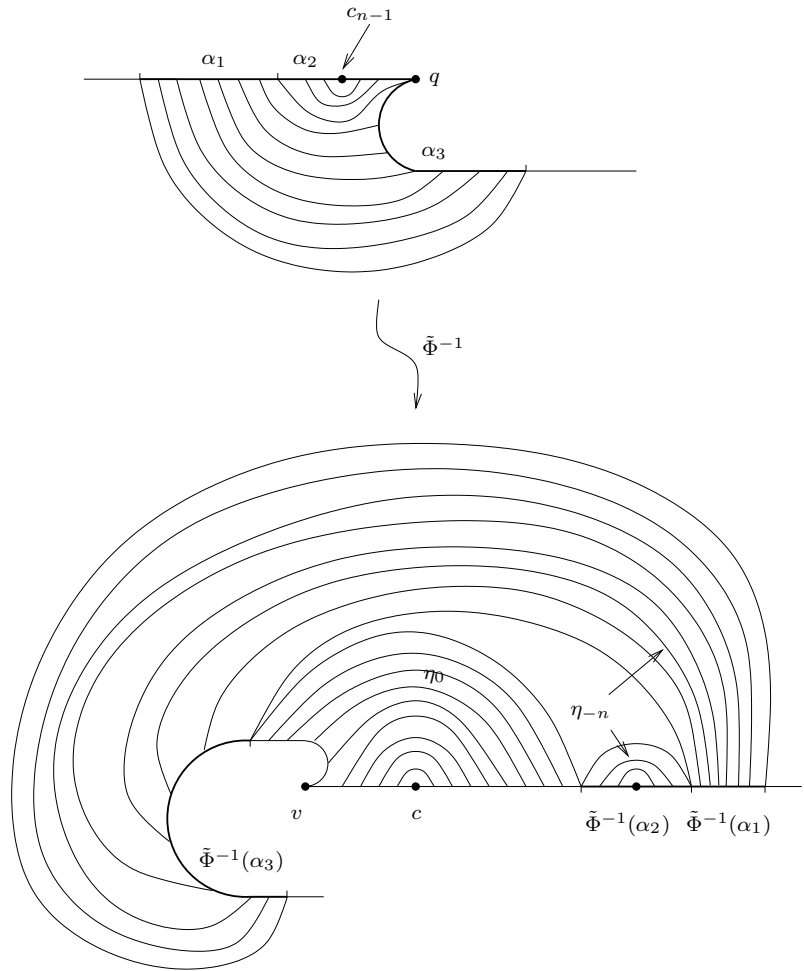}
\caption{The arc bands in case~c)}
\label{fig:constr7}
\end{center}
\end{figure}
The m.u.s.c.\ decomposition consists of the following elements:
\begin{enumerate}[i)]
\item the closures of the (infinitely many) complementary components
of $\cR\cup\bigcup_{j\leq 0}\eta_j$; there is one exterior $n$-gon
containing $\infty$ and having the points in the periodic orbit on
$\pR$ for vertices, and infinitely many 3-gons;
\item the arcs of each arc band in $\eta_j, j\leq 0$, except for those
contained in elements of~i);
\item the points of $\Int(\cR)$.
\end{enumerate}
The complete decomposition for the running example, which satisfies
$\lhe(1/3)\prec s\prec\NBT(1/3)$, is shown in
Figure~\ref{fig:constr8}.
\begin{figure}[ht!]
\begin{center}
\pichere{0.7}{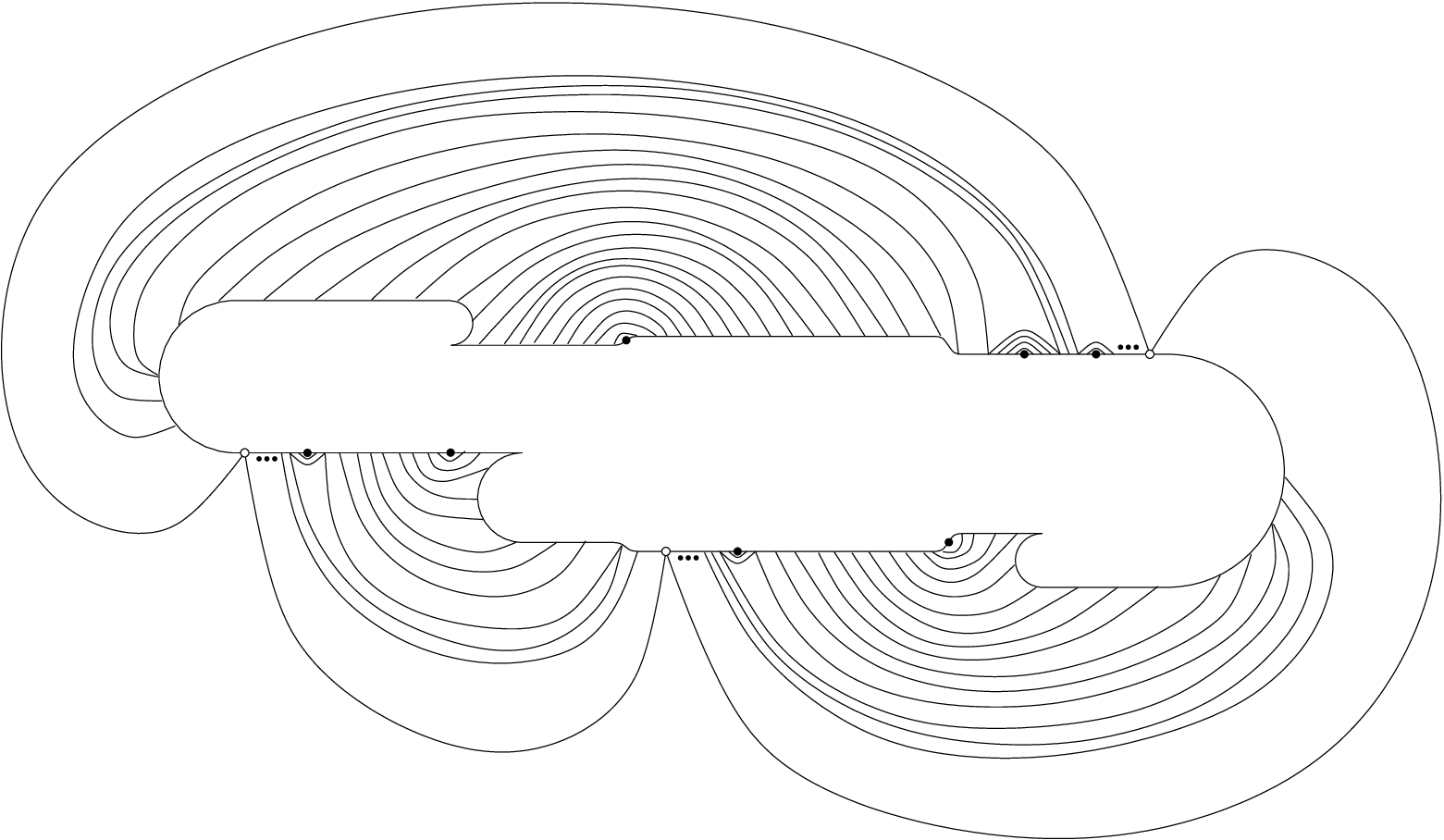}
\caption{The complete second stage decomposition for the running
example} 
\label{fig:constr8}
\end{center}
\end{figure}
\end{enumerate}

\begin{lemma}
\label{lem:musc2}
Each of these decompositions is m.u.s.c.\ and none of its elements
separates $\tilde{S}$.
\end{lemma}
\begin{proof}
 In case~a), where the diameters of the arc bands tends to~$0$, and
case~b), where the arc bands can be regarded as a finite number of
semi-circular arc bands, the proof is straightforward. For case~c),
let $\zeta_n$ be a sequence of decomposition elements, and $x_n\to x$,
$y_n\to y$ be convergent sequences with $x_n,y_n\in\zeta_n$ for
all~$n$. Then it is required to show that $x$ and $y$ lie in the same
decomposition element. As in the proof of Lemma~\ref{lem:musc1}, it
can be assumed that either every $\zeta_n$ is an arc of a distinct
rectangular arc band, or every $\zeta_n$ is a distinct complementary
component. In both cases, both~$x$ and~$y$ lie on the boundary of the
complementary component containing~$\infty$.
\end{proof}

By Theorem~\ref{thm:moore} the quotient space is a topological
sphere~$S$, and by construction $\tPhi$ projects to a homeomorphism
$\Phi\co S\barrow$.  The following theorem asserts that $\Phi$ is a
generalized pseudo-Anosov map, and describes the orbits of 1-pronged
singularities, which are the most important from the point of view of
dynamics. Recall that an orbit $\{f^n(x)\,:\,n\in\Z\}$ of a
homeomorphism $f\co X\barrow$ is {\em homoclinic} to a fixed point~$p$
of~$f$ if $f^n(x)\to p$ as $n\to\pm\infty$.

\begin{thm}
\label{thm:1-prongs}
Let $s\in\MIA$ have height $q(s)=m/n\in(0,1/2)$ and let $\Phi_s\co
S\barrow$ denote the associated homeomorphism constructed above. Then
$\Phi_s$ is a generalized pseudo-Anosov map. Moreover, it has exactly
one orbit of 1-pronged singularities which
\begin{enumerate}[\rm i)]
\item is homoclinic to the point at infinity if and only if $s$ is an
endpoint of the interval $\KS(m/n)$;
\item is finite if and only if $s=\NBT(m/n)$;
\item is backward asymptotic to the point at infinity and
forward asymptotic to the periodic orbit determined by~$s$ if and only
if $\lhe(m/n)\prec s\prec\NBT(m/n)$ or $\NBT(m/n)\prec s\prec\rhe(m/n)$.
\end{enumerate}
\end{thm}

\begin{proof}
Let $\pi\co S^2\to S$ denote the projection map given by the
construction, and let $Q=\bigcup_{i=1}^{N-1}R_i$. Then $\pi(Q)=S$,
since every fibre of~$\pi$ intersects~$Q$; and $\pi|_{\Int(Q)}$ is an
embedding, since every fibre which intersects~$\Int(Q)$ is
trivial. Thus the foliations of~$Q$ by horizontal and vertical line
segments project to regular foliations of $\pi(\Int(Q))$, which are
equipped with transverse measures induced by Euclidean distance
in~$Q$. The arc bands describe isometric identifications of
segments of $\partial Q$, with vertical (respectively horizontal)
segments always identified with other vertical (respectively
horizontal) segments. Thus the measured foliations are also regular at
points of $\pi(\partial Q)$ other than those which are the projections
of vertices of arc bands, or centres of semi-circular arc bands. The
foliations have a $1$-pronged singularity at $\pi(c)$ whenever $c$
is the centre of a semi-circular arc band. If $v$ is the vertex of an
arc band, there are more possibilities.
\begin{enumerate}[a)]
\item If~$v$ is the vertex of an arc band on the interior of 
a vertical segment of $\partial Q$ corresponding to a junction of any
type other than $S^\pm$, then $\pi(v)$ is a $3$-pronged singularity.
\item If~$v$ is the vertex of an arc band on a vertical segment of
$\partial Q$ corresponding to a junction of type $S^\pm$, then
$\pi(v)$ is an essential ($\infty$-pronged) singularity.
\item If~$v$ is a vertex of an arc band on a horizontal segment in the
case where~$s$ is an endpoint of $\KS(m/n)$, then $\pi(v)$ is an
essential singularity.
\item If~$v$ is a vertex of an arc band on a horizontal segment in the
case where $s=\NBT(m/n)$, then~$\pi(v)$ is an $n$-pronged singularity
if~$v$ is a point of the periodic orbit on $\partial Q$, and is a
regular point otherwise.
\item If~$v$ is a vertex of an arc band on a horizontal segment in the
case where $s\neq\NBT(m/n)$ and~$s$ is not an endpoint of~$\KS(m/n)$,
then $\pi(v)$ is an essential singularity if~$v$ is a point of the
periodic orbit on~$\partial Q$, and is a $3$-pronged singularity
otherwise.
\end{enumerate}
By considering the singularities corresponding to junctions of each
given type, it can easily be seen that junctions of types~$BP$,~$B$,
and~$V_0$ give rise to only finitely many singularities, whereas all
other types give rise to infinitely many singularities which
accumulate on the periodic point in the associated junction. If
$s\neq\NBT(m/n)$ then the singularities arising from arc bands on
horizontal segments of~$\partial Q$ accumulate on the point $\pi(P)$,
where~$P$ is the periodic orbit on~$\partial Q$; while if
$s=\NBT(m/n)$ then there are only finitely many such singularities.

Thus the singularities of the foliations of~$S$ accumulate at only
finitely many points. By construction, these foliations are invariant
under~$\Phi$, with the horizontal (unstable) foliation being expanded
by a factor~$\lambda>1$ and the vertical (stable) foliation being
contracted by a factor~$1/\lambda$. Hence~$\Phi$ is a generalized
pseudo-Anosov map as required.

To show that the one-pronged singularities arising from arc bands on
the vertical sides of~$\partial Q$ lie on a single orbit, it is
necessary to consider the action of the train track map
$\phi\co\tau\barrow$ on the bubbles of~$\tau$. If $s=\NBT(m/n)$,
then it is clear that the~$n+2$ bubbles are permuted cyclically
by~$\phi$. In all other cases where~$s$ is periodic, the proof of
Theorem~\ref{thm:pertt} shows that there is a bubble~$B$ in
junction~$N$ which is the image of an infinitesimal edge in
junction~$\pi\I(N)$ joining the two switches of the junction and
passing above the puncture, and all other bubbles are images of this
one. Hence the forwards $\Phi$-orbit of the $1$-pronged singularity
corresponding to~$B$ contains precisely all of the $1$-pronged
singularities arising from arc bands on the vertical sides
of~$\partial Q$. Similarly, if~$s$ is preperiodic then all of the
bubbles of~$\tau$ are images of the unique bubble~$B$ of~$\tau$ in
junction~$N$.

It is immediate from the construction of Section~\ref{sec:hident} that
the one-pronged singularities arising from arc bands on the horizontal
sides of~$\partial Q$ are precisely the points of the backwards
$\Phi$-orbit of the point~$\pi(c)$, since the centres of the arc bands
on $\pR$ are at the points $\tPhi^{-j}(c)$ for $j\ge0$. Since
$\tPhi(c)$ lies on the right hand edge of~$Q$, it follows
that~$\Phi(\pi(c))$ is the unique one-pronged singularity
corresponding to a bubble in junction~$N$ whose preimage has not
already been assigned, and hence the $1$-pronged singularities form a
single orbit as required. That the asymptotics of this orbit are as
given in the statement of the theorem follows directly from the facts
about $1$-pronged singularities established earlier in the proof.
\end{proof}

\begin{remark}
$s=10^\infty$ is the only element of $\MIA$ with height~$0$, and was
excluded from the construction above in order not to have to introduce
exceptions into each statement about $\KS(m/n)$. It can, however, be
treated in exactly the same way as other elements of~$\MIA$: the
invariant generalized train track, depicted in
Figure~\ref{fig:hstrtr}, is described by
\[(S^+\,\,;\,\,B)\]
(and hence is unique among train tracks corresponding to strictly
preperiodic kneading sequences in having a junction of type $S^+$),
and the same construction can be applied starting with this train
track to obtain a `tight horseshoe'. This has been described
in~\cite{dC}.
\end{remark}

\section{The complex structure}
\label{sec:comp}

In this section it is shown that the topological sphere~$S$
constructed above carries a natural complex structure which makes it
into a complex sphere. With respect to this structure the
$\Phi$-invariant foliations become the horizontal and vertical
trajectories of an integrable quadratic differential which is
meromorphic away from finitely many essential singularities, with
respect to which~$\Phi$ is a Teichm\"uller mapping. This is a special
case of results sketched in~\cite{dC}.

\begin{thm}
\label{thm:complexsphere}
Let $s\in\MIA$, and let $\Phi_s\co S\barrow$ be the associated
generalized pseudo-Anosov map. Then the Euclidean structure on the
rectangles used in the construction of~$S$ induces a complex structure
on~$S$, with respect to which it is a complex sphere.
\end{thm}

\begin{proof}
As in the proof of Theorem~\ref{thm:1-prongs}, let $\pi\co S^2\to S$
be the projection map given by the construction of~$S$ and
$Q=\bigcup_{i=1}^{N-1}R_i$, where the~$R_i$ are the rectangles
corresponding to the real edges of the invariant train
track~$\tau$. Since $\pi|_{\Int(Q)}$ is a homeomorphism onto its
image, the Euclidean structure on $\Int(Q)$ induces a natural complex
structure on $\pi(\Int(Q))$. Moreover, the conformal mappings
$z\mapsto z^{2/k}$ (on the slit plane) induce complex structures on
the neighbourhoods of each $k$-pronged singularity. Thus~$S$ has a
complex structure away from the finite set~$\Sigma$ of accumulations
of singularities of the invariant foliations. To establish the
theorem, it is necessary to show that the complex structure on
$S\setminus\Sigma$ regards each point of~$\Sigma$ as a puncture rather
than a hole. Using the theory of extremal length~\cite{Ah,LV}, this
can be accomplished by constructing a nested sequence of annuli
converging to each point of~$\Sigma$, the sum of whose moduli is
divergent. The standard estimate
\[\Mod(A)\ge\frac{\Width(A)^2}{\Area(A)}\] for the modulus of an
annular region~$A$ will be used, where $\Width(A)$ is the minimum
distance between the two boundary components of~$A$.

Let $q(s)=m/n$. If $s=\NBT(m/n)$ then~$\Sigma$ is empty, and there is
nothing to show. If~$\lhe(m/n)\prec s\prec\NBT(m/n)$
or~$\NBT(m/n)\prec s\prec\rhe(m/n)$ then there is one point
of~$\Sigma$ corresponding to each junction of type~$W$, $V_1$, $V_2$,
or~$V_3$ (ie, each junction containing infinitely many infinitesimal
edges), and one point of~$\Sigma$ corresponding to the point at
infinity. Finally, if~$s=\lhe(m/n)$ or~$s=\rhe(m/n)$ then~$\Sigma$
contains a single point. This case is harder, since the
annuli cannot be constructed locally in~$Q$, and is treated first.

\smallskip\noindent
{\em Case 1. $s=\lhe(m/n)$ or $s=\rhe(m/n)$}

Suppose that $s=\lhe(m/n)$: the case $s=\rhe(m/n)$ is entirely
analogous. Let~$M$ be the transition matrix for the train track map
$\phi\co\tau\barrow$,~$\lambda$ be its Perron-Frobenius eigenvalue,
and~$y$ be its Perron-Frobenius eigenvector. A construction similar to
that described below can be found in~\cite{GaEa}.

By Theorem~\ref{thm:pertt},~$\tau$ has~$n-2$ infinitesimal edges which
are not loops, one in each junction other than the leftmost and
rightmost. Carrying out the identifications corresponding to these
edges yields a Euclidean polygon~$E$ with~$n$ vertical and~$n$
horizontal sides. Infinitely many other identifications are applied to
produce~$S$, each of which involves identifying the two halves of an
{\em identification interval} on the boundary of~$E$. The
identification intervals on the vertical sides of~$E$ are given by
(the eigenvector entries corresponding to) the bubbles of~$\tau$,
while those on the horizontal sides are given by the map~$\tilde\Phi$
and its backwards iterates, as described in the proof of
Theorem~\ref{thm:outside2}.

The following notation will be used:
\begin{itemize}
\item $w_v$ is half the length of the longest identification interval
on the vertical sides of~$E$: in other words,~$w_v$ is the largest
entry of~$y$ which corresponds to a bubble of~$\tau$.
\item $w_h$ is half the length of the longest identification interval
on the horizontal sides of~$E$ (this interval is denoted~$\gamma$ in
the proof of Theorem~\ref{thm:outside2}).
\item $w=\min(w_v,w_h)$, $W=\max(w_v,w_h)$.
\item $a^0,a^1,\ldots,a^{n-1}$ are the points of the periodic orbit on
the boundary of~$E$, labelled so that $a^0$ is the rightmost and
$a^{i+1\bmod{n}}$ is the image of $a^i$. Thus $a^1$ is the leftmost
point of the orbit, denoted~$a$ in the proof of
Theorem~\ref{thm:outside2}. Observe that each~$a^i$ is a vertex
of~$E$.
\item Each side of~$E$ is the union of a sequence of adjacent
identification intervals which converge to one of the points
$a^i$. Denote by $(u^i_j)_{j=0}^\infty$ (respectively
$(v^i_j)_{j=0}^\infty$) this sequence on the vertical (respectively
horizontal) side of~$E$ with endpoint~$a^i$. 
\item $a^i_j$ and $a^i_{j+1}$ (respectively~$b^i_j$ and~$b^i_{j+1}$)
are the endpoints of the identification interval~$u^i_j$
(respectively~$v^i_j$).  Thus
\begin{enumerate}[i)]
\item The vertical (respectively horizontal) side of~$E$ containing
the vertex $a^i$ has~$a^i_0$ (respectively~$b^i_0$) as its other
vertex.
\item $(a^i_j)$ and $(b^i_j)$ are both sequences
converging to $a^i$. All of the points of
$\bigcup_{i,j}\{a^i_j,b^i_j\}\cup \bigcup_i\{a^i\}$ are identified to
the point at infinity in~$S$.
\item $|u_0^0|=2w_v$ (where $|I|$ denotes the length of an
interval~$I$), and $|v_0^{i_0}|=2w_h$ for some~$i_0$ (with
the property that $a^{i_0}$ is the topmost point of the periodic orbit
on the boundary of~$E$).
\item \begin{eqnarray*}
|u^i_j|&=&\frac{2w_v}{\lambda^{i+nj}}\qquad\text{ and}\\
|v^i_j|&=&\frac{2w_h}{\lambda^{(i-i_0\bmod{n})+nj}}
\end{eqnarray*}
for $i=0,1,\ldots, n-1$ and $j\ge0$.
\end{enumerate}
\item For each $j\ge 0$, let
\[
r_j=\frac{w}{\lambda^{(n-1)+nj}}. \]
That is, $r_j$ is half the length of the smallest of the~$2n$ identification
intervals $u^i_j$ and $v^i_j$. 
\end{itemize}
A sequence of nested annuli $X_k$, converging to the point at
infinity, will now be constructed. Each~$X_k$ will be given as the
union of finitely many regions of three types.
\begin{description}
\item[Half-annuli] For $0\le i<n$ and $1\le j\le k$, let $A^i_{j,k}$
(respectively $B^i_{j,k}$) be the half-annulus centred at $a_j^i$
(respectively $b_j^i$) with inner radius~$r_{k+1}$ and outer
radius~$r_k$. (These are half-annuli bounded by Euclidean semicircles
and segments of the boundary of~$E$: see
Figure~\ref{fig:halfannulus}.)

\begin{figure}[ht!]
\begin{center}
\lab{a1}{a^i_{j-1}}{}\lab{a2}{a^i_{j}}{}\lab{a3}{a^i_{j+1}}{}
\lab{A1}{A^i_{j-1,k}}{l}\lab{A2}{A^i_{j,k}}{l}\lab{A3}{A^i_{j+1,k}}{l}
\lab{u1}{u^i_{j-1}}{}\lab{u2}{u^i_j}{}
\lab{rk}{r_k}{l}\lab{rk+1}{r_{k+1}}{l}
\pichere{0.3}{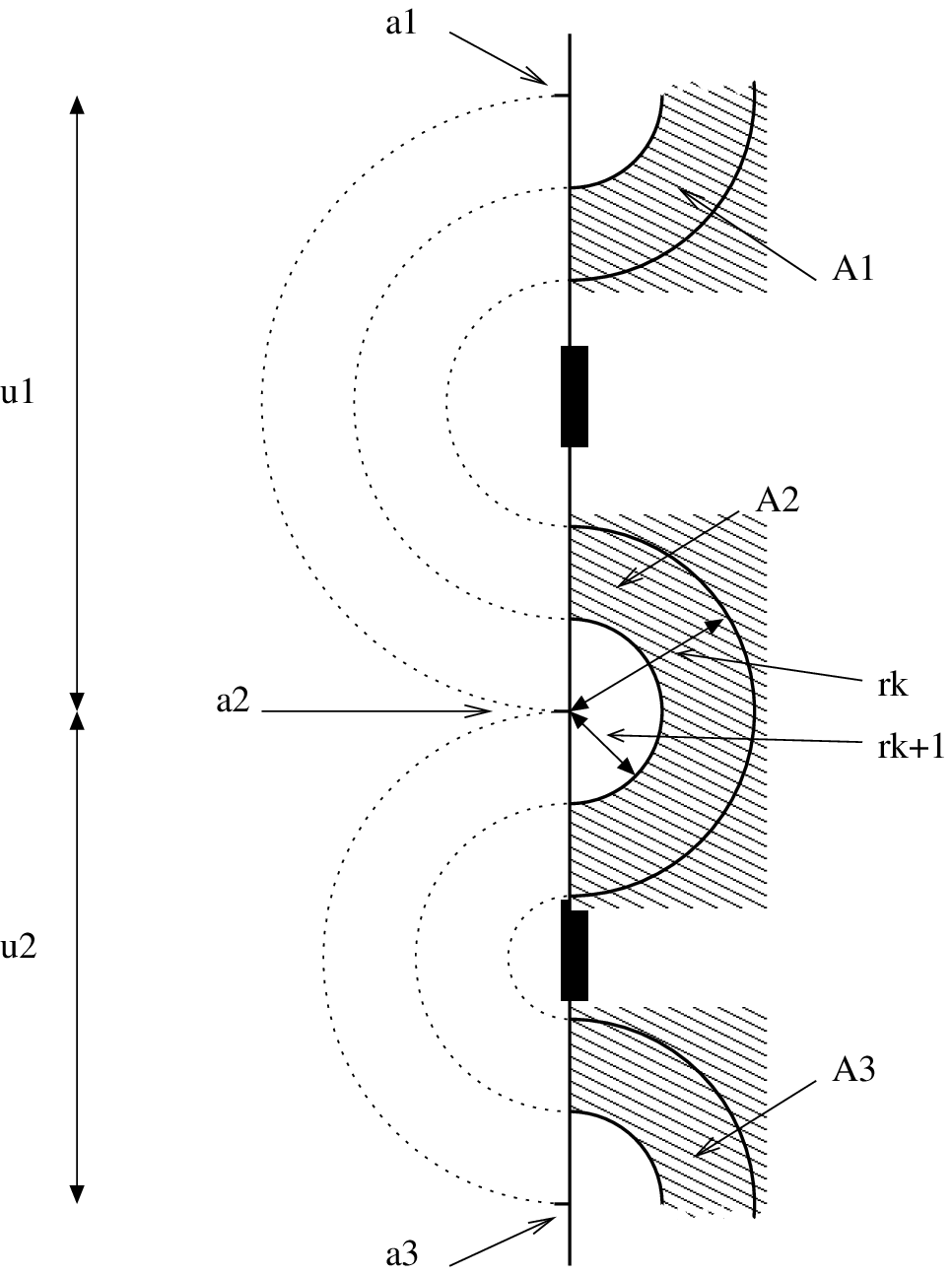}
\caption{The half-annuli $A^i_{j,k}$}
\label{fig:halfannulus}
\end{center}
\end{figure}

The identification intervals which share $a^i_j$ as an endpoint
are~$u^i_{j-1}$ and~$u^i_j$ and $r_k\le|u^i_j|/2 <
|u^i_{j-1}|/2$, so~$A^i_{j,k}$ meets the identification
intervals~$u^i_{j-1}$ and~$u^i_j$ between $a^i_j$ and the centres of
the intervals, as in Figure~\ref{fig:halfannulus}. The same argument
applies to the half-annuli~$B^i_{j,k}$. 

After the identifications on the boundary of~$E$ have been
carried out, each $\bigcup_{j=1}^k A^i_{j,k}$ and each $\bigcup_{j=1}^k
B^i_{j,k}$ is a strip of width $r_k-r_{k+1}$ and area
$k\pi(r_k^2-r_{k+1}^2)/2$.

\item[Rectangular quarter-annuli] Figure~\ref{fig:quarterannulus}
shows the half-annuli $A^i_{k,k}$ and $B^i_{k,k}$. The
endpoints in $u_k^i$ of the boundary semicircles of $A^i_{k,k}$ lie in
the half of $u_k^i$ which is further from $a^i$, and these endpoints
are identified with two points in the half of $u_k^i$ nearer to $a^i$,
which are denoted $z_1$ and $z_2$ in the figure. There are analogous
points $z_3$ and $z_4$ identified with endpoints of the half-annulus
$B^i_{k,k}$. Since the internal angle of~$E$ at~$a_i$ is $\pi/2$, it
is possible to define a rectangular quarter-annulus $Z^i_k$ as shown
in the figure, each of whose boundary components is the union of a
horizontal and a vertical arc, joining either~$z_1$ to~$z_3$ or~$z_2$
to~$z_4$. 

\begin{figure}[ht!]
\begin{center}
\lab{a1}{a^i_{k-1}}{r}\lab{a2}{a^i_{k}}{r}
\lab{A1}{A^i_{k-1,k}}{l}\lab{A2}{A^i_{k,k}}{l}
\lab{u1}{u^i_{k-1}}{r}\lab{u2}{u^i_k}{r}
\lab{b1}{b^i_{k-1}}{}\lab{b2}{b^i_{k}}{}
\lab{B1}{B^i_{k-1,k}}{b}\lab{B2}{B^i_{k,k}}{b}
\lab{v1}{v^i_{k-1}}{t}\lab{v2}{v^i_k}{t}
\lab{z1}{z_1}{bl}\lab{z2}{z_2}{bl}\lab{z3}{z_3}{bl}\lab{z4}{z_4}{bl}
\lab{ai}{a^i}{}\lab{Z}{Z^i_k}{bl}
\pichere{0.5}{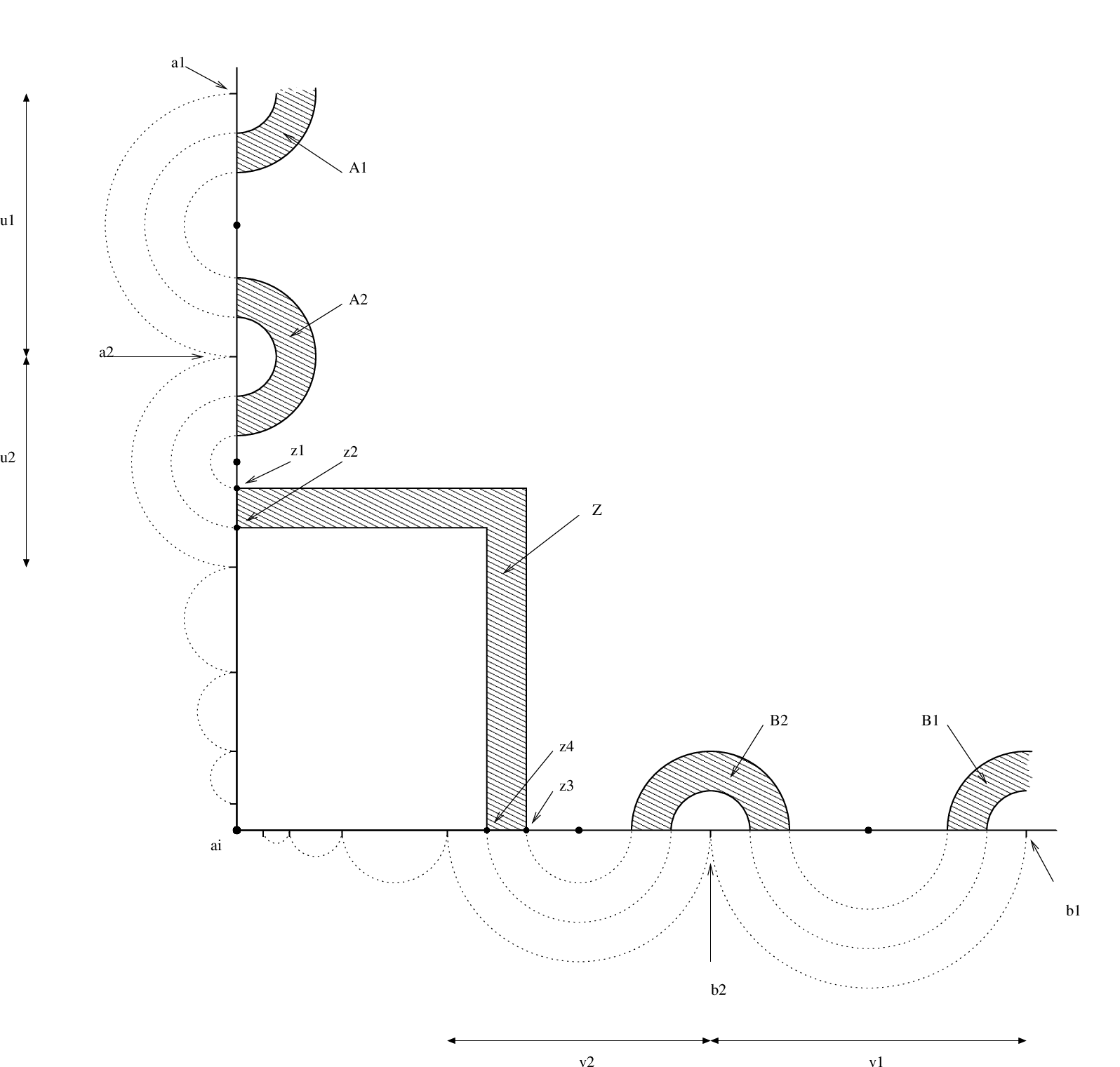}
\caption{The rectangular quarter-annulus $Z^i_k$}
\label{fig:quarterannulus}
\end{center}
\end{figure}

A rectangular quarter-annulus is used since the distance
between~$a^i$ and~$z_1$ need not be equal to the distance
between~$a^i$ and~$z_3$. The former distance is less than
$\sum_{j=k}^\infty |u^i_j|$ and the latter is less than
$\sum_{j=k}^\infty |v^i_j|$. Hence, for each~$i$ and
each~$k$, the region
\[Y^i_k=Z^i_k\cup\bigcup_{j=1}^k A^i_{j,k}\cup \bigcup_{j=1}^k B^i_{j,k}\]
is (once the identifications have been carried out) a strip of width
$r_k-r_{k+1}$ and area bounded above by
\[k\pi(r_k^2-r_{k+1}^2)+(r_k-r_{k+1})\left( \sum_{j=k}^\infty
|u^i_j|+\sum_{j=k}^\infty |v^i_j|\right).\]

\item[Circular vertex annuli] The other vertices $a^i_0=b^{i'}_0$
of~$E$ have internal angle either $\pi/2$ or $3\pi/2$. At each such
vertex, construct a (circular) quarter-annulus or
three-quarter-annulus $D^i_k$ with radii $r_k$ and $r_{k+1}$, whose
boundary components join points identified with endpoints of
$A^i_{1,k}$ to points identified with endpoints of
$B^{i'}_{1,k}$. These partial annuli have areas bounded above by
$\pi(r_k^2-r_{k+1}^2)$.
\end{description}

After the identifications have been carried out,
\[X_k=\bigcup_{i=0}^{n-1} (Y^i_k\cup D^i_k)\]
is an annulus of width $r_k-r_{k+1}$ and area bounded above by
\[n\left((k+1)\pi(r_k^2-r_{k+1}^2)+(r_k-r_{k+1})\left( \sum_{j=k}^\infty
|u^i_j|+\sum_{j=k}^\infty
|v^i_j|\right)\right).\]
The annuli $X_1$ and $X_2$ in the case $s=\lhe(1/3)=(101)^\infty$ are
depicted in Figure~\ref{fig:101infty}.

\begin{figure}[ht!]
\begin{center}
\lab{A212}{A^2_{1,2}}{lb}\lab{A211}{A^2_{1,1}}{}\lab{A222}{A^2_{2,2}}{}
\lab{D22}{D^2_2}{}
\lab{D21}{D^2_1}{}\lab{Z22}{Z_2^2}{}\lab{Z21}{Z^2_1}{}\lab{B222}{B^2_{2,2}}{}
\lab{B212}{B^2_{1,2}}{}\lab{B211}{B^2_{1,1}}{}\lab{B111}{B^1_{1,1}}{bl}
\lab{D01}{D^0_1}{}\lab{D02}{D^0_2}{}
\lab{a0}{a_0}{bl}\lab{a1}{a_1}{}\lab{a2}{a_2}{}
\lab{u00}{u^0_0}{l}
\lab{u01}{u^0_1}{l}\lab{u10}{u^1_0}{r}\lab{u11}{u^1_1}{r}
\lab{v00}{v^0_0}{b}\lab{v01}{v^0_1}{b}\lab{v20}{v^2_0}{}\lab{v21}{v^2_1}{}
\pichere{.94}{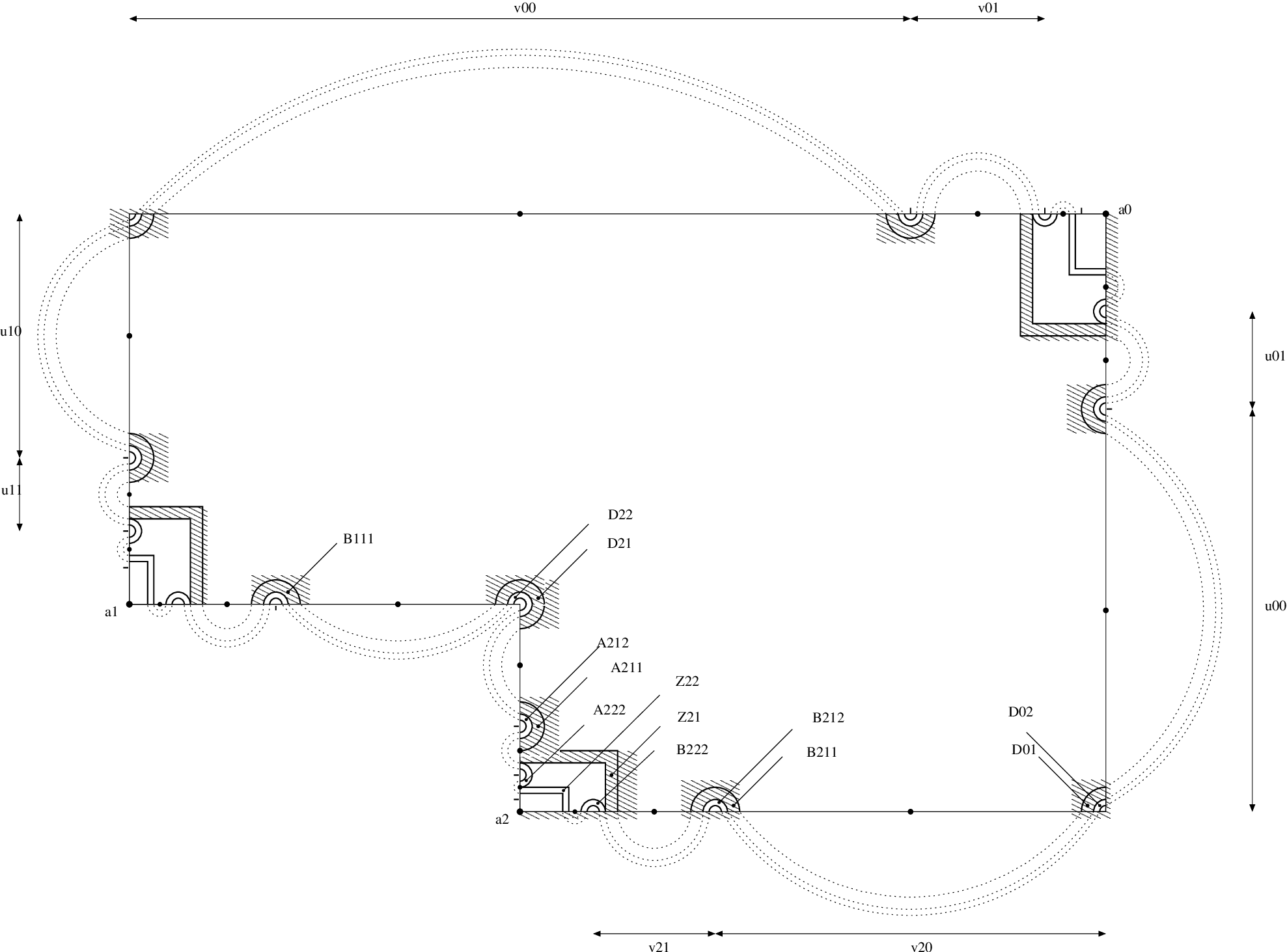}
\end{center}
\caption{Annuli with divergent moduli sum surrounding the 
accumulation of singularities of $\Phi_{\lhe(1/3)}$}
\label{fig:101infty}
\end{figure}

Now 
\[r_k\pm r_{k+1}=\frac{w}{\lambda^{(n-1)+nk}}(1\pm\frac{1}{\lambda^n})\]
and
\[\sum_{j=k}^\infty
|u^i_j|+\sum_{j=k}^\infty |v^i_j| \le
\sum_{j=k}^\infty
\frac{4W}{\lambda^{nj}}=\frac{4W\lambda^n}{\lambda^{nk}(\lambda^n-1)},\] 
so
\begin{eqnarray*}
\Mod(X_k)\ge\frac{\Width(X_k)^2}{\Area(X_k)}&\ge&
\frac{r_k-r_{k+1}}{n\left((k+1)\pi(r_k+r_{k+1})+
\frac{4W\lambda^n}{\lambda^{kn}(\lambda^n-1)}\right)}
\\
&=& \frac{C_1/\lambda^{nk}}{(C_2k+C_3)/\lambda^{nk}}\\
&=&\frac{C_1}{C_2k+C_3},
\end{eqnarray*}
where $C_1$, $C_2$, and $C_3$ depend only on $n$, $\lambda$, $w$, and
$W$.
Hence $\sum_{k\ge1}\Mod(X_k)$ diverges as required.

\smallskip\noindent {\em Case 2. $\lhe(m/n)\prec s\prec\NBT(m/n)$ or
$\NBT(m/n)\prec s\prec\rhe(m/n)$}

Consider first the accumulation of singularities in a junction of
type~$W^\pm$. In Figure~\ref{fig:Wpuncture} the identifications on~$Q$
are indicated with dotted lines, and the shaded half-annuli $A_k$
project to nested annular regions in~$S$ which converge to the
accumulation of singularities~$p$. If the width of the rectangle in
which these half-annuli lie is less than half its height then not all
of them can be constructed there: in this case, simply start with the
first half-annulus $A_{k_0}$ which can be constructed.

\begin{figure}[ht!]
\begin{center}
\lab{c}{c}{}
\lab{A1}{A_0}{}
\lab{A2}{A_1}{}
\lab{A3}{A_2}{}
\pichere{0.5}{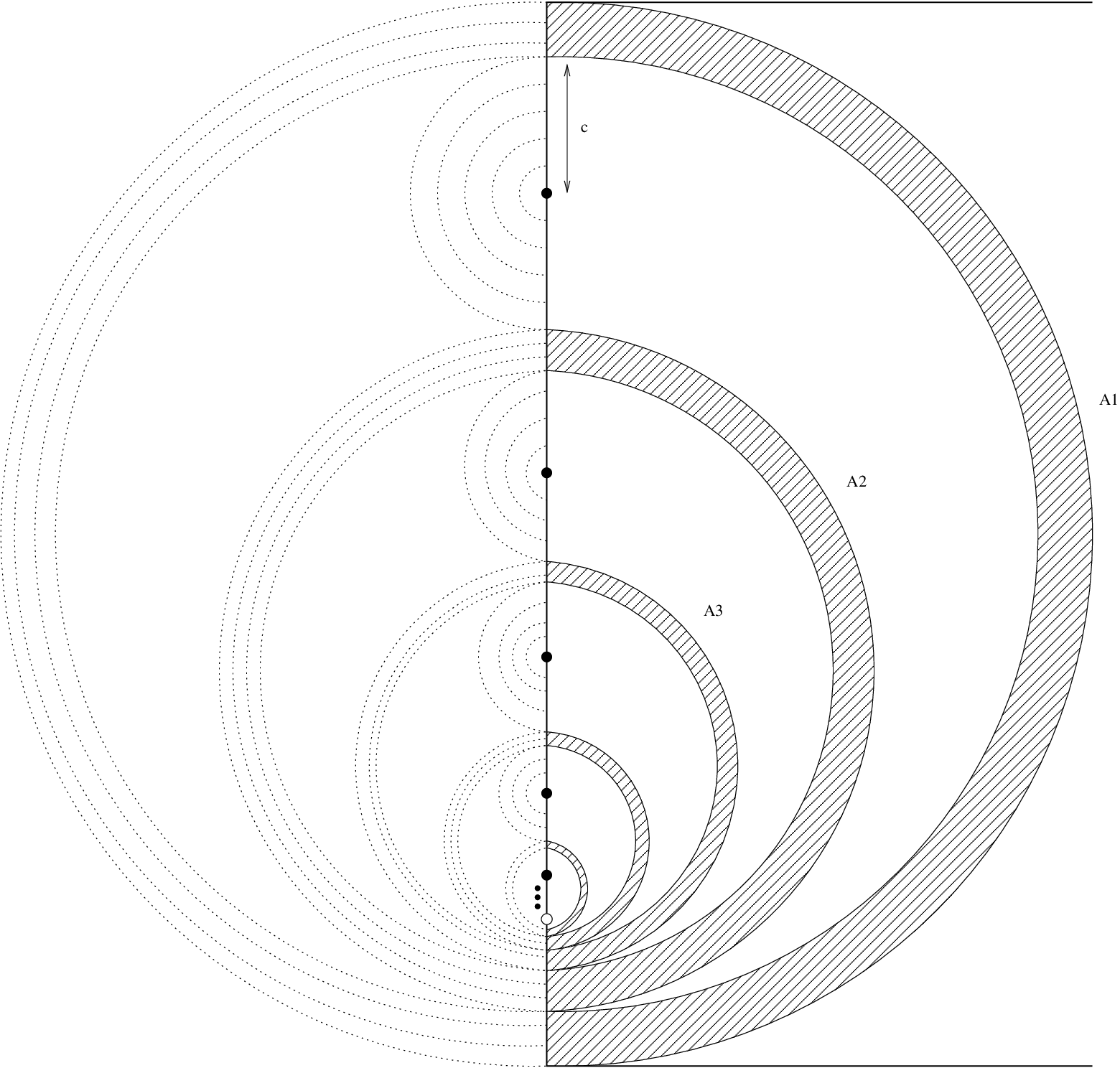}
\caption{Annuli with divergent moduli sum surrounding an 
accumulation of singularities of type $W$.}
\label{fig:Wpuncture}
\end{center}
\end{figure}

Let $\mu=\lambda^{-N}$, where~$N$ is the period of~$p$, let $w$ be the
width of the annulus~$A_0$, and denote by~$c$ the size of the largest
semi-circular arc band. Then $\Width(A_k)=w\mu^k$, and the
$k^{\rm th}$ semi-circular arc band has size $c\mu^k$. 

The outer radius~$R_k$ of~$A_k$ is one half of the distance between
the highest and lowest points of~$A_k$, in other words
\begin{eqnarray*}
R_k&=&\frac{1}{2}\left(2w\sum_{j=k}^\infty
\mu^j+2c\sum_{j=k}^\infty\mu^j\right)\\
&=&\frac{(w+c)\mu^k}{1-\mu}.
\end{eqnarray*}
Thus
\[\Area(A_k)=\frac{\pi}{2}\Width(A_k)(2R_k-\Width(A_k))=C\mu^{2k},\]
where $C$ depends only on $\mu$, $w$, and~$c$. Hence
\[\Mod(A_k)\ge\frac{\Width(A_k)^2}{\Area(A_k)}=\frac{w}{C}\]
so that $\sum_{k\ge k_0}\Mod(A_k)$ diverges as required.

Similar arguments apply to accumulations of singularities of
types~$V$. For a junction of type~$V_1$, let~$R$ and~$R'$ denote the
adjacent rectangles, where~$R$ is the rectangle adjacent to the switch
to which all of the bubbles in the junction are attached. Then a
sequence of half-annuli is constructed in~$R$ as for junctions of
type~$W$, but rather than being completed to annuli by the
identifications they are completed by corresponding half-annuli
in~$R'$. The moduli of these annuli are again bounded below by a
positive constant. Junctions of type~$V_2$ are double covers of
junctions of type~$W$, and a sequence of annuli with moduli bounded
below by a positive constant at an accumulation of type~$W$ thus
induces a similar sequence at an accumulation of type~$V_2$.
Accumulations of type~$V_3^\pm$ (respectively $V_1B$, $V_2B^2$) are
identical to those of type~$W^\pm$ (respectively~$V_1$,~$V_2$) in a
neighbourhood of the accumulation point.

Finally, consider the point at infinity. The pattern of identification
of horizontal sides, described in Section~\ref{sec:hident}, is an
$n$-fold cover of the pattern in junctions of type~$W$: there are~$n$
points on the horizontal boundary of $\cR$ which are identified to the
point at infinity, and the identifications along the horizontal sides
are given by an alternating sequence of rectangular and semi-circular
arc bands whose sizes decay exponentially.
\end{proof}

\def\cprime{$'$}

\end{document}